\documentclass[a4paper,times,3p]{elsarticle}
\usepackage[utf8]{inputenc}

\usepackage[T1]{fontenc}

\usepackage{microtype}

\usepackage{booktabs}

\usepackage{amsmath}
\usepackage{amssymb}
\usepackage{amsthm}
\usepackage{bm}
\usepackage{graphicx}
\usepackage{subfigure}
\usepackage{epstopdf}
\usepackage{psfrag}
\usepackage{color}
\usepackage{hyperref}

\usepackage{mathtools}

\definecolor{cmykcyan}{cmyk}{1,0,0,0}
\definecolor{cmykred}{cmyk}{0,1,1,0}
\definecolor{cmykblack}{cmyk}{0,0,0,1}

\newtheorem*{remark}{Remark}

\usepackage[font=footnotesize,labelfont=bf]{caption}
\usepackage{multirow}

\usepackage{pxfonts}

\usepackage{import}

\graphicspath{./graphics/}

\usepackage[textsize=scriptsize]{todonotes}
\usepackage{setspace}

\journal{Computer Methods in Applied Mechanics and Engineering}

\begin{document}

\begin{frontmatter}

\title{Boundary-Conforming Finite Element Methods for Twin-Screw Extruders using Spline-Based Parameterization Techniques}

\author[add1]{Jochen Hinz\corref{cor1}}
\ead{j.p.hinz@tudelft.nl}
\author[add2]{Jan Helmig\corref{cor1}}
\ead{helmig@cats.rwth-aachen.de}
\author[add1]{Matthias M\"oller\corref{}}
\ead{m.moller@tudelft.nl}
\author[add2]{Stefanie Elgeti\corref{}}
\ead{elgeti@cats.rwth-aachen.de}

\cortext[cor1]{Corresponding author}

\address[add1]{Delft Institute of Applied Mathematics\\
	Delft University of Technology, 2628 XE Delft, the Netherlands}
\address[add2]{Chair for Computational Analysis of Technical Systems (CATS)\\
	RWTH Aachen University, 52056 Aachen, Germany}

\begin{abstract}
This paper presents a novel spline-based meshing technique that allows for usage of boundary-conforming meshes for unsteady flow and temperature simulations in co-rotating twin-screw extruders. Spline-based descriptions of arbitrary screw geometries are generated using Elliptic Grid Generation. They are evaluated in a number of discrete points to yield a coarse classical mesh. The use of a special control mapping allows to fine-tune properties of the coarse mesh like orthogonality at the boundaries. The coarse mesh is used as a \textit{`scaffolding'} to generate a boundary-conforming mesh out of a fine background mesh at run-time. Storing only a coarse mesh makes the method cheap in terms of memory storage. Additionally, the adaptation at run-time is extremely cheap compared to computing the flow solution. Furthermore, this method circumvents the need for expensive re-meshing and projections of solutions making it efficient and accurate. It is incorporated into a space-time finite element framework. We present time-dependent test cases of non-Newtonian fluids in 2D and 3D for complex screw designs. They demonstrate the potential of the method also for arbitrarily complex industrial applications.
\end{abstract}

\begin{keyword}
Co-Rotating Twin-Screw Extruder\sep Elliptic Grid Generation \sep Isogeometric Analysis \sep Space-Time Finite Elements \sep Mesh Update Method\sep Boundary-Conforming Finite Elements \sep Non-Isothermal Flow
\end{keyword}
\end{frontmatter}

\section{Introduction}

Co-rotating twin-screw extruders are widely-used devices in the plastic-producing industry. They are used to distribute and disperse polymers and additives since they provide short residence times and extensive mixing. Typical screw elements are conveying elements that are used to transport the plastic melt from the feed section towards the die as well as kneading and mixing elements. The latter are especially tailored to account for dispersive and distributive mixing. \\
Performing experiments to obtain detailed information about velocity, pressure, and temperature distribution is very complex, time-consuming, and expensive, if not impossible for twin-screw extruders.
This makes Computational Fluid Dynamics (CFD) using finite element analysis an appealing tool to obtain detailed information about the flow inside twin-screw extruders. However, numerical analysis of the flow is not trivial due to constantly moving domains and small gap sizes which makes meshing extremely difficult. To tackle this problem, a wide range of approaches have been proposed.\\

1D models have been developed in order to predict the pressure build-up \citep{chen1991dimensionless} or fill length and specific mechanical energy \cite{vergnes1998global}. Full 3D results have been obtained using commercial software like POLYFLOW which is based on a mesh superposition method that avoids re-meshing \citep{zhang2009numerical,alsteens2004parametric}.
A rather recent development are methods using smoothed particle hydrodynamics \citep{robinson2018effect,wittek2018accuracy}. These methods are extremely useful in case only partly filled extruders are considered. \\
Fictitious domain methods are also quite popular. Meshing is avoided since the actual geometry is embedded into a fixed background mesh. A classical example is the method presented in \citep{valette2009direct} that is used in \citep{durin2014comparison} to compare 3D results with 1D estimates.
However, in \citep{sarhangi2012adaptive} it has been shown that classic fictitious domain methods lack accuracy inside the small gap regions. This drawback can be mitigated by, for example, XFEM. Other examples of fictitious domain methods are the Body Conformal Enrichment method presented in \citep{hetu2013immersed} or a method that tries to concentrate the background mesh at the screw interface using algebraic operations \citep{ianus2014mesh} . 
As a drawback, such methods require additional effort to describe the screw boundary exactly as well as to capture the flow effects in the small gaps. Furthermore, load balancing for highly parallel large scale computations is extremely difficult. \\
A different approach are boundary-conforming methods. The mesh represents the geometry exactly in terms of the underlying finite element discretization, which allows to strongly impose boundary conditions as well as the construction of high-quality boundary-layer meshes.
The drawback is that the generation of a time-dependent boundary conforming mesh can be difficult.
Thus, most of the simulations are based on a ''snapshot'' technique. An individual mesh is generated for each screw orientation. This is a valid approach for flow computations since the flow inside the extruder can be considered to be quasi steady. This approach has been applied in many publications such as \citep{bravo2000numerical,malik2014simulation}.
Also steady temperature results have been obtained in \citep{ishikawa20013,kalyon2007integrated}. However, these methods are limited when simulating time-dependent quantities.
A continuous time-dependent mesh is needed. This requires efficient mesh update techniques to avoid constant re-meshing and in turn expensive projections onto to the new mesh. \\
In \citep{helmig2018boundary}, we present the Snapping Reference Mesh Update Method (SRMUM), an efficient mesh update method for twin-screw extruders. It allows the computation of time-dependent temperature results without any need for re-meshing.
The method uses a structured background mesh that constantly adapts to the current screw configuration in a boundary-conforming manner. As a drawback, SRMUM is limited to convex screw geometries. \\
Instead, \citep{hinz2018elliptic} intends to generate analysis-suitable spline-based parameterizations. The resulting meshes are suitable for numerical analysis using Isogemetric Analayis (IGA) as introduced in \citep{hughes2005isogeometric}. The approach is based on Elliptic Grid Generation, which can handle nonconvex screw geometries.
The method is PDE-based, which makes it suitable for gradient-based shape optimization by simply adding the parameterization PDE-problem as a constraint to the shape optimization formulation, which can then be differentiated with respect to the design variables \citep{hinz2018spline}. Furthermore, classical finite-element meshes have been extracted from the spline parameterization and used to compute first flow results inside idealized twin-screw compressors \citep{moller2018isogeometric}. \\
In this work, we aim to adapt the general concept of spline-based meshing developed for twin-screw compressors to twin-screw extruders and combine it with certain aspects of SRMUM. This allows us to generate flexible, high quality, time-dependent and FEM-suitable meshes. \\
\begin{figure}[t!]
    \centering
    \includegraphics[width = 0.9 \textwidth]{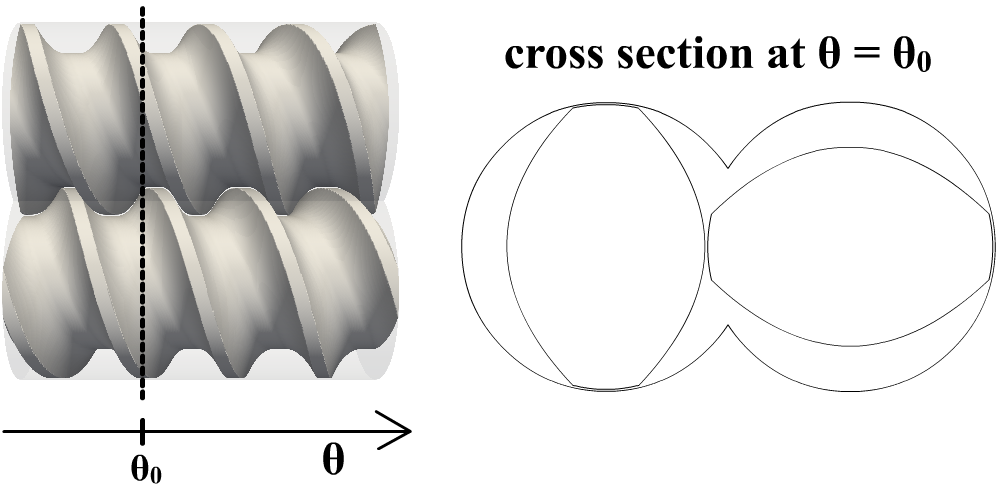}
    \caption{A volumetric mesh for an extruder geometry is built from a large number of planar cross sections parameterized in the plane using Elliptic Grid Generation with spline functions.}
    \label{fig:cross_section}
\end{figure}
\noindent The approach is based on taking a large number of cross sections through a volumetric extruder geometry leading to $2$D extruder contours at a large number of rotational angles $\theta_i$ (see Figure \ref{fig:cross_section}). A spline parameterization $\mathbf{x}_i$ is generated for the $2$D profiles at every $\theta_i$ using Elliptic Grid Generation techniques. \\ The spline parameterizations are evaluated in a large number of points to generate a mesh \textit{`scaffolding'} to adapt a SRMUM-like background mesh to the current screw shape via interpolation. A volumetric mesh is then built from the planar meshes by stacking them in the $z$-direction. \\
Once a spline parameterization has been generated, it is relatively cheap to generate finer meshes, e.g. for refinement studies. Furthermore, function evaluations in the spline mappings $\mathbf{x}_i$ can be tuned using a \textit{'control mapping'} in order to account for special properties of the mesh, like orthogonality at the screw surface. Since only a \textit{`scaffolding'} is extracted, it is not necessary to save the whole $3$D mesh for all possible time instances which would be very time- and memory-consuming. Furthermore, the computational cost of the mesh update at each time step during the simulation is only a fraction of the time spent for solving the flow solution. \\
The paper is structured as follows: Section \ref{sec:gridgeneration} explains the meshing concept. This includes the generation of a spline parameterization up to the analysis-suitable mesh. The underlying physics and equations required to compute the flow inside twin-screw extruders are stated in Section \ref{sec:equations}, as well as the space-time finite element solution method. Section \ref{sec:numericalexamples} contains $2$D and $3$D examples of flows through twin-screw extruders that aim to validate as well as highlight the advantages of the presented method.
\section{Grid Generation} \label{sec:gridgeneration}
Since the geometry pipeline provides no more than a description of the twin-screw profile (typically in terms of a point cloud), the first step towards numerical simulation is the generation of a planar computational mesh. For the mesh generation, we employ the following two-step approach: at every discrete angle $\theta$, we first generate a spline-based geometry description which we then evaluate it in a large number of discrete points with known connectivity in order to yield a classical mesh. The advantages of a spline-based description is the possibility to fine-tune the properties of the mesh by appropriately choosing the evaluation abscissae, which follow from a precomputed control mapping. \\
In the following, we give a brief introduction to spline-surfaces after which we discuss the basics of the numerical algorithm employed in the meshing process and its applications to the target geometry. 

\subsection{B-Splines}
\label{sect:B_Splines}
\begin{figure}[t!]
    \centering
    \includegraphics[scale = 0.45]{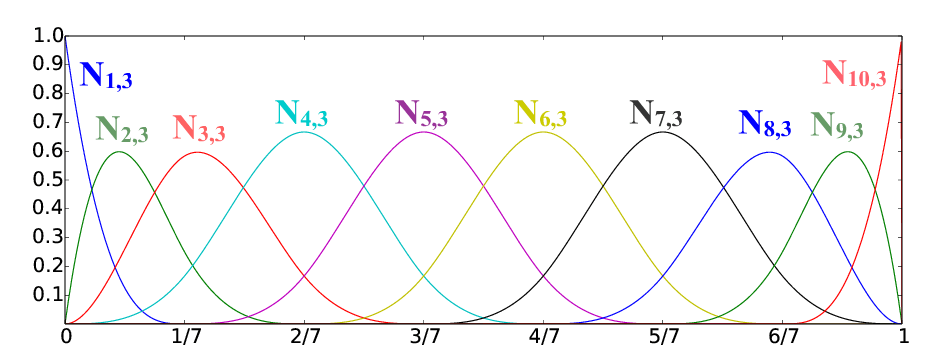}
    \caption{The univariate B-spline basis resulting from the knot-vector $\Xi_3$.}
    \label{fig:open_kv_p_3}
\end{figure}

As described in \citep{piegl2012nurbs}, B-splines are piecewise-polynomial functions that can be constructed so as to satisfy various continuity properties at the places where the polynomial segments connect. Their properties are determined by the entries of the so-called \textit{knot vector}
\begin{align}
    \Xi = \{ \xi_1, \xi_2, \ldots, \xi_{n + p + 1} \}.
\end{align}
The knot vector is a nondecreasing sequence of parametric values $\xi_i \subset [0,1]$ that determine the boundaries of the segments on which the spline-basis is polynomial. Selecting some polynomial order $p$, the $p$-th order spline-functions $N_{i,p}$ are constructed recursively, utilizing the relation (with $\tfrac{0}{0} \equiv 0$)
\begin{align}
    N_{i,s}(\xi) = \frac{ \xi - \xi_i }{ \xi_{i+1} - \xi_i } N_{i, s-1}(\xi) + \frac{ \xi_{i+s+1} - \xi }{ \xi_{i+s+1} - \xi_{i+1} } N_{i+1, s-1}(\xi),
\end{align}
starting from
\begin{align}
    N_{i,0} = \left \{ \begin{array}{ll} 1 & \text{if } \xi_i \leq \xi \leq \xi_{i+1} \\
                                         0 & \text{otherwise} \end{array} \right. ,
\end{align}
and iterating until $s = p$. The support of basis function $N_{i,p}$ is given by the interval $\mathcal{I}_{i,p} = [\xi_i, \xi_{i+p+1}]$ and the amount of continuous derivatives across knot $\xi_j$ is given by $p - m_j$, where $m_j$ is the multiplicity of $\xi_j$ in $\mathcal{I}_{i,p}$. In practice, $\xi_1 = 0$ is repeated $p+1$ times as well as $\xi_{n+p+1}$ such that $\xi_1 = \ldots = \xi_{p+1} = 0$ and $\xi_{n+1} = \ldots = \xi_{n+p+1} = 1$. As a result, the corresponding basis $\sigma = \{ N_{1,p}, \ldots, N_{n,p} \}$ forms a non-negative partition of unity on the entire parametric domain $[0,1]$, that is:
\begin{align}
    \sum_{i = 1}^n N_{i,p}(\xi) = 1,
\end{align}
with
\begin{align}
    N_{i,p}(\xi) \geq 0,
\end{align}
for all spline functions $N_{i,p}$ \cite{hughes2005isogeometric}. Figure \ref{fig:open_kv_p_3} shows the $p = 3$ B-spline basis resulting from the knot-vector
\begin{align}
    \Xi = \left \{0, 0, 0, 0, \tfrac{1}{7}, \tfrac{2}{7}, \tfrac{3}{7}, \tfrac{4}{7}, \tfrac{5}{7}, \tfrac{6}{7}, 1, 1, 1, 1 \right \}.
\end{align}

\noindent The extension to bivariate spline bases is now straight-forward: given two univariate bases $\sigma_\Xi = \{N_{1, p}, \ldots, N_{n, p}\}$ and $\sigma_{\mathcal{H}} = \{M_{1, q}, \ldots M_{m, q} \}$, we build a bivariate basis $\Sigma = \{ w_{i,j} \}_{(i,j) \in \{1, \ldots, n \} \times \{ 1, \ldots, m \} }$, by means of a tensor-product, where
\begin{align}
    w_{i, j}(\xi, \eta) = N_{i, p}(\xi) M_{j, q}(\eta).
\end{align}
The values contained in $\Xi$ and $\mathcal{H}$ (without knot-repetitions) hereby become the boundaries of the polynomial segments, which can be regarded as the counterparts of classical elements. \\
We construct the mapping of a B-spline surface as follows:
\begin{align}
\label{eq:mapping_tensorial}
    \mathbf{x} = \sum_i \sum_j \mathbf{c}_{i,j} w_{i,j},
\end{align}
where the $\mathbf{c}_{i,j} \in \mathbb{R}^2$ are referred to as the \textit{control points}. We refer to the $\mathbf{c}_{i,j}$ with $i \in \{1,n \}$ or $j \in \{1, m\}$ as the \textit{boundary control points} while the remaining $\mathbf{c}_{i,j}$ are called \textit{inner control points}. Bluntly put, the local density of spline basis functions is determined by the local number of knots. As such, a cleverly chosen knot-vector can be utilized to properly resolve important features of the geometry. Furthermore, as $\mathbf{x}$ is a linear combination of the $w_{i, j} \in \Sigma$, it will inherit the local continuity properties of the basis. Therefore, many geometrical features can be better captured by a clever choice of the knot multiplicities in $\Xi$ and $\mathcal{H}$.

\subsection{Spline-Based Meshing Techniques}
In the following, we shall drop the tensor-product indexing used in (\ref{eq:mapping_tensorial}) and replace it by a single global index. As such, the mapping operator will be of the form
\begin{align}
\label{eq:mapping}
    \mathbf{x} = \sum_{i \in \mathcal{I}_\text{inner} } \mathbf{c}_i w_i + \sum_{j \in \mathcal{I}_\text{boundary} } \mathbf{c}_j w_j,
\end{align}
where $\mathcal{I}_\text{inner}$ and $\mathcal{I}_\text{boundary}$ refer to the index-set of inner and boundary basis functions, respectively. The mapping $\mathbf{x}$ is a function from the parameter domain $\hat{\Omega} = [0, 1]^2$ to the physical geometry $\Omega$, whose boundaries $\partial \Omega$ follow from a collocation of the employed spline-basis $\Sigma = \{w_1, \ldots, w_N\}$, resulting from the tensor product knot vector $\boldsymbol{\Xi}$, to the input point cloud. \\
For the proper selection of $\boldsymbol{\Xi}$ and the subsequent collocation, we utilize the approach illustrated in \cite{hinz2018spline}. Our methodology consists of a stabilized least-squares fit of the points against the basis and an adaptive reselection of $\boldsymbol{\Xi}$ based on the local magnitude of the projection residual. The process is repeated until the residual is deemed sufficiently small. Upon completion, the restriction of $\mathbf{x}$ to $\partial \hat{\Omega}$ parameterizes $\partial \Omega$. Note that the above procedure essentially selects the boundary control points in (\ref{eq:mapping}) such that $\mathbf{x}\vert_{\partial \hat{\Omega}}$ is a proper collocation of the input point cloud (hence the word \textit{boundary control points}). However, the inner control points are unknown at this stage. As such, the task of any parameterization algorithm is to properly select the $\mathbf{c}_i$ in (\ref{eq:mapping}), such that the resulting parameterization is folding-free. \\
A computationally inexpensive, yet often sufficiently powerful method is transfinite interpolation \cite{gordon1973transfinite}. The inner control points follow from an interpolation of the southern and northern boundaries in combination with an interpolation from east to west. Let $\boldsymbol{\gamma}_e, \boldsymbol{\gamma}_w, \boldsymbol{\gamma}_s$ and $\boldsymbol{\gamma}_n$ parameterize the four segments of $\partial \Omega$. In a spline-based setting, the mapping operator is constructed as follows:
\begin{align}
\label{eq:transfinite_interpolation}
    \mathbf{x}(\xi, \eta) &= (1 - \xi) \boldsymbol{\gamma}_w + \xi \boldsymbol{\gamma}_e + (1 - \eta) \boldsymbol{\gamma}_s + \eta \boldsymbol{\gamma}_n - (1 - \xi)(1 - \eta) \mathbf{p}_{0, 0} - \xi \eta \mathbf{p}_{1, 1} - \xi (1 - \eta) \mathbf{p}_{1, 0} - (1 - \xi) \eta \mathbf{p}_{0, 1},
\end{align}
where the $\mathbf{p}_{i, j}$ denote the corners at $\mathbf{x}(i, j)$. The symbolic parameterization from (\ref{eq:transfinite_interpolation}) constitutes a recipe for determining the $\mathbf{c}_i$ in (\ref{eq:mapping}). They can be computationally inexpensively determined by a $L_2$-projection or a collocation over the \textit{Greville-abscissae} \cite{johnson2005higher} corresponding to $\boldsymbol{\Xi}$. Whenever $\mathbf{x}$ is an O-type parameterization, we simply remove the nonexistent pair of boundary contours ($\boldsymbol{\gamma}_n, \boldsymbol{\gamma}_s$ or $\boldsymbol{\gamma}_w, \boldsymbol{\gamma}_e$) from (\ref{eq:transfinite_interpolation}), as well as the $\mathbf{p}_{i, j}$ and perform a unidirectional interpolation. \\
Transfinite interpolation is an $\mathcal{O}(N)$ operation but does not guarantee a folding-free mapping. In order to handle the complex characteristics of twin-screw extruders, it is desirable to have a more powerful parameterization technique in our arsenal. \\
A second class of parameterization methods are PDE-based, notably approaches based on the principles of \textit{Elliptic Grid Generation} (EGG). The main purpose of EGG is to compute a mapping $\mathbf{x}: \hat{\Omega} \rightarrow \Omega$ such that the components of $\mathbf{x}^{-1}: \Omega \rightarrow \hat{\Omega}$ are harmonic. This is accomplished by imposing the following system of PDEs on $\mathbf{x}$:
\begin{align}
\label{eq:EGG_not_discretized}
    g_{22} \mathbf{x}_{\xi \xi} - 2 g_{12} \mathbf{x}_{\xi \eta} + g_{11} \mathbf{x}_{\eta \eta} &= 0, \quad \text{s.t. } \mathbf{x} \vert_{\partial \hat{\Omega}} = \partial \Omega,
\end{align}
where the $g_{ij}$ denote the entries of the metric tensor corresponding to $\mathbf{x}$. We have
\begin{align}
    \begin{pmatrix} g_{11} & g_{12} \\ g_{21} & g_{22} \end{pmatrix} & = \begin{pmatrix} \mathbf{x}_\xi \cdot \mathbf{x}_\xi & \mathbf{x}_\xi \cdot \mathbf{x}_\eta \\ \mathbf{x}_\eta \cdot \mathbf{x}_\xi & \mathbf{x}_\eta \cdot \mathbf{x}_\eta \end{pmatrix}.
\end{align}
A justification of this approach is based on the observation that if $\partial \Omega$ satisfies certain regularity conditions, $\mathbf{x}$ will be bijective (and hence, folding-free) \cite{azarenok2009generation}. As such, any sufficiently accurate approximation of $\mathbf{x}$ will be bijective as well.

\subsection{Numerical Implementation}
A basic numerical algorithm to approximately solve (\ref{eq:EGG_not_discretized}) with globally $C^{\geq1}$-continuous spline basis functions is discussed in \cite{hinz2018elliptic}. Here, we follow a similar approach with auxiliary variables, first introduced in \cite{hinz2019iga}, that allows for $C^0$-continuities in $\xi$-direction (whose purpose shall become apparent in subsection \ref{subsect:Applications_Geometry}). We approximately solve (\ref{eq:EGG_not_discretized}) with a Galerkin approach, introducing the auxiliary variable $\mathbf{u}$ that serves the purpose of reducing the continuity requirements of the basis in $\xi$-direction from $C^1$ to $C^0$. \\
Let
\begin{align}
\label{eq:auxilliary_mapping}
    \mathbf{u}(\xi, \eta) &= \sum_{i=1}^{\tilde{N}} \mathbf{d}_i \tilde{w}_i(\xi, \eta).
\end{align}
where the $\tilde{w}_i \in \tilde{\Sigma}$ are basis functions taken from the auxiliary spline basis $\tilde{\Sigma}$. \\
In order to discretize (\ref{eq:EGG_not_discretized}), we express the $\xi$-derivatives in terms of $\mathbf{u}$ and introduce a scaling. Consider
\begin{align}
    U( \mathbf{u}, \mathbf{x} ) &= \frac{ g_{22} \mathbf{u}_\xi - g_{12} \mathbf{u}_\eta - g_{12} \mathbf{x}_{\xi \eta} + g_{11} \mathbf{x}_{\eta \eta} }{ \underbrace{g_{11} + g_{22}}_{\geq 0} + \epsilon },
\end{align}
with $\epsilon >0$ a small term that ensures numerical stability. In order to determine the $\mathbf{c}_i$ in (\ref{eq:mapping}) and $\mathbf{d}_i$ in (\ref{eq:auxilliary_mapping}), we solve the system
\begin{align}
\label{eq:EGG_discretized}
    \forall \boldsymbol{\sigma}_i \in [\tilde{\Sigma}]^2 \times [\Sigma_0]^2: \int_{\hat{\Omega}} \boldsymbol{\sigma}_i \cdot \begin{pmatrix} \mathbf{x}_\xi - \mathbf{u} \\ U(\mathbf{u}, \mathbf{x}) \end{pmatrix} \mathrm{d} \boldsymbol{\xi} = 0, \quad \text{s.t. } \mathbf{x} \vert_{\partial \hat{\Omega} } = \partial \Omega,
\end{align}
where $\Sigma_0 = \{ w_i \in \Sigma \enskip \vert \enskip i \in \mathcal{I}_\text{inner} \}$. For a memory-efficient algorithm to tackle the root-finding problem (\ref{eq:EGG_discretized}), we refer to the aforementioned publication \cite{hinz2019iga}. \\ 
Upon solving (\ref{eq:EGG_discretized}), we are in the possession of $\mathbf{x}$ and $\mathbf{u}$. Here, $\mathbf{x}$ parameterizes $\Omega$ and $\mathbf{u}$ serves no further purpose and can be discarded. Unlike the exact solution of (\ref{eq:EGG_not_discretized}), its numerical approximation may fold due to numerical inaccuracies. Since a defect is a direct cause of insufficient numerical accuracy, it can be repaired by increasing the local density of spline functions wherever the defect is located and recomputing a mapping operator from the enriched spline basis. In practice, we have observed folding to be a rare occurrence.

\subsection{Applications to the Geometry}
\label{subsect:Applications_Geometry}
\begin{figure}[t!]
    \centering
    \includegraphics[width = 0.7 \textwidth]{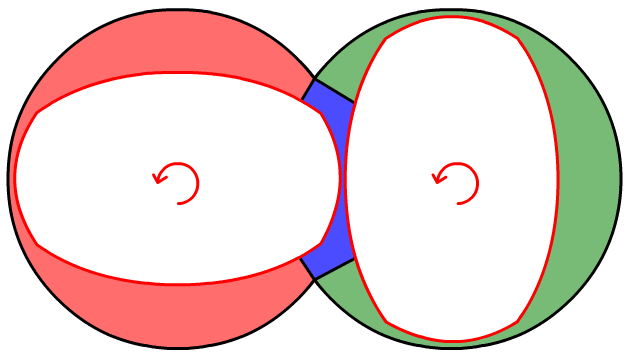}
    \caption{The topology we use for the parameterization of the geometry. Here, black boundaries are held fixed while red boundaries slide along the grid. The two $C$-grids (red, green) are parameterized (if possible) with transfinite interpolation while we use EGG (see (\ref{eq:EGG_discretized})) for the separator (blue). }
    \label{fig:topology}
\end{figure}
\begin{figure}[t!]
    \centering
    \includegraphics[width = 0.7 \textwidth]{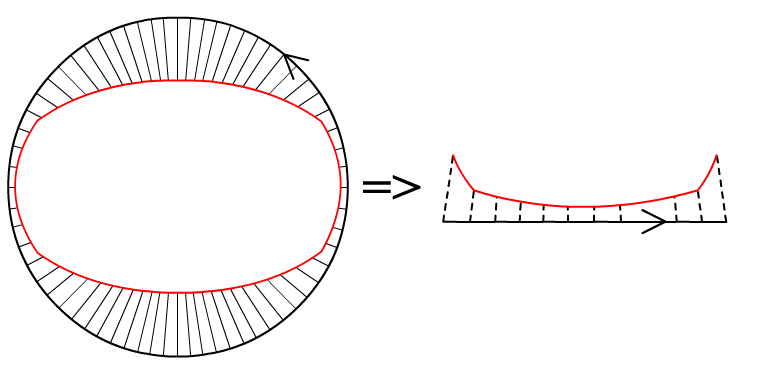}
    \caption{The first step towards a three-patch parameterization is finding a valid parameterization of a single rotor-casing $O$-grid. A valid parameterization is one in which vertical isolines do not cross upon a pullback of the circular casing onto a reference line segment.}
    \label{fig:valid_O_grid}
\end{figure}
\begin{figure}[t!]
    \centering
    \includegraphics[width = 0.7 \textwidth]{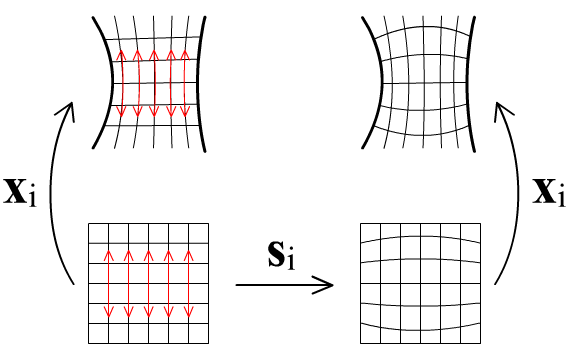}
    \caption{The control mapping serves the purpose of altering the parametric properties of the composite mapping. Here, we have chosen to orthogonalize the grid lines.}
    \label{fig:control_mapping}
\end{figure}
\begin{figure}[t!]
    \centering
    \includegraphics[width = 0.55 \textwidth]{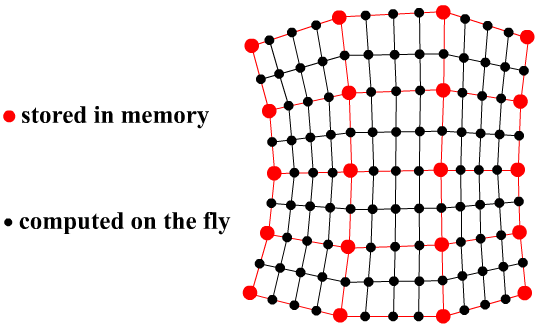}
    \caption{The composite mapping $\mathbf{x}_i \circ \mathbf{s}_i$ is evaluated in a number of uniformly-spaced abscissae and the resulting points (red) serve as a \textit{`scaffolding'} for the remaining mesh vertices. Hereby it is important to make a good trade off between memory requirements and grid quality (more evaluation points tends to increase grid quality)}
    \label{fig:evaluation}
\end{figure}

\noindent For the parameterization of the main geometry, we choose the topology illustrated in Figure \ref{fig:topology}. The topology consists of two $C$-type parameterizations connected by a third parameterization (blue), henceforth referred to as the \textit{'separator'}. The topology is designed in such a way that upon rotation, the rotor lobes slide along the grid in the parametric domain while the grid is held fixed at the casings. Furthermore, the casings are held arc-length parameterized while the parametric properties of the rotor lobes depend on the rotation angle.\\
As such, the first step in building a three-patch parameterization is finding a valid $O$-grid parameterization of a single rotor plus casing. We take the circular casing and map it onto a straight line segment. A valid parameterization is one in which the vertical isolines do not intersect under the mapping (see Figure \ref{fig:valid_O_grid}).
To find a valid lobe parameterization, we employ the matching algorithm from \cite{hinz2018spline} to the rotor and casing point clouds. The points are matched based on their Euclidean distance in a hierarchical fashion and matched points receive the same parametric value. Upon completion, we have a monotone reparameterization function $\xi^\prime: [0, 1] \rightarrow [0, 1]$ and we assign the parametric value $\xi_i = \xi^\prime( \hat{\xi}_i )$ to the $i$-th point in the rotor point cloud. Here, the $\hat{\xi}_i$ result from a chord-length parameterization of the point cloud. The steps described above are carried out for $\theta = 0$ only. The reparameterization function can be reused for different rotational angles by performing an appropriate periodic shift that is based on the value of $\theta$. \\
The coordinates of the points along with their parametric values are utilized for a collocation of the left and right rotor $O$-grids to a suitable spline basis. Here, the knot vector is chosen such that the collocation residual is below a user-defined threshold (see \cite{hinz2018spline}). Upon completion, the boundary curves are known and we utilize interpolation (see (\ref{eq:transfinite_interpolation})) to form folding-free $O$-grid parameterizations for both rotor casing pairs.
Next, we cut the left and right rotor $O$-grids at the parametric values that correspond to the casing CUSP points to form the left and right rotor $C$-grids. \\
We repeat above steps for a large number of discrete angles $\Theta = \{\theta_1, \ldots, \theta_n\}$ taken from the interval that corresponds to a rotation after which the initial position is again obtained (typically the interval $[0, \pi]$). Upon completion, we have a database that is filled with $C$-grid parameterizations for all discrete angles $\Theta$. We combine the $C$-grid boundaries along with the rotor lobe parts that were cut from the rotor-casing $O$-grids in order to generate one boundary description of the separator, as shown in Figure \ref{fig:topology}. Clearly, a spline collocation of the northern and southern boundaries has to be performed with a locally $C^0$-continuous spline basis, due to the spiked nature of the input point clouds. Therefore, we assign the parametric value $\xi = 0.5$ to the CUSP-point and utilize a knot vector with $p$-fold internal knot repetition at $\xi=0.5$ for the $\xi$-direction. \\
The next step is to assign suitable parametric values to the points contained in the eastern and western input point clouds. Unfortunately, we have observed a chord-length parameterization to lead to unsatisfactory results. In order to ensure that similar parametric values are assumed on either sides of small gaps in the geometry, we again apply the matching algorithm from \cite{hinz2018spline}. In contrast to the $C$-grid case, we do not keep one side chord-length parameterized but let both sides float. If two points $p_i^l$ and $p_j^r$ have been matched, we assign the parametric value 
\begin{align*}
    \eta = \tfrac{1}{2}( \hat{\eta}_i^l + \hat{\eta}_j^r )
\end{align*}
to both of them, where the values of $\hat{\eta}_{i}^{l}$ and $\hat{\eta}_{j}^{r}$ correspond to chord-length parameterizations. Upon completion, we have two reparameterization functions that we utilize to assign parametric values to the input point clouds, similar to the $C$-grid case. \\
We perform the above steps for each $k$-th discrete angle $\theta_i$. Typically we take $k = 5$. A global reparameterization function is constructed by blending each $k$-th reparameterization function over the $\theta$-interval. Blending serves the purpose of achieving a degree of smoothness in the parametric properties of the separator as a function of $\theta$, which would be lost if we reparameterized at every $\theta_i$, due to the discrete nature of the reparameterization algorithm. \\
Upon completion, we start filling our database with separator parameterizations utilizing the EGG algorithm from (\ref{eq:EGG_discretized}). The auxiliary spline space $\bar{\Sigma}$ is of Raviart-Thomas type \cite{buffa10isogeometric} on each of the macro elements \cite{bressan2013isogeometric} $\boldsymbol{\xi} \in [0, 0.5] \times [0, 1]$ and $\boldsymbol{\xi} \in [0.5, 1] \times [0, 1]$ with $C^0$ interface coupling. To be precise, let 
\begin{align}
    \Sigma = \sigma_\Xi^{p_1} \times \sigma_\mathcal{H}^{p_2} \quad \text{and} \quad \bar{\Sigma} = \sigma_{\bar{\Xi}}^{q_1} \times \sigma_{\bar{\mathcal{H}}}^{q_2}
\end{align} be built from the knot vectors 
\begin{align}
    \boldsymbol{\Xi}^{p_1, p_2} = \Xi^{p_1} \times \mathcal{H}^{p_2} \quad \text{and} \quad \boldsymbol{\bar{\Xi}}^{q_1, q_2} = \bar{\Xi}^{q_1} \times \bar{\mathcal{H}}^{q_2},
\end{align} 
respectively. Here, the $p_i$ and $q_i$ denote the polynomial orders used in each direction. If
\begin{align}
    \Xi^{p_1} = \{\underbrace{0, \ldots, 0}_{p_1 + 1 \text{ times}}, \xi_1, \ldots, \xi_{q}, \underbrace{0.5, \ldots, 0.5}_{p_1 \text{ times}}, \xi_{r}, \ldots, \xi_s, \underbrace{1, \ldots, 1}_{p_1 + 1 \text{ times}} \},
\end{align}
we take $\boldsymbol{\bar{\Xi}}^{q_1, q_2} = \bar{\Xi}^{p_1 + 1} \times \mathcal{H}^{p_2}$, with
\begin{align}
    \bar{\Xi}^{p_1 + 1} = \{\underbrace{0, \ldots, 0}_{p_1 + 2 \text{ times}}, \xi_1, \ldots, \xi_{q}, \underbrace{0.5, \ldots, 0.5}_{p_1 + 1 \text{ times}}, \xi_{r}, \ldots, \xi_s, \underbrace{1, \ldots, 1}_{p_1 + 2 \text{ times}} \}.
\end{align}
In our examples, $\Sigma$ is bicubic, i.e., $p_1 = p_2 = 3$. \\
The root-finding problem is initialized with a transfinite mapping (see (\ref{eq:transfinite_interpolation})). Convergence can be further accelerated using multigrid techniques (see \cite{hinz2018elliptic}). We fill our database in a hierarchical fashion. As soon as the database contains separator parameterizations for a sufficiently large subset of angles from $\Theta$, we utilize interpolation in order to generate initial guesses for the remaining angles to further speed up the process. \\
After the algorithm has finished, we have a three-patch parameterization for the geometry at every $\theta_i$. A classical mesh can be generated by evaluating the spline mapping in a large number of uniformly-spaced points. The properties of the mesh can be further tuned by using a control mapping. We have observed this to be especially important for the mesh quality in the separator. By the tuple $(\mathbf{x}_i, \mathbf{s}_i)$, we denote the mapping operator and the corresponding control mapping for the separator at $\theta = \theta_i$. The control mapping $\mathbf{s}_i(\mu, \nu) = \sum_i \tilde{\mathbf{c}}_i \tilde{w}_i(\mu, \nu)$ is a function from the unit square onto itself. We choose its control points by solving the following minimization problem:
\begin{align}
\label{eq:minimization_control_mapping}
& \frac{1}{2} \int \limits_{ [0,1]^2 } \left( \frac{ \partial (\mathbf{x}_i \circ \mathbf{s}_i) }{ \partial \mu } \cdot \frac{ \partial (\mathbf{x}_i \circ \mathbf{s}_i) }{ \partial \nu } \right)^2 \mathrm{d} \mu \mathrm{d} \nu \rightarrow \min \limits_{ \{ \tilde{ \mathbf{c}_i } \}_i }, \quad \text{subject to boundary conditions.}
\end{align}
The cost function from (\ref{eq:minimization_control_mapping}) maximizes orthogonality in the composite mapping $\mathbf{x}_i \circ \mathbf{s}_i$. The integral is approximated with a discrete Riemann sum. The boundary conditions follow from the requirement that $\mathbf{s}_i$ be a mapping from the unit square onto itself. We only let $\mathbf{s}_i$ slide in $\eta$-direction (see Figure \ref{fig:control_mapping}) as this leads to better results for tube-like shaped geometries such as the separator. A sufficient linear condition for the aforementioned requirement is easily formulated and imposed as a constraint on (\ref{eq:minimization_control_mapping}). The minimization problem is tackled with an SLSQP \cite{kraft1988software} algorithm while the gradient of the cost function is approximated with finite differences. To speed up convergence, control mappings of previous slices can be utilized to initialize the next control mapping. We compute a control mapping for the separator at all discrete angles $\theta_i \in \Theta$. We do not compute control mappings for the $C$-grids (i.e., the control mapping is the identity).
\noindent Upon completion, we build a mesh \textit{'scaffolding'} for the separator from the composite mapping $\mathbf{x}_i \circ \mathbf{s}_i$ by evaluating it in a large number of uniformly-spaced points. By $n_{\mu, \nu}$, we denote the number of elements in $\mu$ and $\nu$ direction, respectively. An adaptive background mesh is generated by adding additional vertices via linear interpolation within the mesh scaffolding point cloud. The vertex positions are computed on the fly and do therefore not have to be stored explicitly (see Figure \ref{fig:evaluation}). For the background mesh, we denote the number of elements in radial and screw direction by $n_r$ and $n_s$ with $n_{r} \geq n_{\mu}$ and $n_{s} \geq n_\nu$, respectively. The ratios $n_{r}/ n_{\mu} \geq 1$ and $n_s / n_{\nu} \geq 1$ are hence based on a trade off between memory requirements and grid quality, where ratios closer to $1$ typically lead to better meshes (for a given number of elements) while increasing the memory requirements. The converse holds for larger ratios, whereby it is important to take $n_{\mu, \nu}$ sufficiently large so that the resulting mesh \textit{'scaffolding'} is sufficiently close to the folding-free $\mathbf{x}_i \circ \mathbf{s}_i$ and therefore does not fold itself. \\
\noindent For the C-grids, the mesh scaffolding cannot fold due to an insufficient number of elements in $\mu$-direction, thanks to the precomputed reparameterization function (see Figure \ref{fig:valid_O_grid}). Hence, we take $n_{\mu} = 1$ and add vertices until $n_r(C) = n_r( \text{separator}) / 2$ to yield a boundary-conforming mesh on both sides of the CUSP points. The storage requirements for the C-grids are thus fully determined by $n_\nu$, whose value is chosen to properly resolve the rotor boundaries. Typically, for the C-grids, we simply take $n_s = n_\nu$.\\
In order to construct a 3D mesh, we linearly connect the 2D meshes, where $n_a$ determines the number of elements in $z$-direction. Every slice in $z$-direction knows its initial angle, which is used to interpolate between the stacked planar slices of the \textit{`scaffolding'}.
\section{Governing Equations and Solution Method} \label{sec:equations}

\subsection{Governing Equations for Flow and Temperature of Plastic Melt}

 We choose to represent the molten polymer in the extruder as viscous, incompressible and temperature-dependent fluid. A time-dependent computational domain, denoted by $\Omega _t \; \subset \; \mathbb{R}^{n_{sd}}$, is considered. It is enclosed by its boundary $\Gamma _t$, where $t \in (0,T)$ is an instance of time and $n_{sd}$ the number of space dimensions. The velocity ${\bf u}({\bf x},t)$ and pressure $p({\bf x},t)$ are governed by the incompressible Navier-Stokes equations:

\begin{align}
\boldsymbol{\nabla \cdot} {\bf{u}} = 0 \quad \mbox{on} \ \Omega_t, \quad \forall t \in (0,T), \label{eq:cont} \\
 \rho \left( \dfrac{\partial \bf{u}}{\partial t} + \bf{u} \cdot \boldsymbol{\nabla} \bf{u} \right) - \boldsymbol{\nabla \cdot \sigma} = \boldsymbol{0} \quad
 \mbox{on} \ \Omega_t, \quad \forall t \in (0,T), \label{eq:momentum}
\end{align}

where $\rho$ is the fluid density. The stress tensor $\boldsymbol{\sigma}$ is used to close the set of equations:

\begin{align}
    \boldsymbol{\sigma} ( {\boldsymbol{u} }, p) = -p {\boldsymbol{I} } + 2 \eta \left( \dot{\gamma }, T  \right) \boldsymbol{\varepsilon}( {\boldsymbol{u} }), \\
    \boldsymbol{\varepsilon}( {\boldsymbol{u} }) = \frac{1}{2} \left( \boldsymbol{\nabla u} + \left( \boldsymbol{\nabla u} ^T \right) \right),
\end{align}

 with $\eta$ being the dynamic viscosity. It is constant for a Newtonian fluid and varies for Generalized Newtonian models with respect to temperature $T$ and shear rate $\dot{\gamma}$. \\

Using single-phase Navier-Stokes equations implies the assumption that the extruder is fully filled.
However, in practice it might occur that the flow channels inside the twin-screw extruder are only partially filled. In order to account for this, one would have to extend the model to also include the air phase, which results in solving a multi-phase flow problem. Popular numerical methods for multi-phase flow are the levelset \cite{osher1988} or the volume of fluid method \cite{hirt1981a}. For simplicity, only single-phase flow will be considered within this work. \\

The temperature $T({\bf x},t)$ inside the extruder is governed by the heat equation:

\begin{align}
\rho c_p \left( \dfrac{\partial T}{ \partial t} + \boldsymbol{u} \cdot \boldsymbol{\nabla} T \right) - \kappa \boldsymbol {\Delta} T - 2 \eta \boldsymbol{\nabla u} \colon \boldsymbol{\varepsilon} \left( \boldsymbol{u} \right)  \; = 0 \quad
 \mbox{on} \ \Omega_t, \quad \forall t \in (0,T). \label{eq:heat}
\end{align}

The Dirichlet and Neumann boundary conditions for temperature and flow are defined as:

\begin{align}
    {\bf u} = {\bf g}^f \; \mbox{on} \; \left( \Gamma_t \right) ^f _g, \\
    {\bf n} \cdot \boldsymbol{\sigma} = {\bf h}^f \; \mbox{on} \; \left( \Gamma_t \right) ^f _h, \\
    T = g^T \; \mbox{on} \; \left( \Gamma_t \right) ^T _g, \\
    {\bf n} \cdot \kappa \boldsymbol{\nabla} T = h^T \; \mbox{on} \; \left( \Gamma_t \right) ^T _h.
\end{align}

$ \left( \Gamma_t \right) ^i _g$ and $\left( \Gamma_t \right) ^i _h$ are complementary portions of $\Gamma _t ^i$, with $i=f\mbox{ (Fluid)},\;T\mbox{ (Temperature)}$. \\

It is of course true that plastic melts exhibit viscoelastic behavior. We assume, however, that this is negligible in this context, since the residence time in the twin-screw extruder is small.
We use Generalized Newtonian models to account for shear-thinning and shear-thickening behavior of the polymer melt. Using Generalized Newtonian models implies that the viscosity depends on the invariants of the rate of strain tensor $\boldsymbol{\varepsilon}$, such as the shear rate $\dot{\gamma}$:

\begin{align}
    \dot{\gamma} = \sqrt{2 \boldsymbol{\varepsilon} \left( \bf u \right)  \colon \boldsymbol{\varepsilon} \left( \bf u \right) }.
\end{align}

Within this work we use two different models, namely the Carreau and the Cross-WLF model.\\

The Carreau model is a very popular shear-thinning model in the plastic community \citep{carreau1979review}. It is defined as:

\begin{align}
    \eta \left( \dot{\gamma} \right) = \eta _{\infty}  + \left( \eta_0 - \eta _{\infty} \right)  \left( 1+ \left( \lambda \dot{\gamma} \right) ^2 \right)^{\frac{n-1}{2}},
\end{align}

where $\lambda$ is the relaxation time, $n$ is the power index, $\eta_0 $ is the viscosity at zero shear rate and $\eta _{\infty} $ is the viscosity at infinite shear rate. \\

The Cross-WLF model considers along with the effects of the shear rate also the influence of temperature on the viscosity \citep{rudolph2014polymer}. We neglect the infinite viscosity such that the Cross model is defined as:

\begin{align}
    \eta \left( \dot{\gamma}, T \right) = \frac{\eta_0 \left( T \right)}{1+\left( \frac{\eta_0  \left( T \right) \dot{\gamma}}{\tau ^*} \right)^{(1-n)}},
\end{align}

where $\tau ^*$ is the critical shear stress at the transition from the Newtonian plateau. $\eta \left( \dot{\gamma}, T \right) $ depends on the temperature now. This relation is modeled via the WLF equation:

\begin{align}
    \eta_{0} (T) = D_1 \; exp \left( - \frac{A_1  \left( T - T_{ref} \right) }{A_2 + \left( T - T_{ref} \right) } \right),
\end{align}

with $D_1$ being the viscosity at a reference temperature $T_{ref}$ and $A_1$ and $A_2$ are parameters that describe the temperature dependency. \\

\begin{remark}
In most of the available literature, Stokes equations are used to model the fluid inside a twin-screw extruder. However, Reynolds numbers around 0.1 occur in our applications s.t. the Stokes equations are not necessarily a valid assumption any more. Thus, we model the fluid using the Navier-Stokes equations which also include the advective term. It is noteworthy that the nonlinearity introduced due to the advective term is small compared to the one introduced by using Generalized Newtonian models.
\end{remark}

\subsection{Space-time Finite Element Discretization}

We need to discretize the equations in space and time in order to capture the transient fluid flow. Additionally, the rotating screws introduce a constantly moving and deforming domain. A natural approach accounting for all of these requirements is the DSD/SST (Deforming Spatial Domain / Stabilized Space-Time) method \citep{Tezduyar92a}. In DSD/SST, the weak form is written not only over the spatial domanain, but instead the space-time domain. Thus, we do not need to modify the equations to account for the deforming domain.

In the following, we will define the finite element function spaces for the DSD/SST method. The time interval $(0,T)$ is divided into subintervals $I_n = (t_n, t_{n+1})$, where $n$ defines the time level. We set $\Omega _n = \Omega_{t_n}$ and $\Gamma _n = \Gamma _{t_n}$.  Thus, a space-time slab $Q_n$ is defined as the volume enclosed by the two surfaces $\Omega_n$, $\Omega_{n+1}$ and the lateral surface P$_n$. P$_n$ is described by $\Gamma_t$ as it traverses $I_n$.

We use first-order interpolation for all degrees of freedom. Thus, a SUPG/PSPG stabilization technique is used in order to fulfill the LBB condition \citep{donea2003finite}. The finite-element interpolation and weighting function spaces for velocity, pressure and temperature for every space-time slab are defined as

\begin{align}
&(\mathcal{S}^h_{\bf u})_n = \{   {\bf u}^h \in [H^{1h}(Q_n)]^{n_{sd}} \; | \; {\bf u }^h \doteq {\bf g}^{f,h} \; \text{on} \; (P_n)_{\bf g} \}, \\
&(\mathcal{V}^h_{\bf u})_n = \{   {\bf w}^h \in [H^{1h}(Q_n)]^{n_{sd}} \; | \; {\bf w }^h \doteq {\bf 0} \; \text{on} \; (P_n)_{\bf g} \}, \\
&(\mathcal{S}^h_p)_n = (\mathcal{V}^h_p)_n = \{ p^h \in H^{1h}(Q_n)\}, \\
&(\mathcal{S}^h_T)_n = \{   T^h \in H^{1h}(Q_n) \; | \; T^h \doteq g^{T,h} \; \text{on} \; (P_n)_g \}, \\
&(\mathcal{V}^h_T)_n = \{   v^h \in H^{1h}(Q_n) \; | \; v^h \doteq 0 \; \; \; \text{on} \; (P_n)_g \}.
\end{align}

The stabilized space-time formulation for equations \eqref{eq:cont} and \eqref{eq:momentum} then reads:

Given $({\bf u}^h)^-_n$, find $ {\bf u}^h \in (\mathcal{S} _{\bf u} ^h)_n$ and $p^h \in (\mathcal{S} _p ^h)_n$ such that:

\begin{align}
\label{eq:weakflow}
\begin{split}
\int_{Q_n} {\bf w}^{h} \cdot \rho \left( \frac{\partial {\bf u}^{h}}{\partial t} + {\bf u}^{h} \cdot \boldsymbol{\nabla u}^{h} \right) \; dQ + \int_{Q_n} \boldsymbol{\varepsilon} ( {\bf w}^{h} ) \colon \boldsymbol{ \sigma}^{h} ( p^{h}, {\bf u} ^{h} ) \; dQ \\
+ \int_{Q_n} q^{h}\boldsymbol{\nabla} \cdot {\bf u}^{h}\;dQ + \int_{\Omega_n} ({\bf w}^{h})^+_n \cdot \rho \left( ({\bf u}^{h})^+_n - ({\bf u}^{h})^-_n \right) \;d\Omega \\
+ \sum_{e=1}^{({n_{el})}_n} \int_{Q^e_n} \tau_{\mbox{\tiny{MOM}}} \frac{1}{\rho} \left[ \rho \left( \frac{\partial {\bf w}^{h}}{\partial t}+{{\bf u}^{h} \cdot \boldsymbol{\nabla w}^{h}} \right) + \boldsymbol{\nabla}q^{h} \right] \\
\cdot \left[ \rho \left(\frac{\partial {\bf u}^{h}}{\partial t} + {{\bf u}^{h} \cdot \boldsymbol{\nabla u}^{h}}  \right) - \boldsymbol{\nabla} \cdot {\boldsymbol{\sigma}}^{h}(p^{h}, {\bf u}^{h}) \right]\;dQ \\
+ \sum_{e=1}^{({n_{el})}_n} \int_{Q^e_n} \tau_{\mbox{\tiny{CONT}}} \boldsymbol{\nabla} \cdot {\bf w}^{h} \rho \boldsymbol{\nabla } \cdot {\bf u}^{h}\;dQ
= \int_{P_n} {\bf w}^{h} \cdot {\bf h}^{f,h} \;dP,
\end{split}
\end{align}

holds for all $ {\bf w}^h \in ( \mathcal{V}  _{\bf u} ^h)_n$ and $ q^h \in  (\mathcal{V} _p ^h)_n $.\\

The stabilized space-time formulation for the heat equation \eqref{eq:heat} reads:

Given $ (T^h)^{-}_n $, find $ {T}^h \in (\mathcal{S} _{T} ^h)_n$ such that:

\begin{align}
\label{eq:weaktemp}
\begin{split}
\int_{Q_n} {v}^{h} \cdot \rho c_p \left( \frac{\partial {T}^{h}}{\partial t} + {{\bf u}^{h} \cdot \boldsymbol{\nabla} T^{h}} \right)\;dQ + \int_{Q_n} \boldsymbol{\nabla}v^{h} \cdot \kappa \boldsymbol{\nabla} T^{h} \;dQ \\
- \int_{Q_n} v^{h} \; \phi \;dQ + \int_{\Omega_n} (v^{h})^+_n \rho c_p \left( (T^{h})^+_n - (T^{h})^-_n \right)\;d\Omega\\
+ \sum_{e=1}^{({n_{el})}_n} \int_{Q^e_n} \tau_{\mbox{\tiny{TEMP}}} \frac{1}{\rho c_p} \left[ \rho c_p \left( \frac{\partial v^{h}}{\partial t} + {\bf u}^{h} \cdot \boldsymbol{\nabla}v^{h} \right) \right] \\
\cdot \left[ \rho c_p \left(\frac{\partial T^{h}}{\partial t} + {{\bf u}^{h} \cdot \boldsymbol{\nabla}T^{h}}  \right) - \boldsymbol{\nabla} \cdot \kappa \boldsymbol{\nabla}T^{h} - \phi \right]\;dQ \\
= \int_{P_n} v^{h} h^{T,h} \;dP,
\end{split}
\end{align}

holds for all $ v^h \in (\mathcal{V} _T ^h)_n$.

The following notation is used:

\begin{align}
\begin{split}
\left( {\bf u}^h \right) ^{\pm} _ n = \lim\limits_{\zeta \to 0} {\bf u}^h \left( t_n \pm \zeta \right) \\
\int_{Q_n} .\;.\;.\; dQ = \int_{I_n} \int_{\Omega _{t_n}} .\;.\;.\; d\Omega dt \\
\int_{P_n} .\;.\;.\; dP = \int_{I_n} \int_{\Gamma _{t_n}} .\;.\;.\; d\Gamma dt
\end{split}
\end{align}

The stabilization parameters $\tau_{\mbox{\tiny{MOM}}}$, $\tau_{{\mbox{\tiny{CONT}}}}$ and $\tau_{\mbox{\tiny{TEMP}}}$ are based on expressions given in \citep{pauli2017stabilized}. We use a Newton-Raphson method to solve equation \eqref{eq:weakflow} and \eqref{eq:weaktemp}. The flow and temperature fields are coupled strongly using a fixed-point iteration until convergence is reached.
The stress contributions in the SUPG/PSPG stabilization terms (fifth term in equation \eqref{eq:weakflow} and fourth term in equation \eqref{eq:weaktemp} ) are zero, since they involve second-order derivatives. We improve the consistency of our method by employing a least-squares recovery technique for these terms \citep{jansen1999better}.
In order to solve the resulting linear system of equations within each Newton iteration, we use a GMRES solver with an ILUT preconditioner.

\section{Numerical Examples} \label{sec:numericalexamples}

\subsection{2D Isothermal Flow}

We want to show that simulations based on the new method for grid generation produces the same high-quality results as already established methods like SRMUM or XFEM. Therefore, we simulate the isothermal flow of a plastic melt, described by the Carreau model, in a 2D cross section of a twin-screw extruder. The Carreau parameters are given in Table \ref{table:carreau2D}. The screw geometry is generated based on an adapted version of Booy's description as presented in \cite{sarhangi2012adaptive,helmig2018boundary}. The screw parameters are given in Table \ref{table:screw2D}. The rotation speed of the screws is $\omega _s = 60$ rpm in mathematically positive direction.
We set the rotational velocity on the screws as a Dirichlet boundary condition and a no-slip condition on the barrel. In \cite{sarhangi2012adaptive}, Stokes flow was used. Since we use Navier-Stokes equations to model the flow, we employ a density of $\rho = 1 \; kg/m^3$ in order to make the results comparable.

\begin{table}[h]
  \begin{minipage}[t]{0.5\linewidth}
  \centering
  \begin{tabular}{l c c}
  \hline
  $\eta_0$  & 1290 & $Pa \; s$ \\
  $\eta_{\infty}$ & 0 & $Pa \; s$ \\
  $n$ & 0.559 & - \\
  $\lambda$ & 0.112 & $s$ \\
  \hline
\end{tabular}
\caption{Carreau parameters.}
\label{table:carreau2D}
\end{minipage}
\begin{minipage}[t]{0.5\linewidth}
\centering
\begin{tabular}{l r}
\hline
Screw radius $R_s$  & 15.275 $mm$ \\
Center line distance $C_l$ & 26.2 $mm$ \\
Screw-screw clearance $\delta _s$ & 0.2 $mm$ \\
Screw-barrel clearance $\delta _b$ & 0.15 $mm$ \\
\hline
\end{tabular}
\caption{Screw geometry parameters.}
\label{table:screw2D}
\end{minipage}
\end{table}

We employ a time step size of $\Delta t = 0.00625 s$ or $2.25 \; ^\circ / s$. The resulting flow can be considered as quasi-steady or instantaneous because the time-scales of the momentum diffusion are very small compared to the process itself. The resulting Reynolds number inside the small gap region is $Re = 0.000003$.  Therefore, it is sufficient to compare the solution at individual screw orientations. Two screw orientations, namely $\theta = 0 ^\circ $ and $\theta = 112.5 ^\circ $ have been used in \cite{helmig2018boundary}. It is demonstrated that $\theta = 112.5 ^\circ $ is the more complex orientation. Thus, we will only compare results for that orientation.

\begin{table}
\centering
\begin{tabular}{l c c c c c c c}
 \hline
 mesh & $n_s$ C-grid & $n_s$ separator & $n_s$ total & $n_r$ & $\#$ elements & $n_{\mu}$ & $n_{\nu}$ total \\
 1 & 200 & 80 & 280  & 6 & 3360 & 12 & 280\\
 2 & 300 & 160 & 460 & 10 & 9200 & 12 & 460\\
 3 & 300 & 300 & 600 & 12 & 14400 & 12 & 600 \\
 4 & 300 & 600 & 900 & 18 & 32400 & 12 & 900\\
 \hline
 \end{tabular}
 \caption{Mesh discretization for 2D convergence study.}
 \label{table:mesh2DConvergence}
\end{table}

\begin{figure}[h]
  \centering
  \subfigure[Mesh and plot over line.]{\includegraphics[width=.43\linewidth]{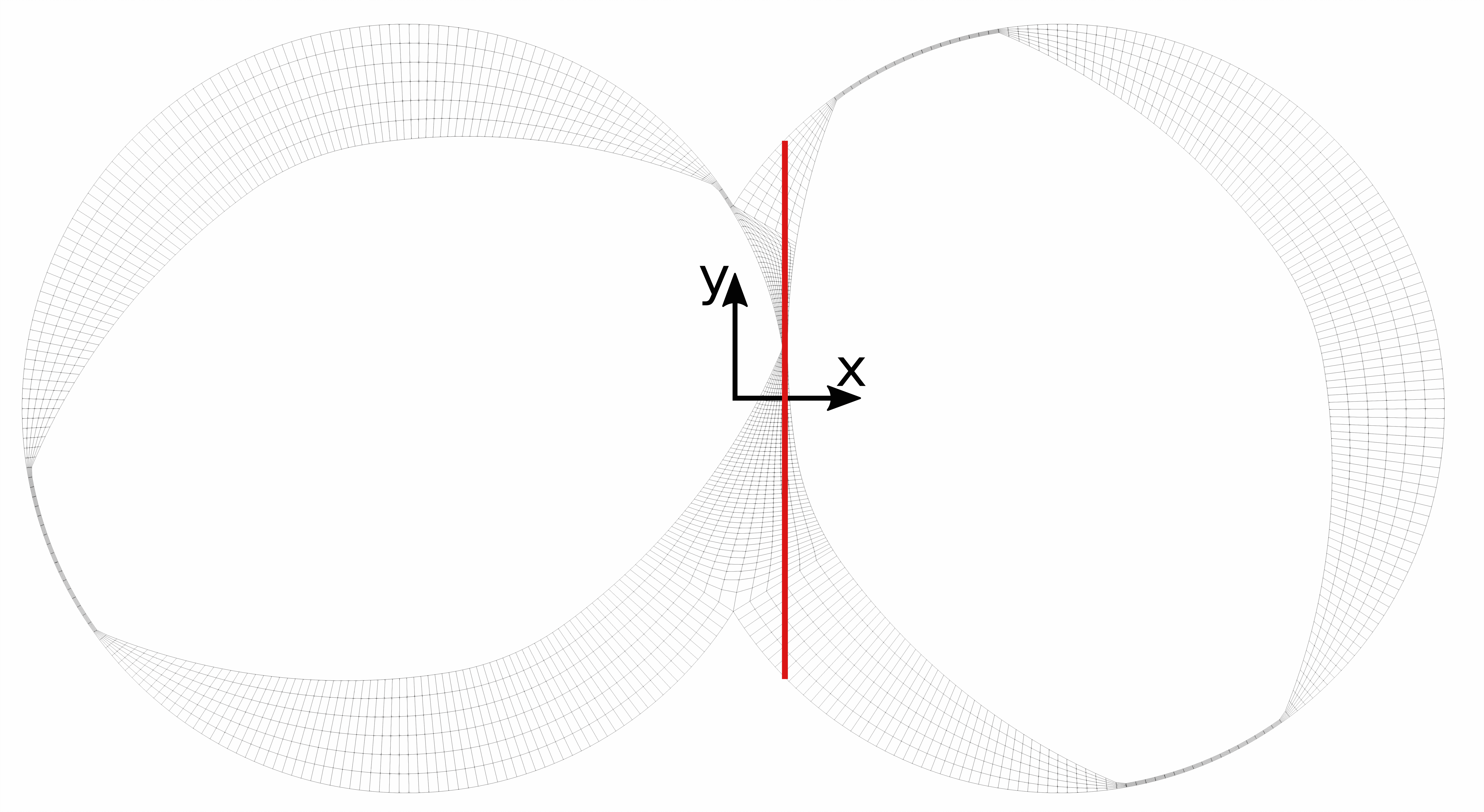}\label{fig:domain2d}}
  \centering
  \subfigure[Pressure distribution.]{\includegraphics[width=.43\linewidth]{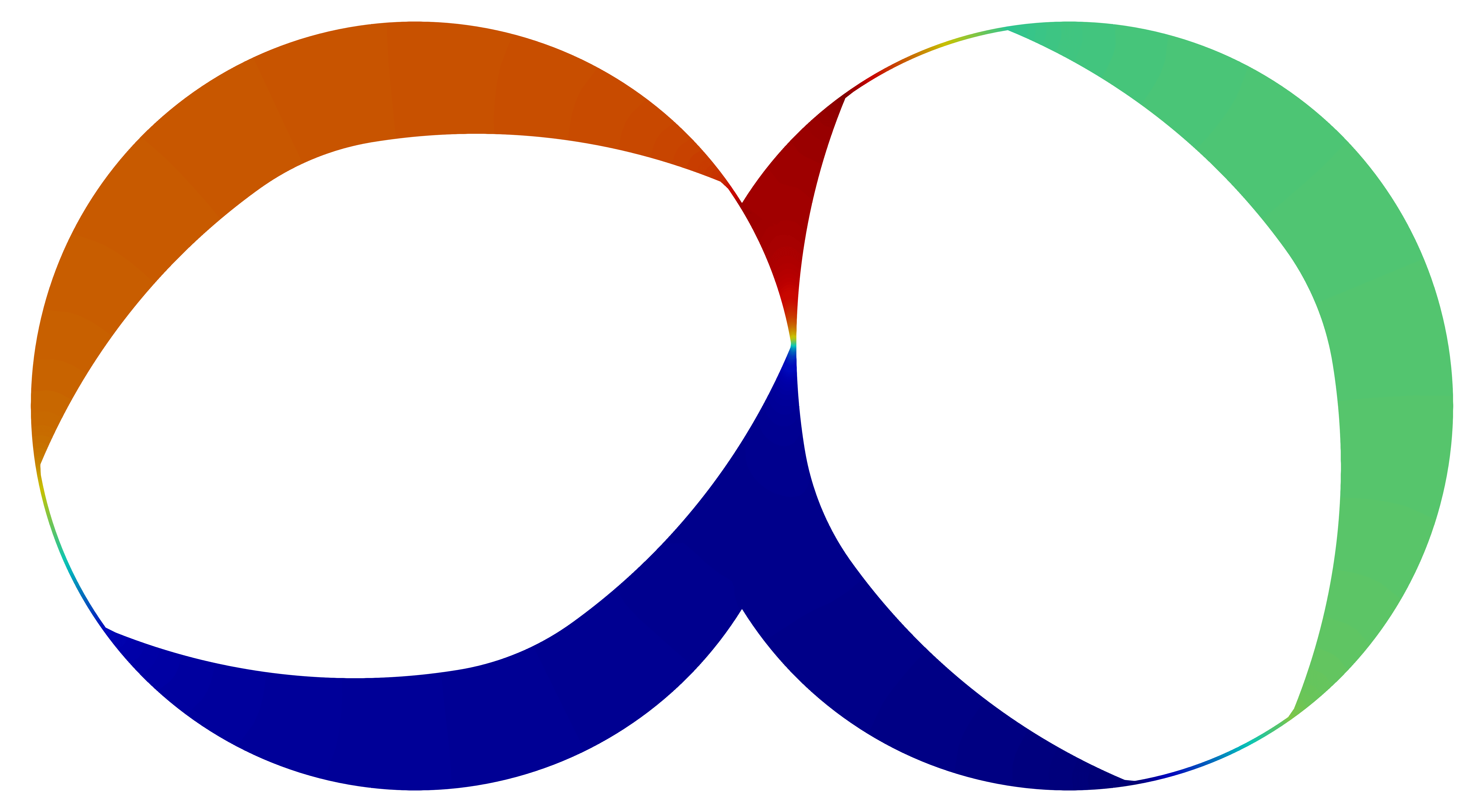}\label{fig:pressure2d}}
  \centering
  \subfigure{\includegraphics[width=.085\linewidth]{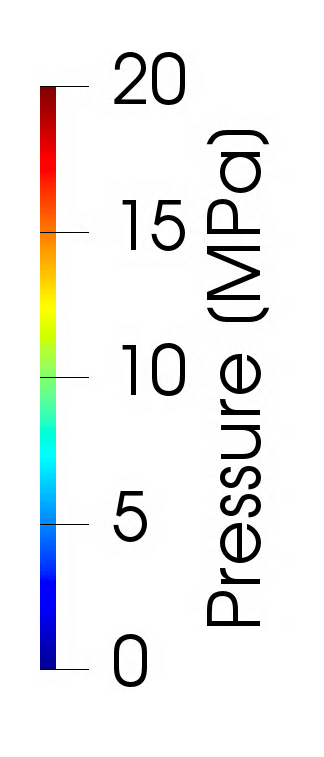}}
  \caption{Mesh as well as pressure results for orientation $\theta = 112.5^{\circ}$.}
  \label{fig:angle1125}
\end{figure}

\begin{figure}[h]
  \centering
  \subfigure[]{\includegraphics[width=.45\linewidth]{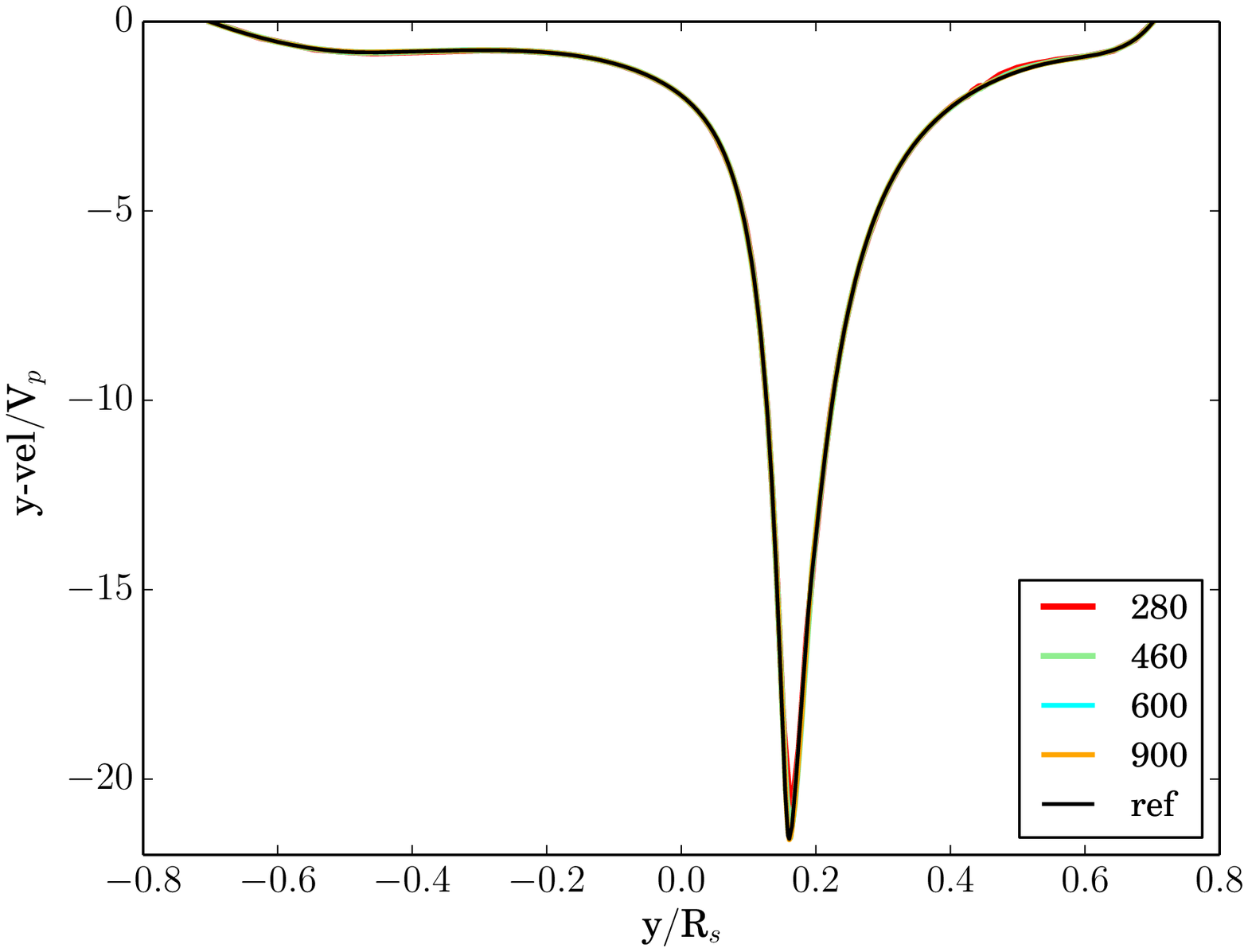}\label{fig:angle1125Veloverall}}
  \centering
  \subfigure[]{\includegraphics[width=.45\linewidth]{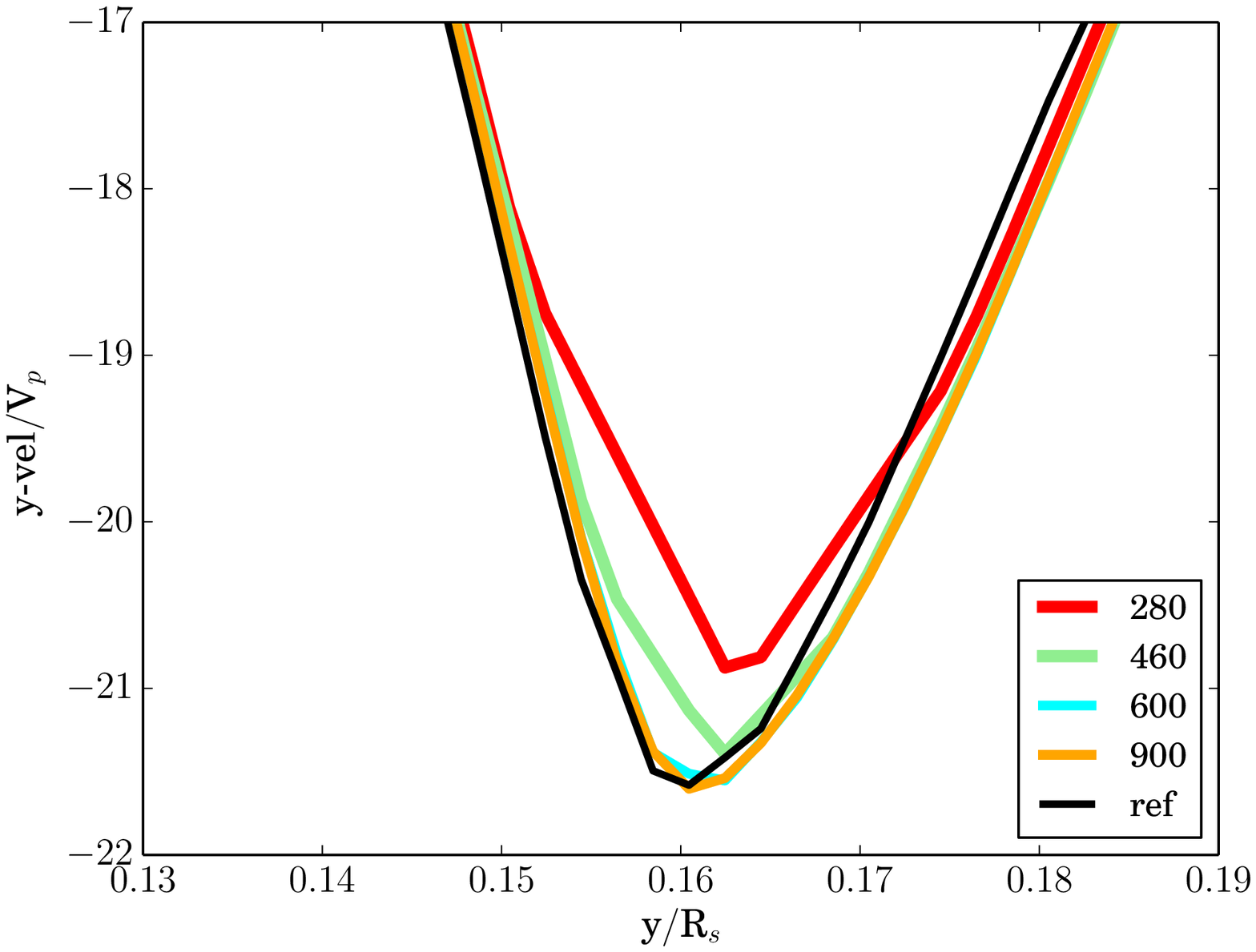}\label{fig:angle1125Velclose}}
  \caption{Velocity plots over line for orientation $\theta = 112.5^{\circ}$ -- (a) velocity profile over entire line (b) close view. The numbers in the legend denote meshes based on $n_s$ total, see Table \ref{table:mesh2DConvergence}. The reference solution $ref$ is computed with the time step based on a SRMUM mesh with 1000 elements in screw and 20 elements in radial direction.}
  \label{fig:angle1125Vel}
\end{figure}

In order to verify the general applicability of the method, we use 4 different mesh resolutions. One advantage of the presented method is that we have to construct the spline parameterization only once and can then extract the \textit{'scaffolding'} for all desired mesh resolutions.
The region of special interest is the intermeshing area between the two screws. A high pressure drop drives the flow solution resulting in high velocities and high shear rates, see Figure \ref{fig:pressure2d}. Thus, resolving this part correctly is extremely important. Therefore, we only refine our mesh inside this region and keep the mesh resolution inside the two C-grids constant. The only exception is the first mesh, where we slightly reduce the resolution for the C-grid in order to avoid highly stretched elements. We evaluate the spline-parameterization of the C-grids as well as the separator for all mesh points on the screw, meaning $n_s = n_{\mu}$. Inside the separator we use a \textit{'scaffoling'} with 12 elements in $\mu$-direction. This number is kept constant throughout the refinement of the background mesh.
All mesh quantities are given in Table \ref{table:mesh2DConvergence}. In the following, we will refer to them by the total number of elements in screw direction '$n_s$ total'.

We compare the velocity results inside the intermeshing area along a line in $y$-direction at $x = 2.077 \;mm$. The line is visualized in red for mesh 280 in Figure \ref{fig:domain2d}. As a reference solution, we use flow results computed with the same time step size on a SRMUM mesh with 1000 elements in screw and 20 elements in radial direction.
The plot over line for the normalized y-velocity $y_{vel}/V_p$, with $V_p = 2 \pi R_s \omega_s$ is shown in Figure \ref{fig:angle1125Vel}. Looking at the overall plot over line in Figure \ref{fig:angle1125Veloverall}, we can not observe major differences between the results computed on different meshes compared to the reference solution.
The velocities differ by less than $1 \; \%$. However, the peak velocities differ more. Figure \ref{fig:angle1125Velclose} shows that this difference decreases with increasing mesh resolution. The results are also in accordance with the results presented in \cite{sarhangi2012adaptive}.

\subsection{2D Temperature-Dependent Flow inside Mixing Elements}

\begin{figure}[h]
  \centering
  \subfigure[Config. 1]{\includegraphics[width=.45\linewidth]{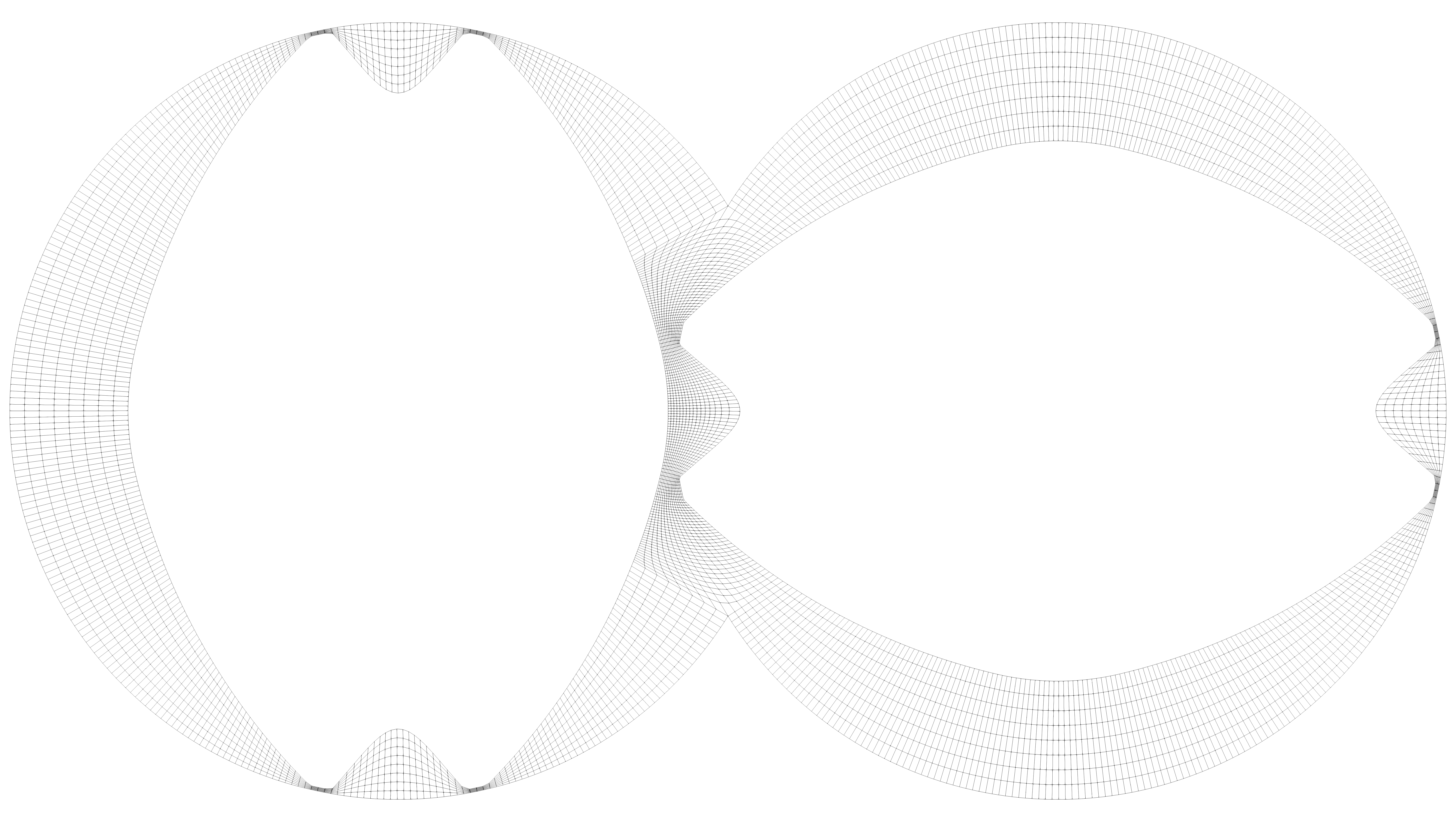}\label{fig:domainMixing}}
  \centering
  \subfigure[Config. 2]{\includegraphics[width=.45\linewidth]{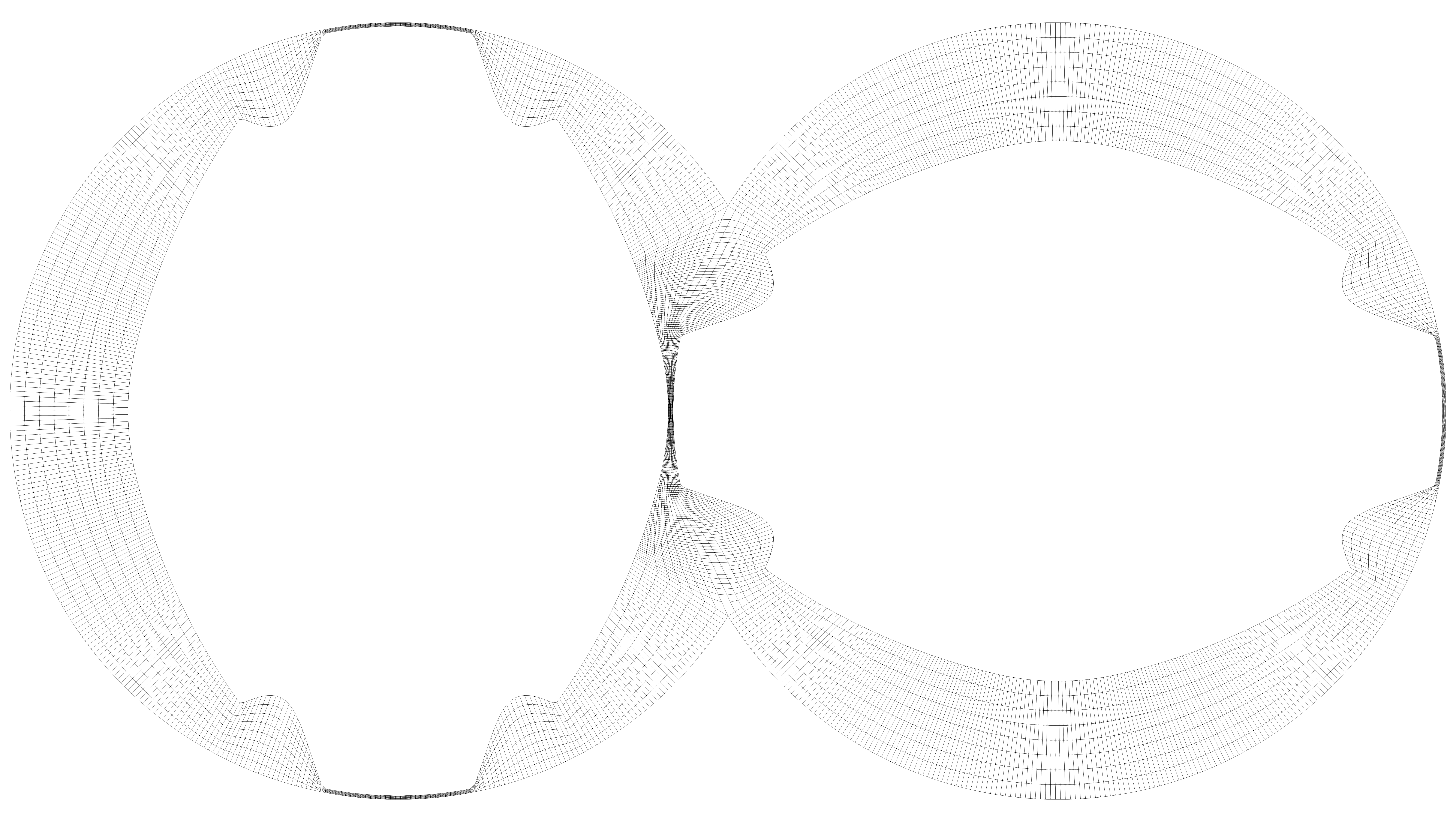}\label{fig:domainMixingCut}}
  \caption{Two 2D screw configurations for different mixing elements.}
  \label{fig:domainMixing2D}
\end{figure}

In the previous section, we have shown that the meshes based on the spline-based parameterization technique produce valid results. However, the true advantage over SRMUM only becomes apparent in the context of non-convex screw shapes - shapes that SRMUM cannot handle at all. We use two different screw configurations that are inspired by the screw design that has already been given in Table \ref{table:screw2D}. The resulting screw geometries are given in Figure \ref{fig:domainMixing2D}.

\begin{table}[h]
\centering
\begin{tabular}{l c c c c c c c c c}
 \hline
 mesh & $n_{s}$ C-grid & $n_{s}$ separator & $n_s$ total & $n_r$ & $\#$ elements & $n_{\mu}$ & $n_{\nu}$ total & $n_{slices}$ $\pi$ & memory savings \\
 1 & 150 & 60 & 210  & 4 & 1680 & 12 & 210& 101 & 22 $\%$ \\
 2 & 300 & 120 & 420  & 8 & 6720 & 12 & 420& 101 & 61 $\%$\\
 3 & 600 & 240 & 840  & 16 & 26880 & 12 & 840& 101 & 81 $\%$\\
 \hline
 \end{tabular}
 \caption{Mesh discretizations for 2D mixing elements for Config. 1.}
 \label{table:mesh2DMixingConfig1}
\end{table}

\begin{table}[h]
\centering
\begin{tabular}{l c c c c c c c c c}
 \hline
 $n_{s}$ C-grid & $n_{s}$ separator & $n_s$ total & $n_r$ & $\#$ elements & $n_{\mu}$ & $n_{\nu}$ total & $n_{slices}$ $\pi$ & memory savings \\
 400 & 140 & 540 & 8 & 8640 & 12 & 540& 101 & 63 $\%$ \\
 \hline
 \end{tabular}
 \caption{Mesh discretization for 2D mixing elements for Config. 2.}
 \label{table:mesh2DMixingConfig2}
\end{table}

\begin{figure}[h]
  \centering
  \subfigure[$0.0^{\circ}$]{\includegraphics[width=.22\linewidth]{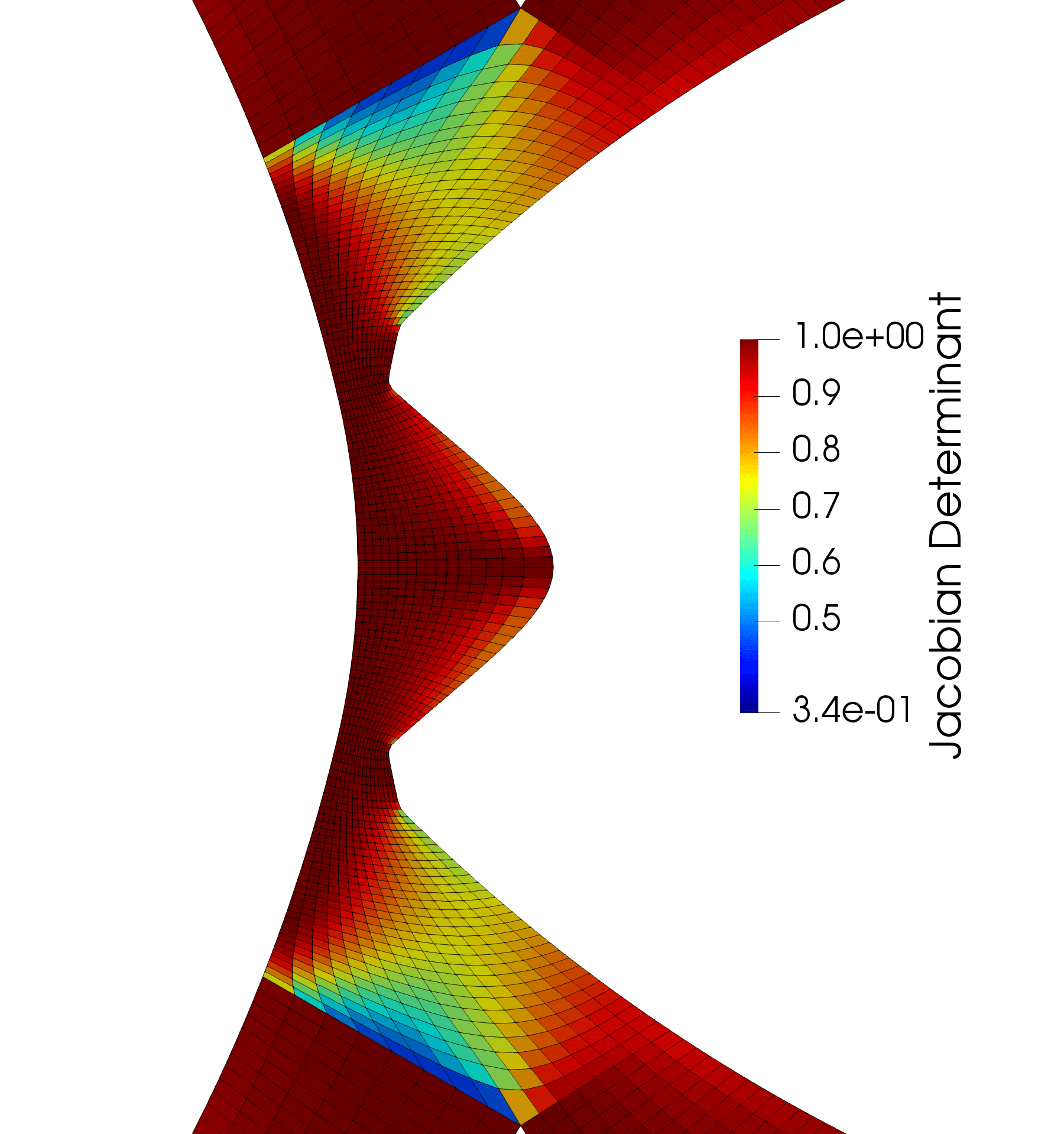}}
  \centering
  \subfigure[$32.4^{\circ}$]{\includegraphics[width=.22\linewidth]{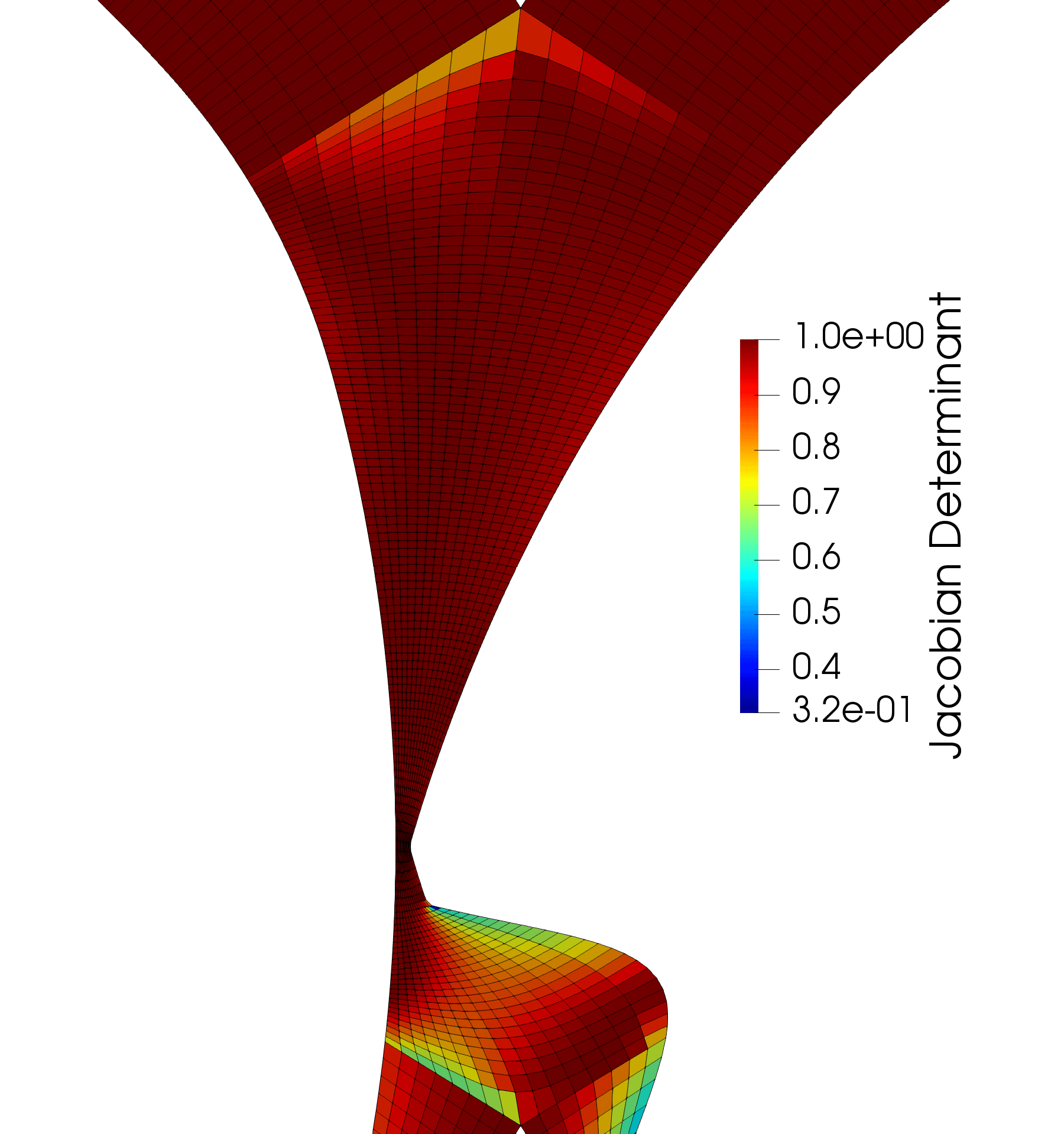}}
  \centering
  \subfigure[$129.6^{\circ}$]{\includegraphics[width=.22\linewidth]{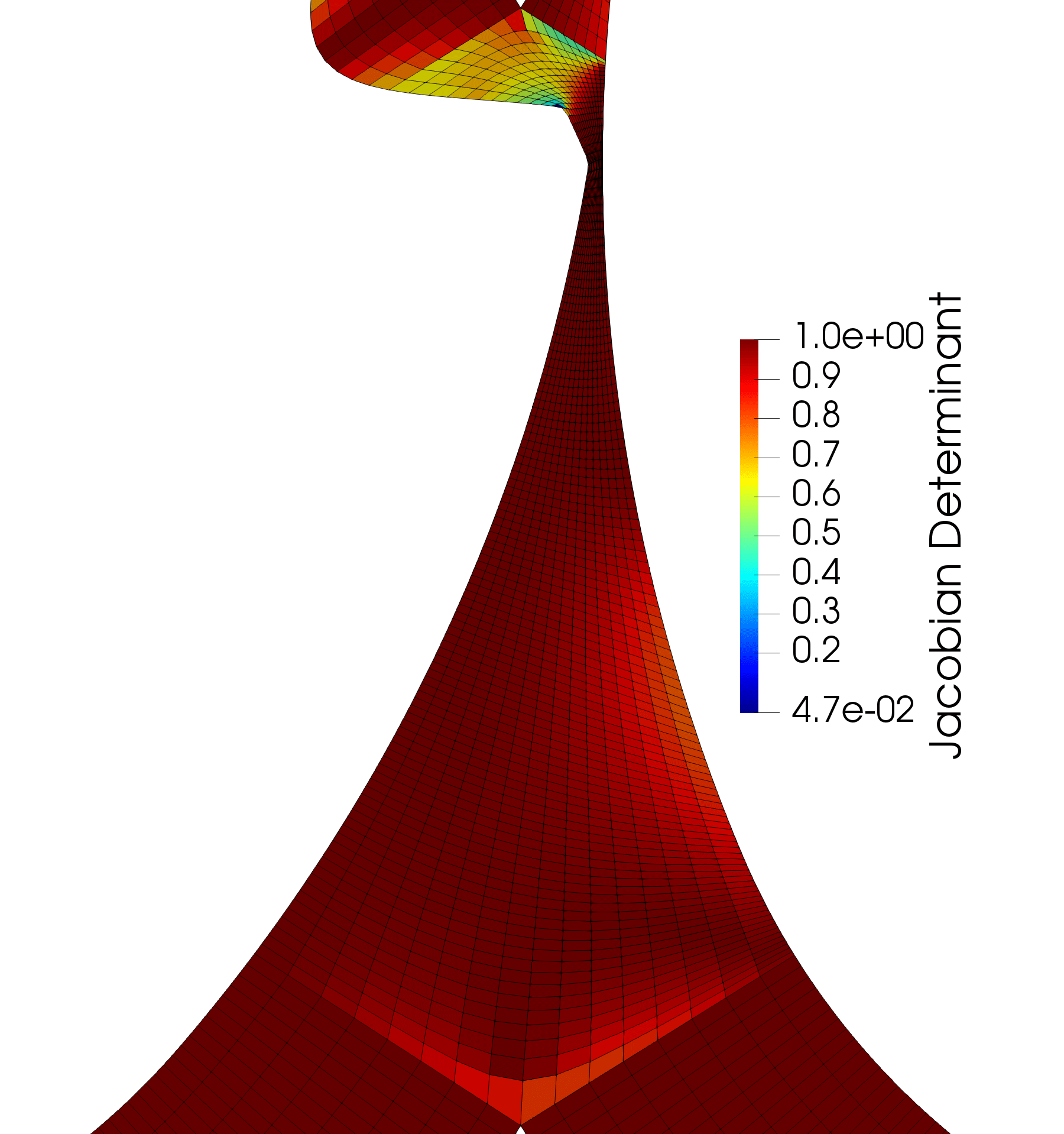}}
  \centering
  \caption{Scaled Jacobian determinant inside the separator or intermeshing area for Config. 1 for different angles on mesh 2.}
  \label{fig:jacdetconfig1}
\end{figure}
\begin{figure}[h]
  \centering
  \subfigure[$0.0^{\circ}$]{\includegraphics[width=.22\linewidth]{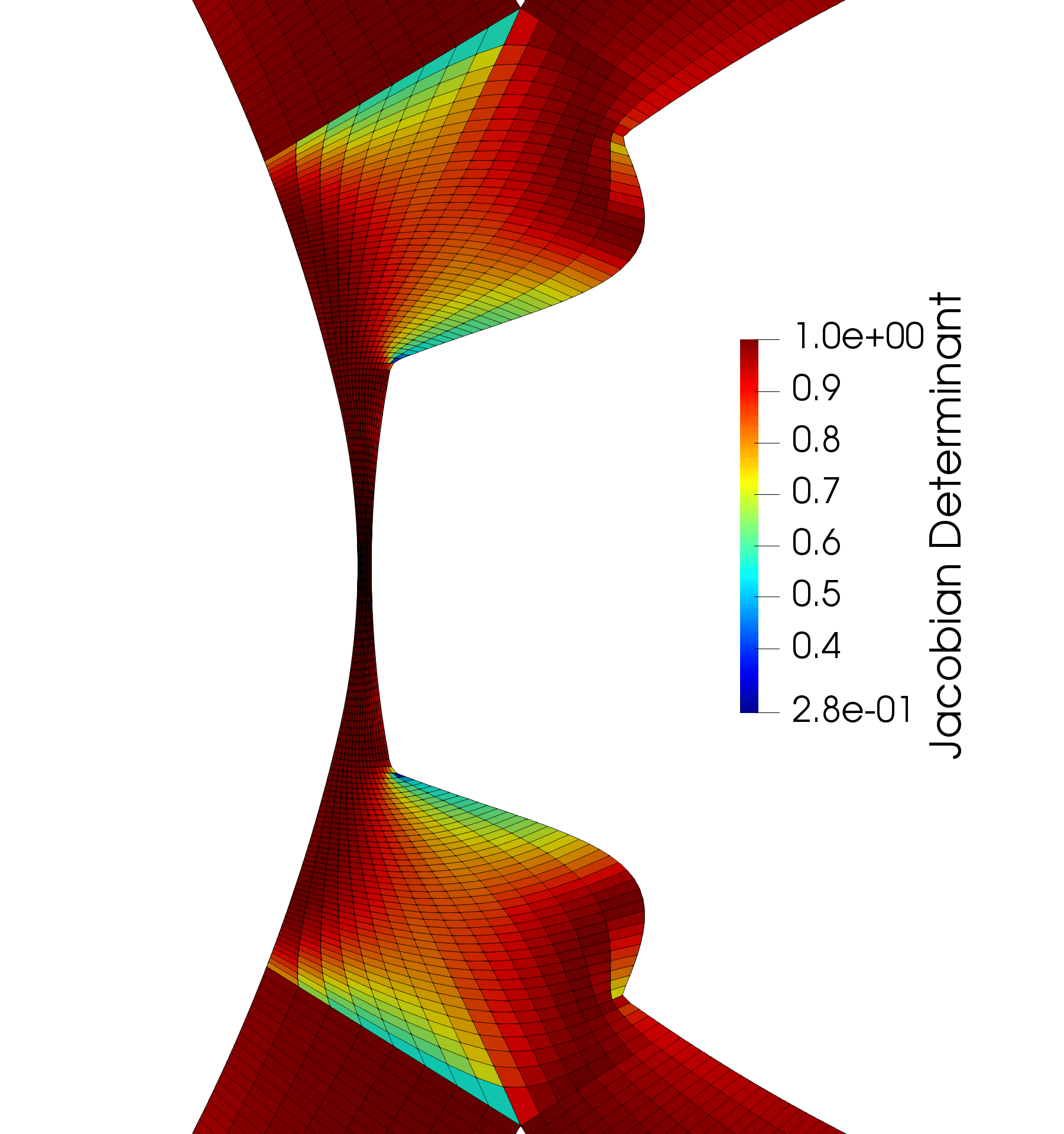}}
  \centering
  \subfigure[$32.4^{\circ}$]{\includegraphics[width=.22\linewidth]{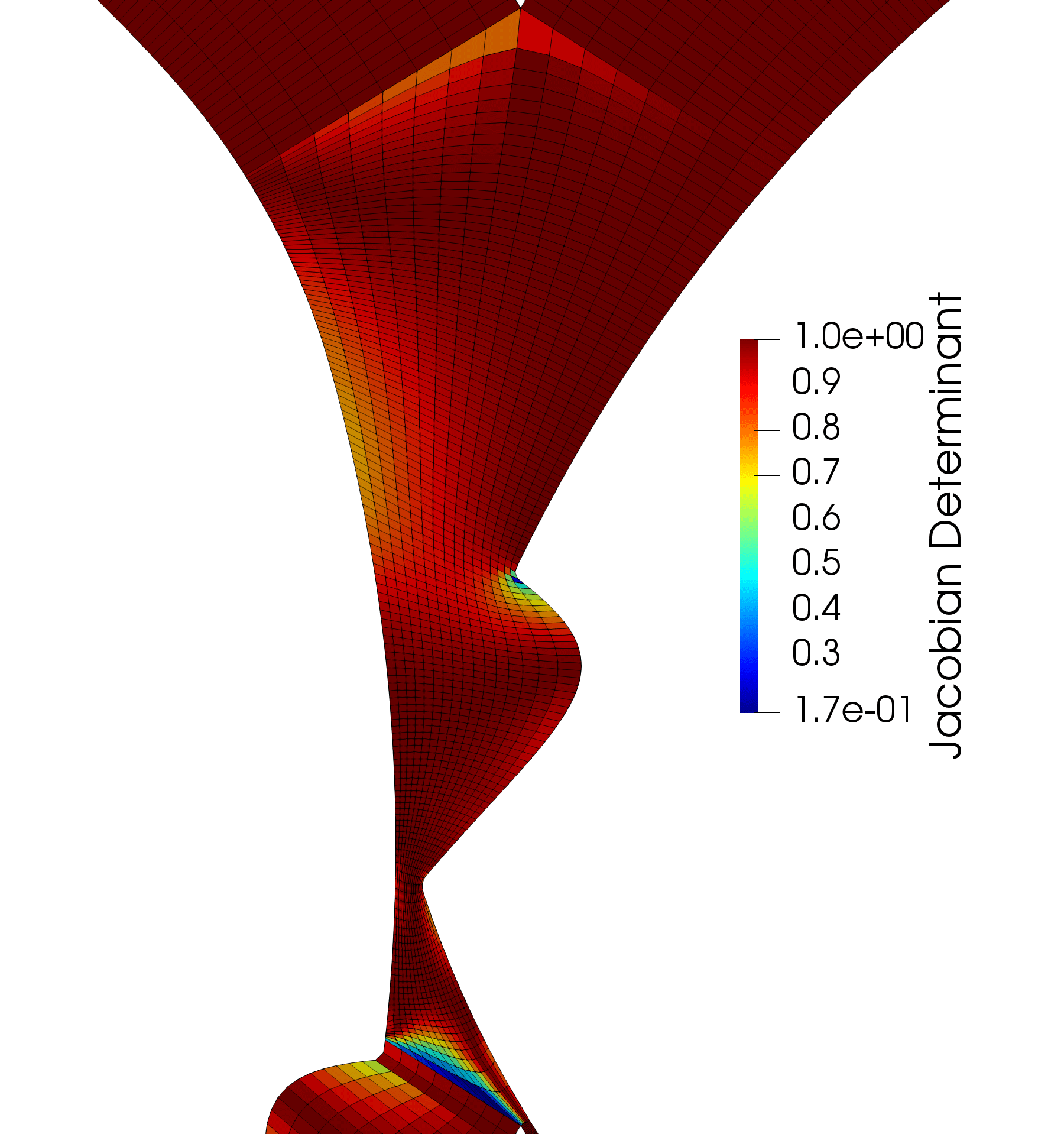}}
  \centering
  \subfigure[$129.6^{\circ}$]{\includegraphics[width=.22\linewidth]{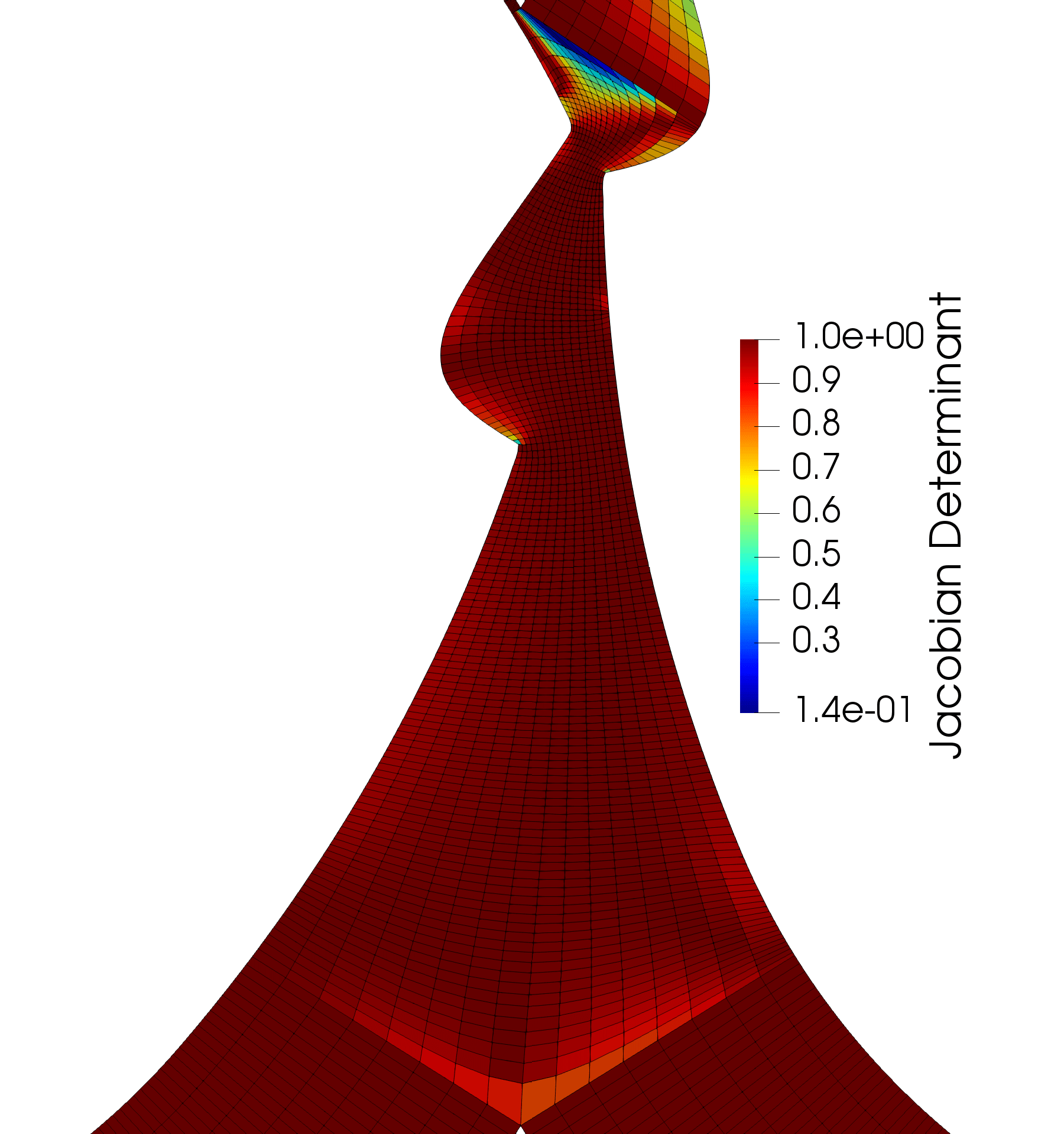}}
  \centering
  \caption{Scaled Jacobian determinant inside the separator or intermeshing area for Config. 2 for different angles.}
  \label{fig:jacdetconfig2}
\end{figure}

Config. 2 can be considered slightly more complex since the screw geometry has sharper edges compared to Config. 1. We employ several meshes for Config. 1 that will be later used in order to perform a mesh convergence study. Mesh details are given in Table \ref{table:mesh2DMixingConfig1} and \ref{table:mesh2DMixingConfig2}, including the memory savings obtained by only storing a \textit{'scaffolding'}.
Note, that the mesh for Config. 2 is slightly finer than the second mesh of Config. 1.
This is due to the more complex shape of the screw profile, which necessitates using more grid points to properly capture the sharpe edges. \\ Similar to the previous section, we evaluate the spline parameterization along with the computed control mapping. Also in this non-convex case, the reparameterization is used to ensure that the C-grid parameterization with transfinite interpolation does not fold. Inside the separator, we again evaluate the spline control mapping composition in a uniform grid comprised of 12 elements in $\mu$-direction and $n_{\nu}=n_s$ in $\nu$-direction. The resulting interpolated fine meshes inside the separator are given for different angles for Config. 1 in Figure \ref{fig:jacdetconfig1} and Config. 2 in Figure \ref{fig:jacdetconfig2}.
In order to evaluate the mesh quality, we use the scaled Jacobian determinant calculated by the \textit{'Mesh quality'} filter in paraview \cite{Ayachit:2015:PGP:2789330}. We obtain highly distorted elements in particular in the vicinity of concave corners, but the scaled Jacobian determinant is never negative.
This is extremely important because the finite-element discretization can cope with highly distorted elements but fails in case of negative determinants. \\

\begin{table}
  \centering
  \begin{tabular}{l c c}
  \hline
  $D1$  & 5.0e+13 & $Pa \; s$ \\
  $\tau^{*}$ & 2500.0 & $Pa$ \\
  $n$ & 0.29 & - \\
  $T_{ref}$ & 263.15 & $K$ \\
  $A1$ & 28.32 & - \\
  $A2$ & 51.60 & $K$ \\
  \hline
\end{tabular}
\caption{Cross-WLF parameters.}
\label{table:CrossWLF2D}
\end{table}

In the following, we aim to demonstrate the advantage of the presented meshing method. In case of isothermal flow, it is not crucial to have matching discretization for consecutive screw orientations due to the quasi-steady behavior of the flow. However, taking temperature effects into account results in equations where a solution at the next time step strongly depends on the previous one. Therefore, we simulate a temperature-dependent plastic melt inside the two 2D mixing elements. The melt is modeled using the Cross-WLF model. The parameters are given in Table \ref{table:CrossWLF2D}. Note, that the parameters have been selected for testing purposes but are in inspired by those of polypropylene.
The plastic melt has density $\rho = 700 \; kg / m^3$, specific heat $c_p = 2400 \; J/ (kg \;K)$ and thermal conductivity $\kappa = 10.0 \; W/(m\;s)$. The screws rotate in mathematically positive direction with $\omega = 60$ rpm. For the flow, we set a no-slip condition on the barrel and the rotational velocity as Dirichlet condition on the screws.

Concerning the temperature, the screws are considered to be adiabatic. On the barrel we set a Dirichlet temperature condition as $T_{barrel} = 473 K + (10x)/0.03 K/m$ and by that, generate a linear increase of the temperature from the left to the right barrel.
The initial condition is $T_0 = 473 K$. It is important to note that this constitutes a test case that has no connection to a real twin-screw extruder application.
It merely demonstrates the validity of the method, as cold/hot melt is constantly pushed from left to right and vice versa.
We use different time step sizes to show that the solutions are independent of the time step. Thus, we have to interpolate linearly between slices for most of the resulting screw orientation. Assuming that each individual slice is bijective, a linear interpolation in $z$-direction will lead to a bijective volumetric mesh.\\

\subsubsection*{Config. 1:}

\begin{figure}
  \centering
  \subfigure[Refinement in space with $\Delta t = 0.0625 s$.]{\includegraphics[width=.45\linewidth]{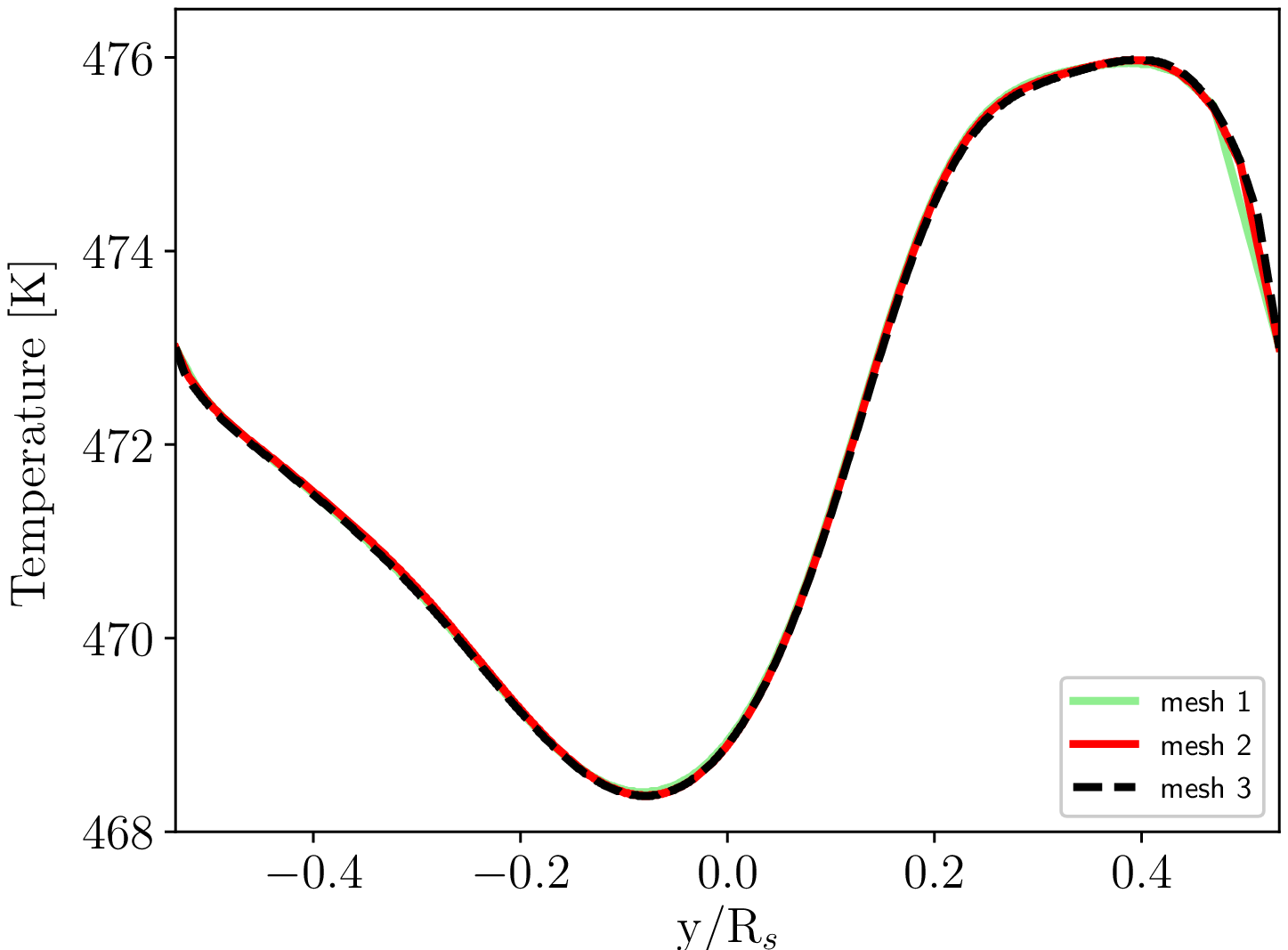}\label{fig:refinmentSpaceNoViscDiss}}
  \centering
  \subfigure[Refinement in time on mesh 2.]{\includegraphics[width=.45\linewidth]{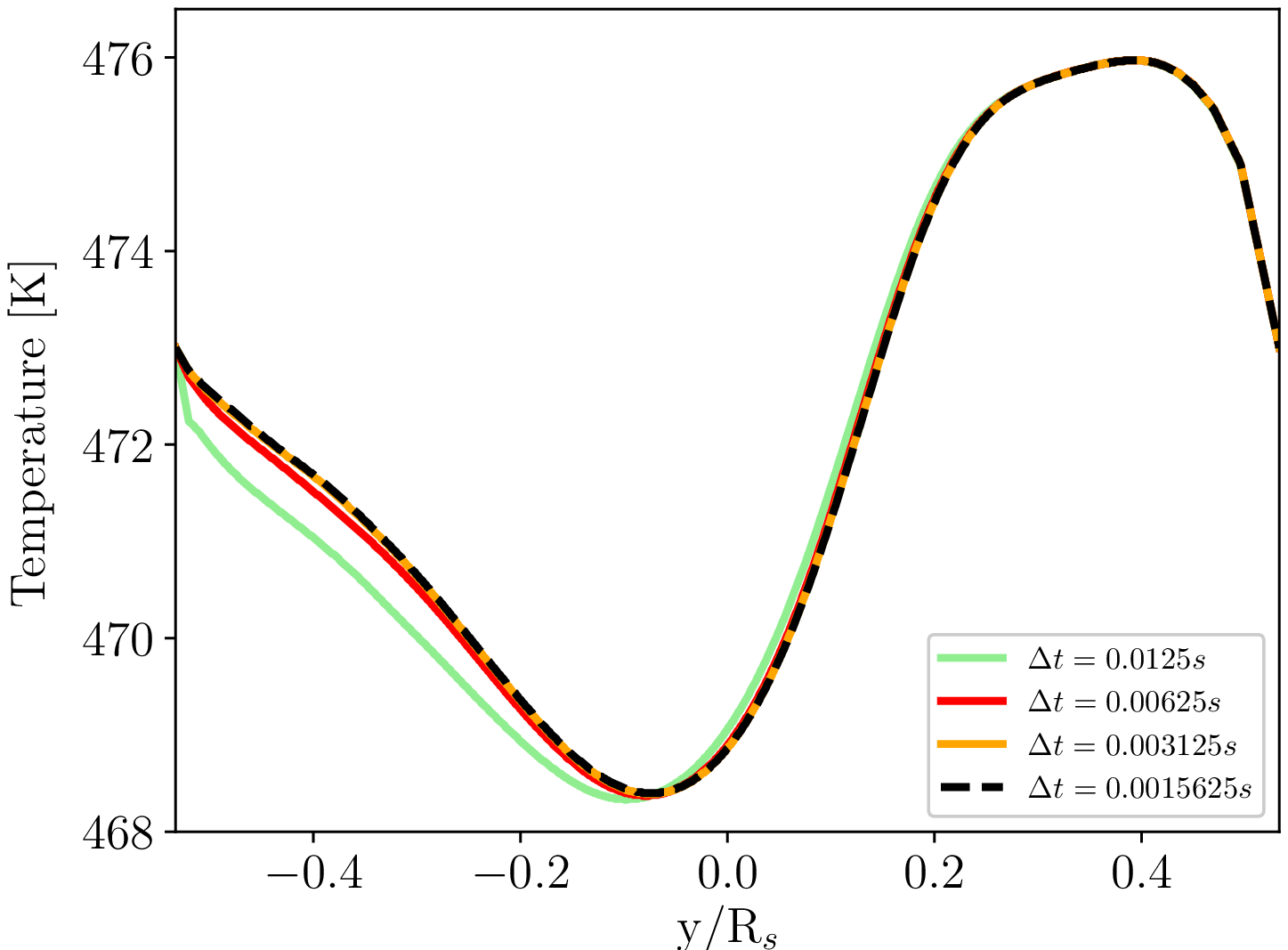}\label{fig:refinmentTimeNoViscDiss}}
  \centering
  \caption{Mesh and time step refinement study for Config. 1 without viscous dissipation. The temperature results are compared along a line from the lower to upper cusp point at time $t = 1.625 s$.}
  \label{fig:tempMixCut2DNoViscDissRef}
\end{figure}

\begin{figure}
  \centering
  \subfigure[$t=1.625 s$]{\includegraphics[width=.38\linewidth]{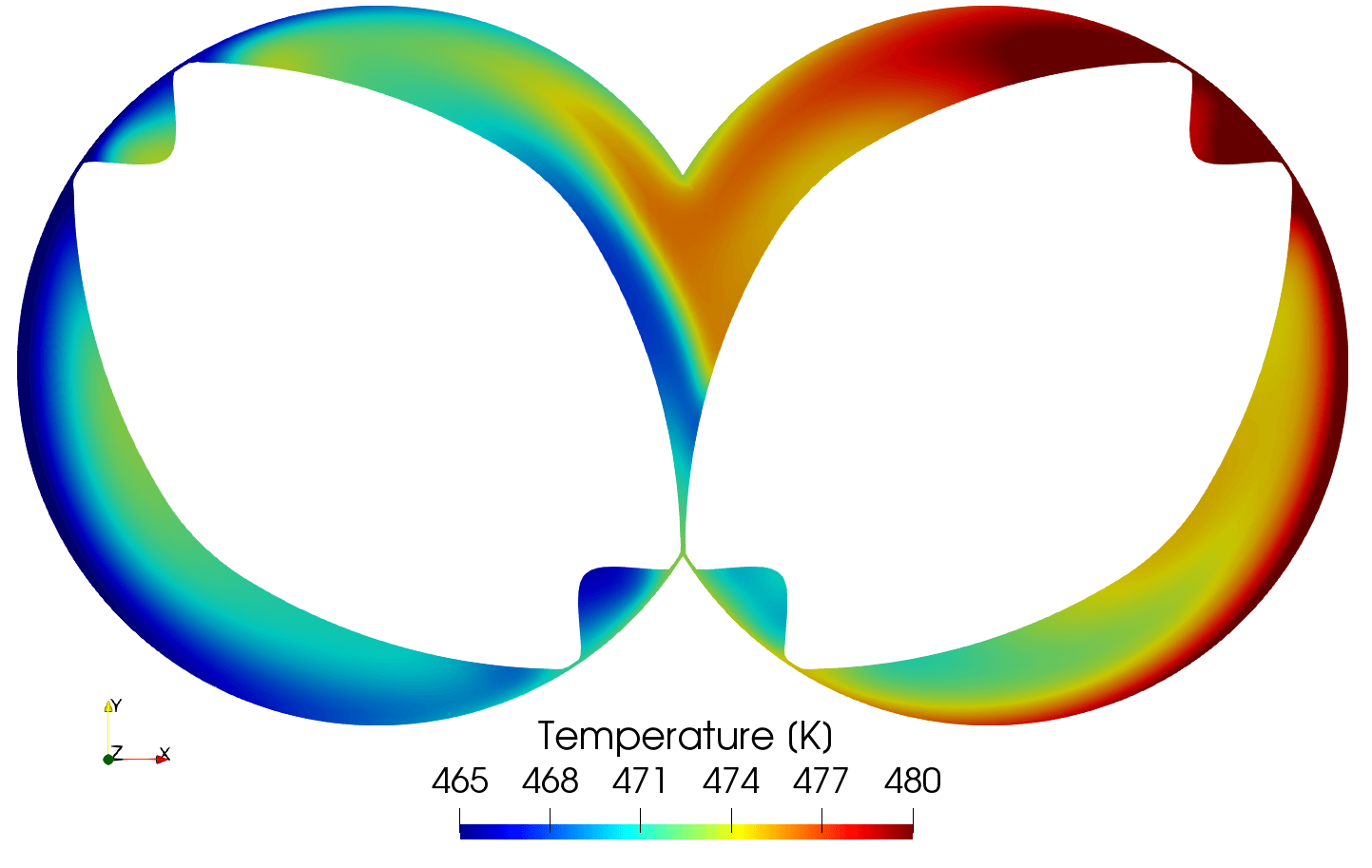}\label{fig:tempMixCut2DNoViscDissa}}
  \centering
  \subfigure[$t=1.675 s$]{\includegraphics[width=.38\linewidth]{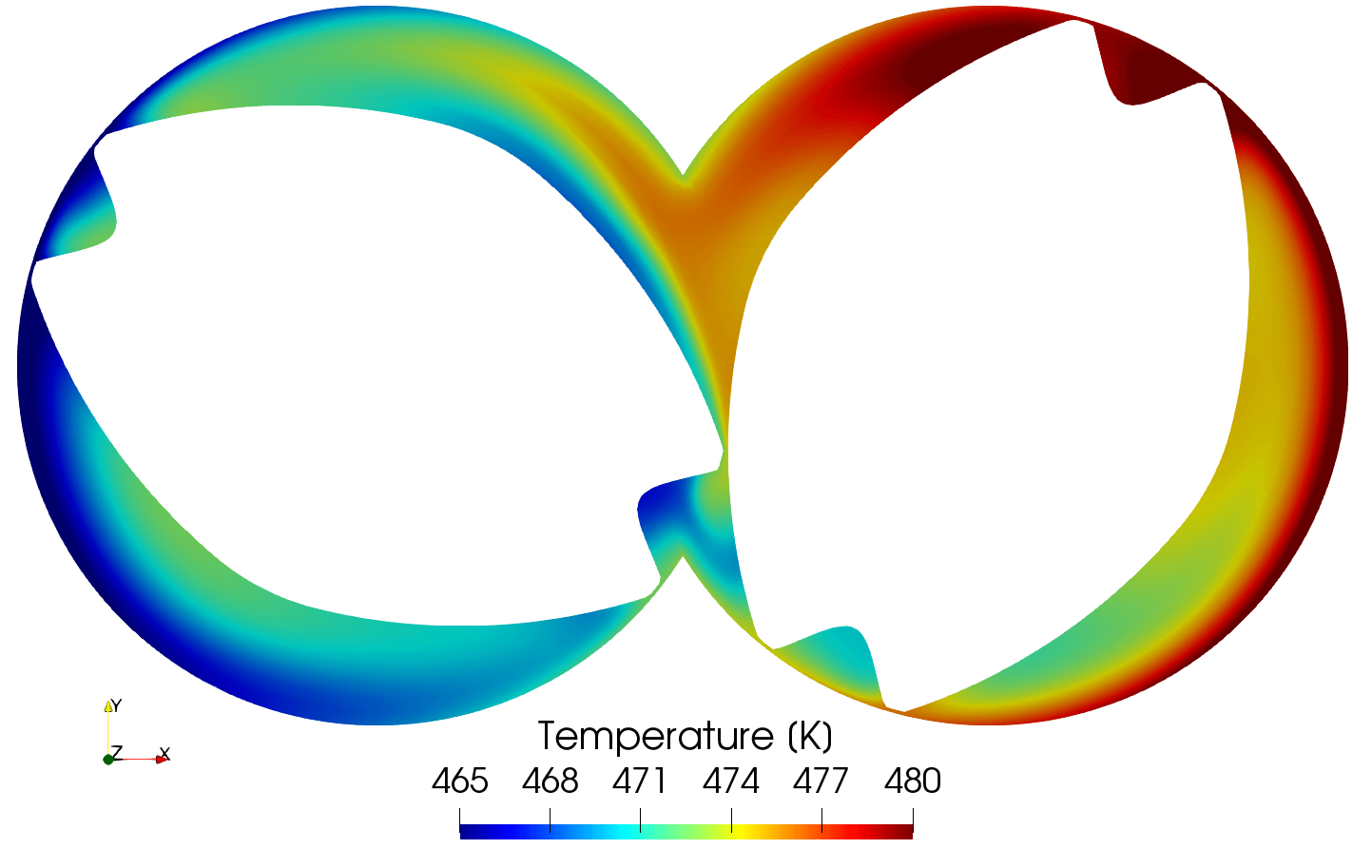}\label{fig:tempMixCut2DNoViscDissb}}
  \centering
  \subfigure[$t=1.725 s$]{\includegraphics[width=.38\linewidth]{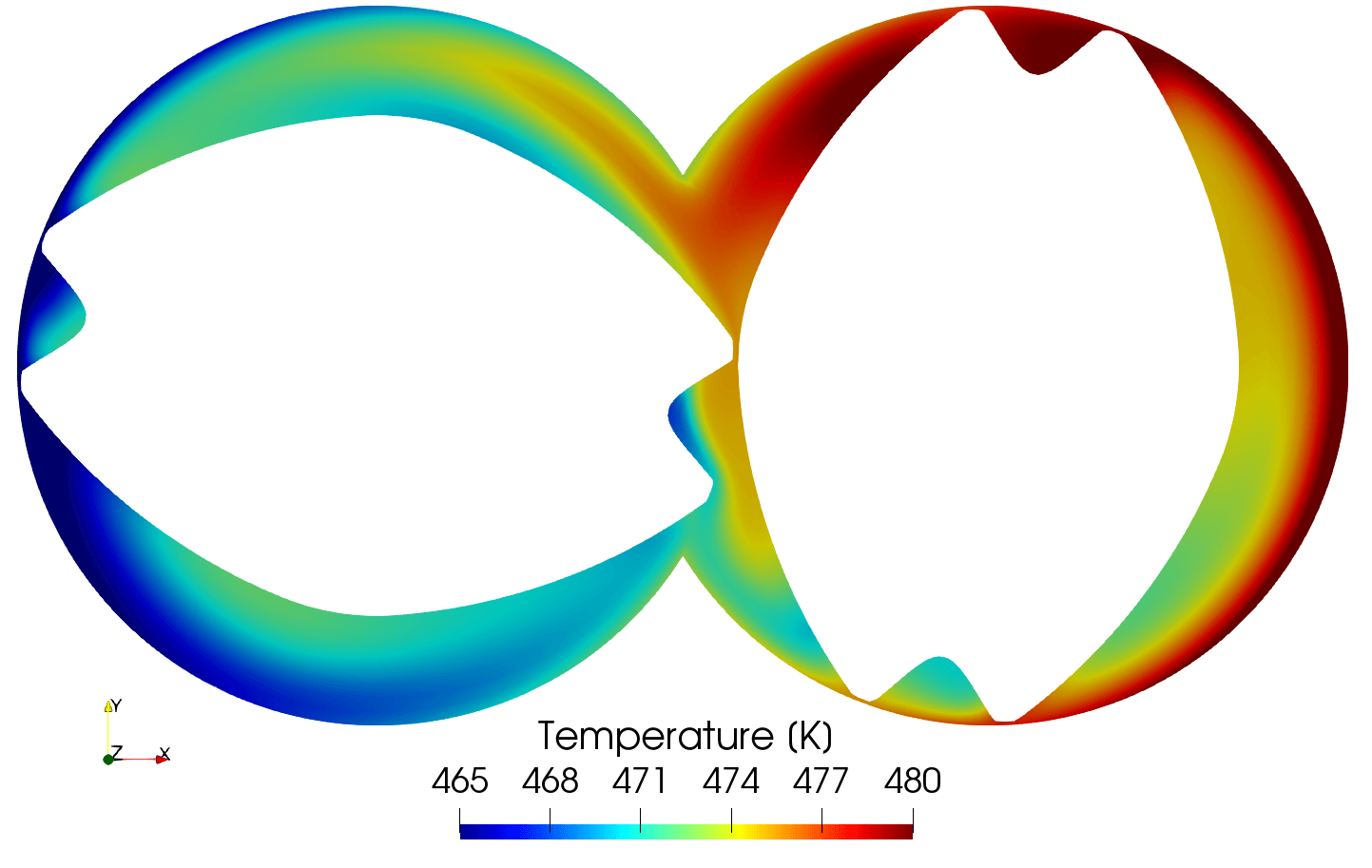}\label{fig:tempMixCut2DNoViscDissc}}
  \centering
  \subfigure[$t=1.775 s$]{\includegraphics[width=.38\linewidth]{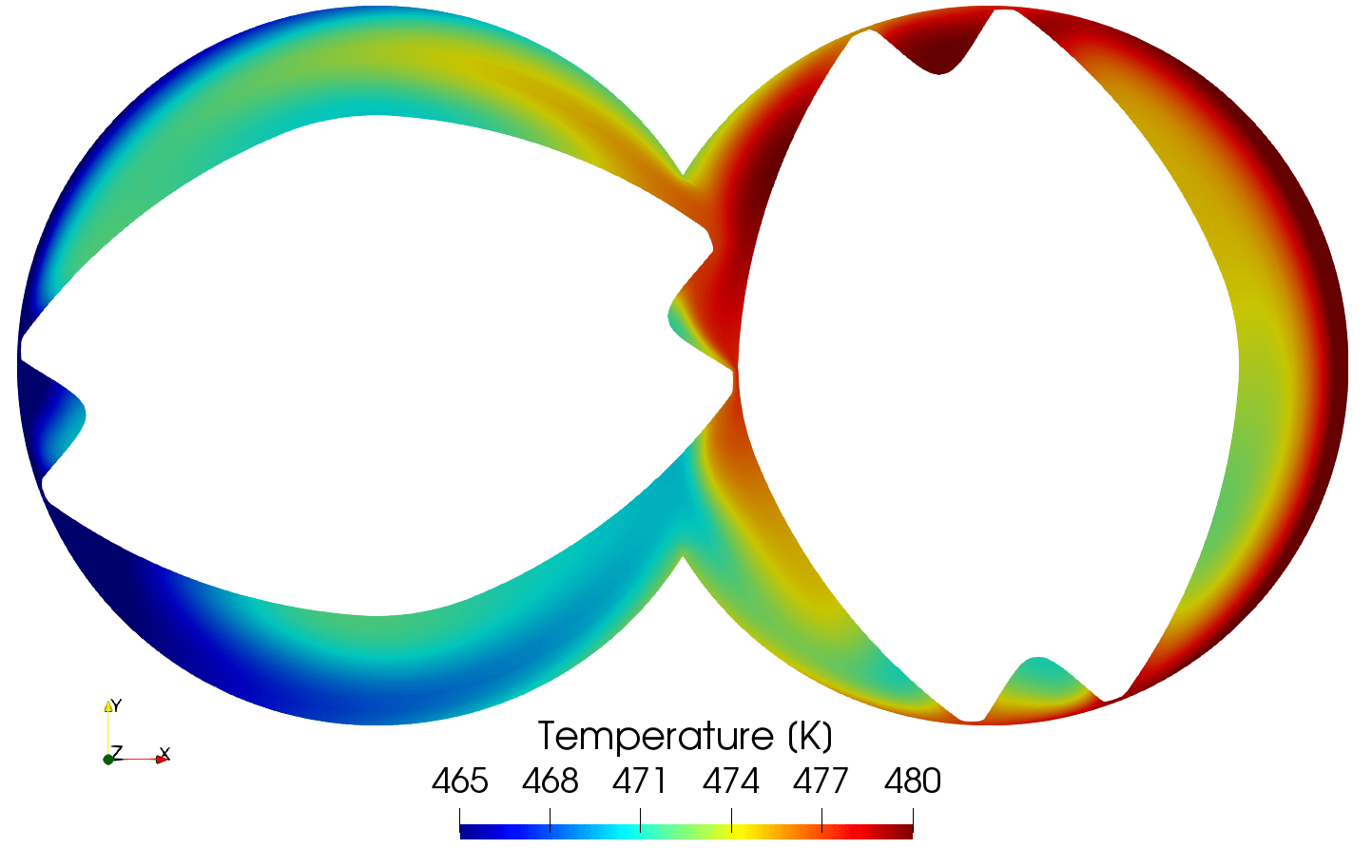}\label{fig:tempMixCut2DNoViscDissd}}
  \centering
  \caption{Temperature field for Config. 1 for selected time steps demonstrating relevant physical effects without viscous dissipation. The results have been computed on mesh 2 using a time step of $\Delta t = 0.003125 s$.}
  \label{fig:tempMixCut2DNoViscDiss}
\end{figure}

We investigate the temperature field for Config. 1 for different time steps on different meshes. In a first step we neglect the effects of viscous dissipation to solely demonstrate how the hot and cold melt is pushed from left to right and vice versa. We use three different meshes to show that the results are mesh independent. Therefore, we compare the temperature results along a line between the lower and upper cusp point at time $t = 1.625s$. The time step size is $\Delta t = 0.00625 s$. The results are given in Fig. \ref{fig:refinmentSpaceNoViscDiss}. Mesh convergence is clearly visible. The results on the coarsest mesh already show very good agreement to the ones obtained on the finest mesh.
In the following, we investigate the choice of an appropriate time step size. We compute the temperature results on mesh 2 for four different time steps, \ref{fig:refinmentTimeNoViscDiss}. We can observe a difference between results for the different time step sizes especially inside the small gap region.
However, even for $\Delta t = 0.0125 s$ we obtain results with a reasonable error for industrial applications. The difference between results obtained with $\Delta t = 0.0625 s$ and the finest time step size $\Delta t = 0.0015625 s$ are less than $1 \; \%$.

Figure \ref{fig:tempMixCut2DNoViscDiss} shows the temperature field at different angles computed on mesh 2 with $\Delta t = 0.003125 s$.
As already described, cold melt is transported by the screws from left to right.
Inside the recessed portion of the screw in the lower part of the left barrel, cold melt is transported into the intermeshing area, see Figure \ref{fig:tempMixCut2DNoViscDissa}. At time $t=1.675 s$, the cold fluid has been transported into the intermeshing area and starts to be convected into the lower part of the right barrel. The flow direction inside the small gaps between the screw is negative, pushing warmer melt into this area.
In the following time steps, $t=1.725 s$ and $t=1.775 s$, more and more warm melt is pushed through the small gap into the lower part of the intermeshing area, leading to a temperature increase there. The temperature field emerges smoothly in time showing the good quality of the underlying meshes as well as numerical methods. \\

\begin{figure}
  \centering
  \subfigure[Refinement in space with $\Delta t = 0.0625 s$.]{\includegraphics[width=.45\linewidth]{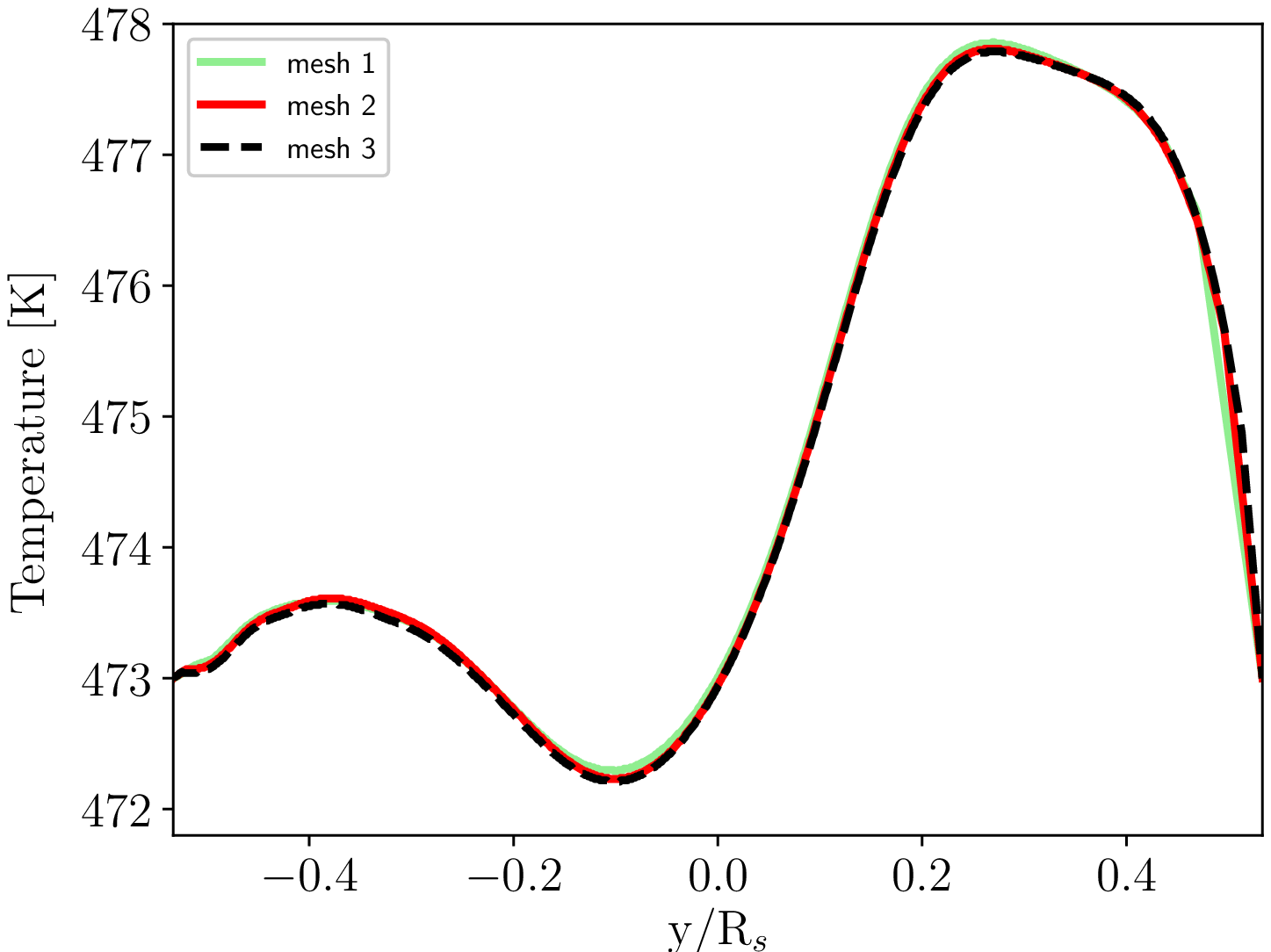}\label{fig:refinmentSpaceViscDiss}}
  \centering
  \subfigure[Refinement in time on mesh 2.]{\includegraphics[width=.45\linewidth]{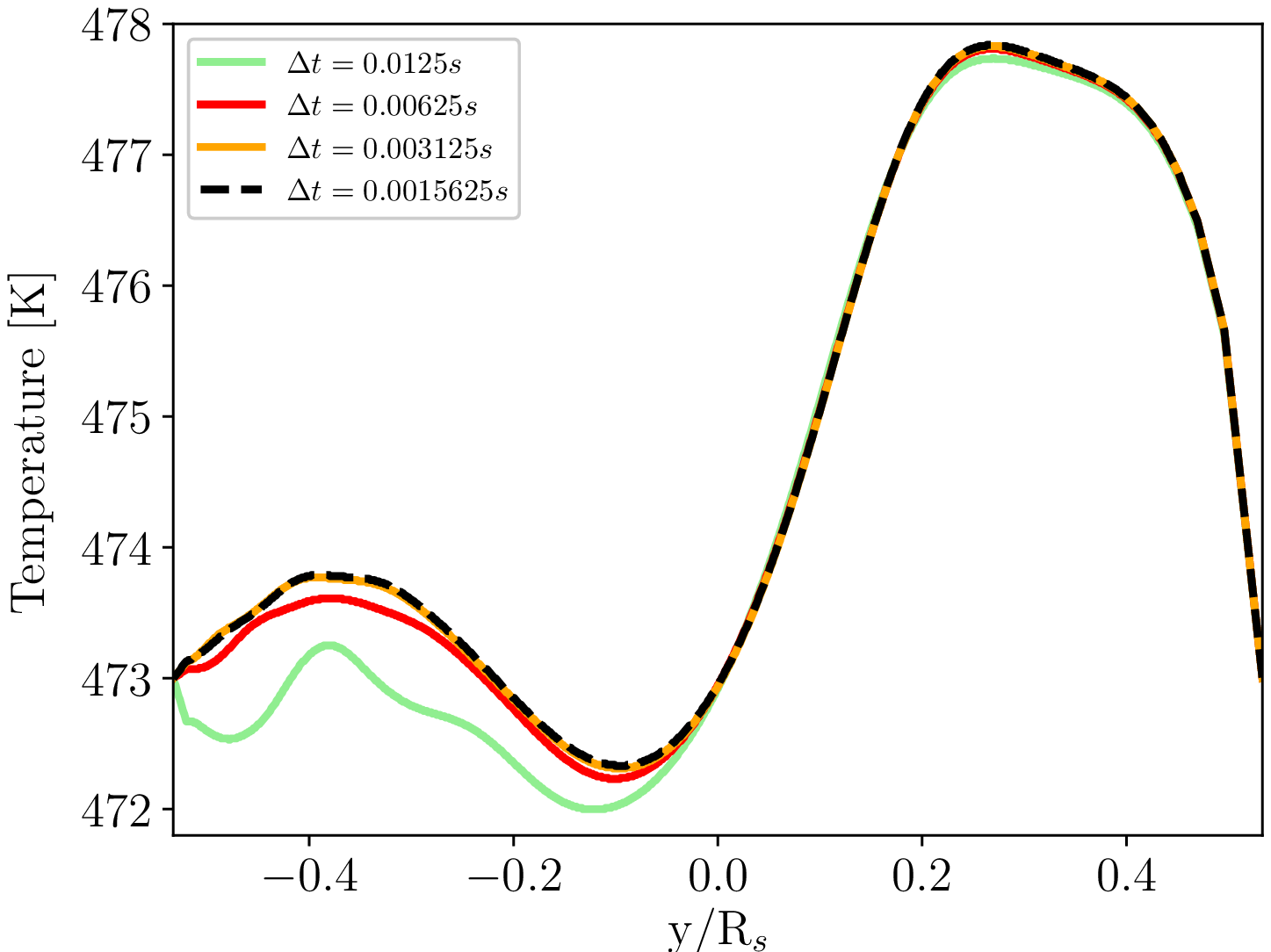}\label{fig:refinmentTimeViscDiss}}
  \centering
  \caption{Mesh and time step refinement study for Config. 1 including viscous dissipation. The temperature results are compared along a line from the lower to upper cusp point at time $t = 1.625 s$.}
  \label{fig:tempMixCut2DViscDissRef}
\end{figure}

\begin{figure}
  \centering
  \subfigure[$t=1.625 s$]{\includegraphics[width=.38\linewidth]{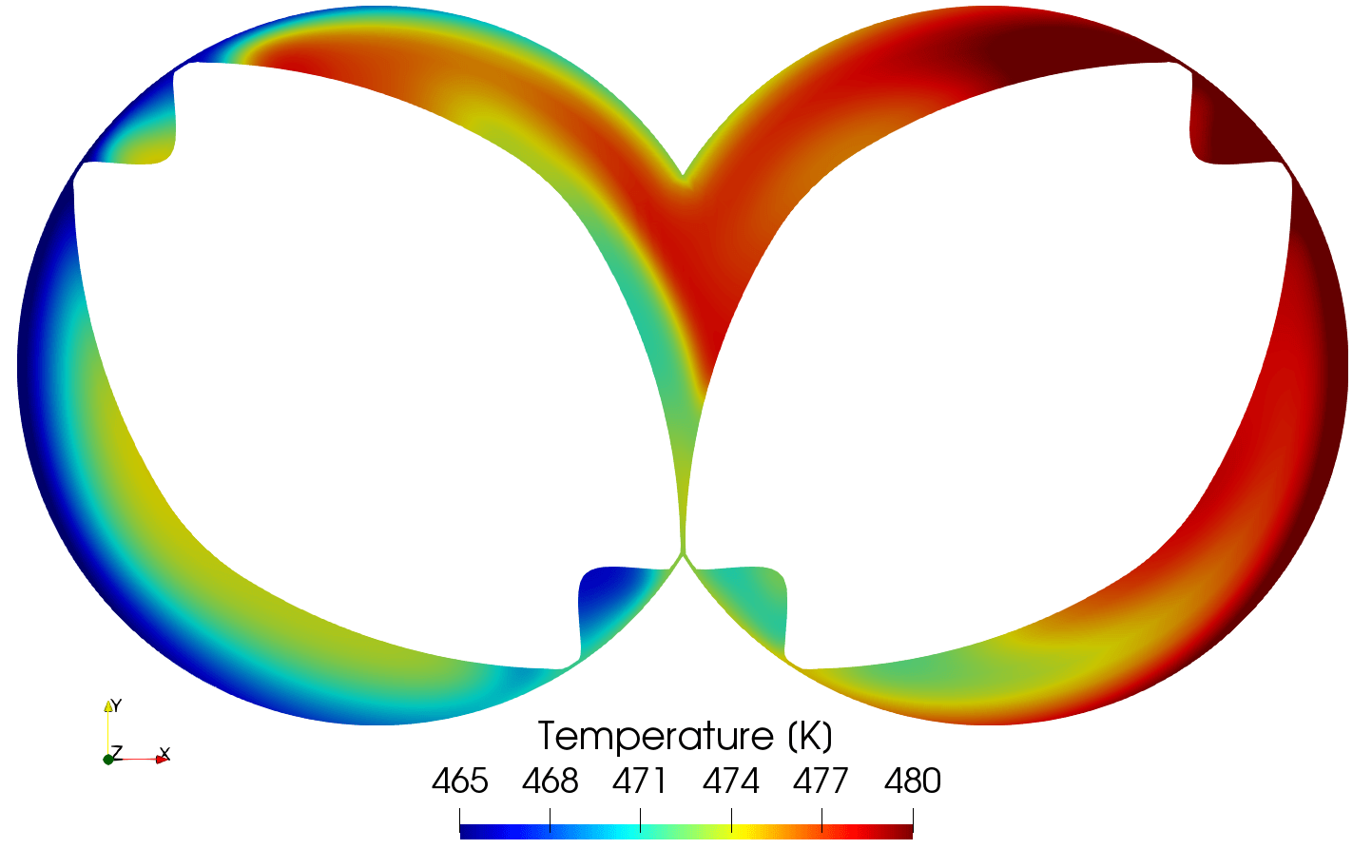}\label{fig:tempMixCut2Da}}
  \centering
  \subfigure[$t=1.675 s$]{\includegraphics[width=.38\linewidth]{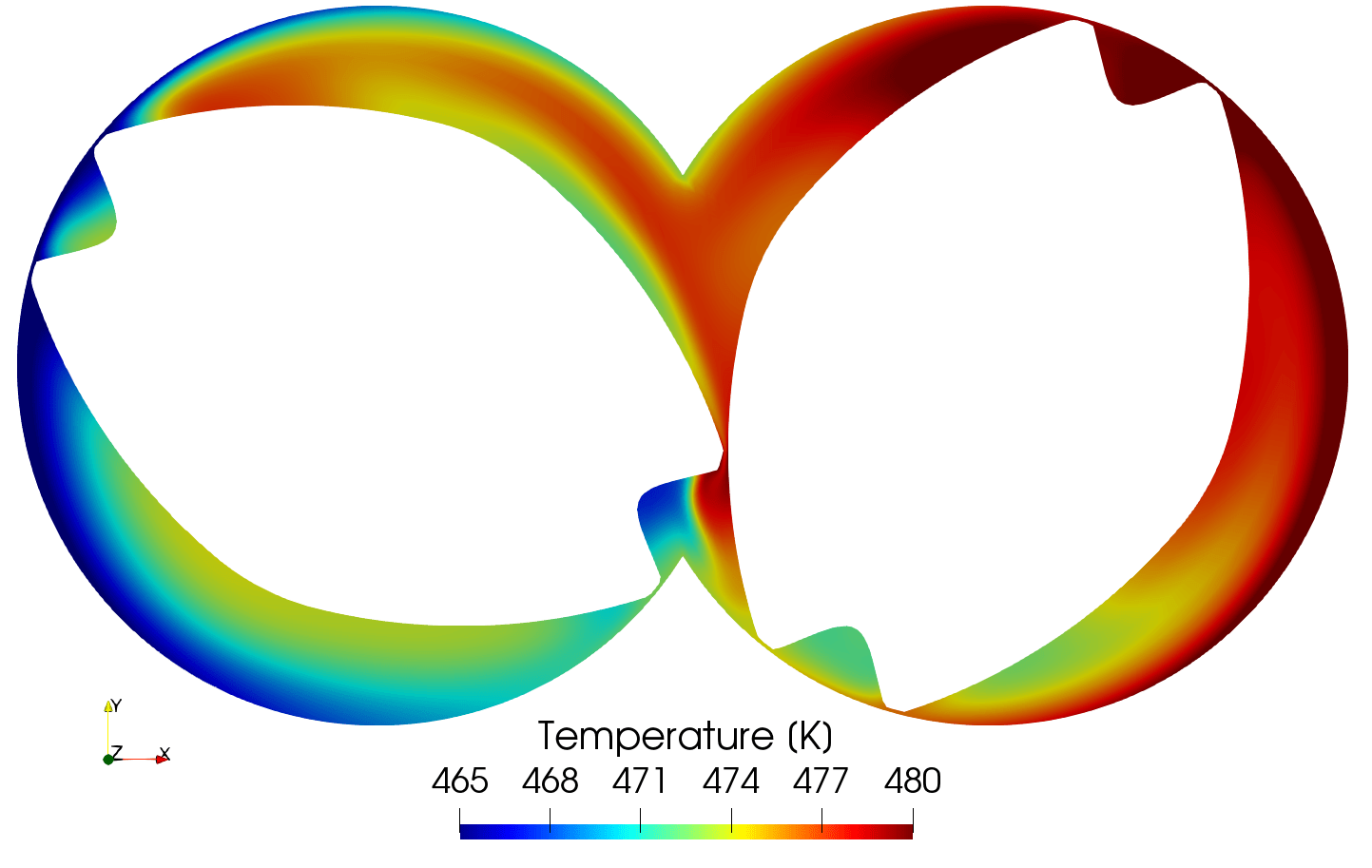}\label{fig:tempMixCut2Db}}
  \centering
  \subfigure[$t=1.725 s$]{\includegraphics[width=.38\linewidth]{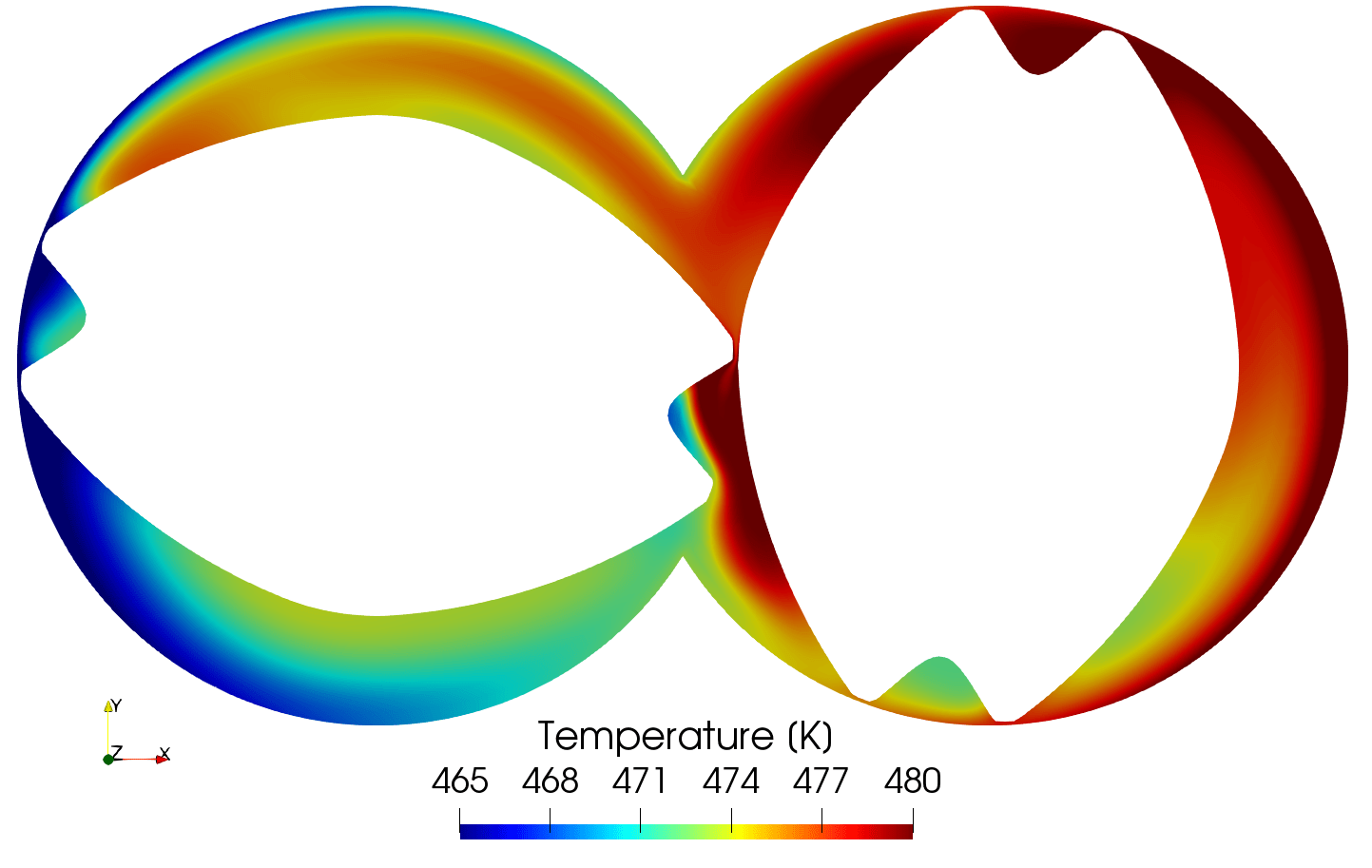}\label{fig:tempMixCut2Dc}}
  \centering
  \subfigure[$t=1.775 s$]{\includegraphics[width=.38\linewidth]{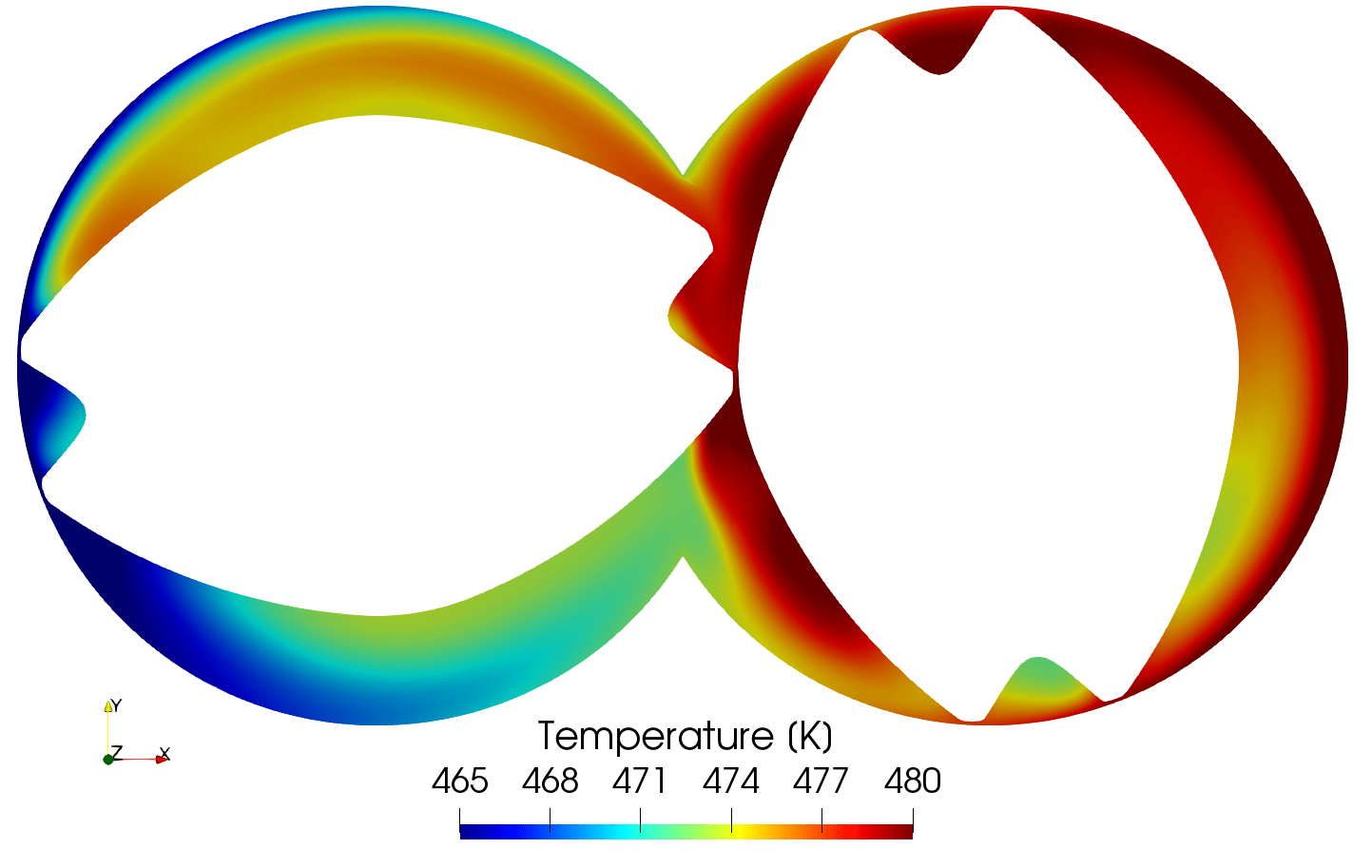}\label{fig:tempMixCut2Dd}}
  \centering
  \caption{Temperature field for Config. 1 for selected time steps demonstrating relevant physical effects including viscous dissipation. The results have been computed on mesh 2 using a time step of $\Delta t = 0.003125 \; s$.}
  \label{fig:tempMixCut2DViscDiss}
\end{figure}

In the following, we include viscous dissipation. Thus, we expect the melt to heat up due to high shear rates inside the small gaps. Again, we analyze the effect of the mesh and time step size. The results are shown in Fig. \ref{fig:tempMixCut2DViscDissRef}. Similar to the case without viscous dissipation, we see a good convergence in time and space. However, the difference of the solution between the coarsest and finest time step increases especially inside the small gap region. This shows the additional complexity introduced to the simulation by viscous dissipation. The temperature field at different angles computed on mesh 2 with time step $\Delta t = 0.003125 s$ is shown in Figure \ref{fig:tempMixCut2DViscDiss}.
Comparing the results with and without viscous dissipation, we can clearly observe the effects of the heating of the melt due to viscous dissipation. At time $t=1.675 s$, heated melt inside the small gap is pushed into the lower part of the intermeshing area. This continues for the following time steps, leading to a major increase of temperature in that area. The temperatures are more than four degrees above those resulting from neglecting viscous dissipation. However, the convective behavior of the melt is similar to the case with no viscous dissipation. \\

\subsubsection*{Config. 2:}

\begin{figure}
  \centering
  \subfigure[$t=1.5 s$]{\includegraphics[width=.38\linewidth]{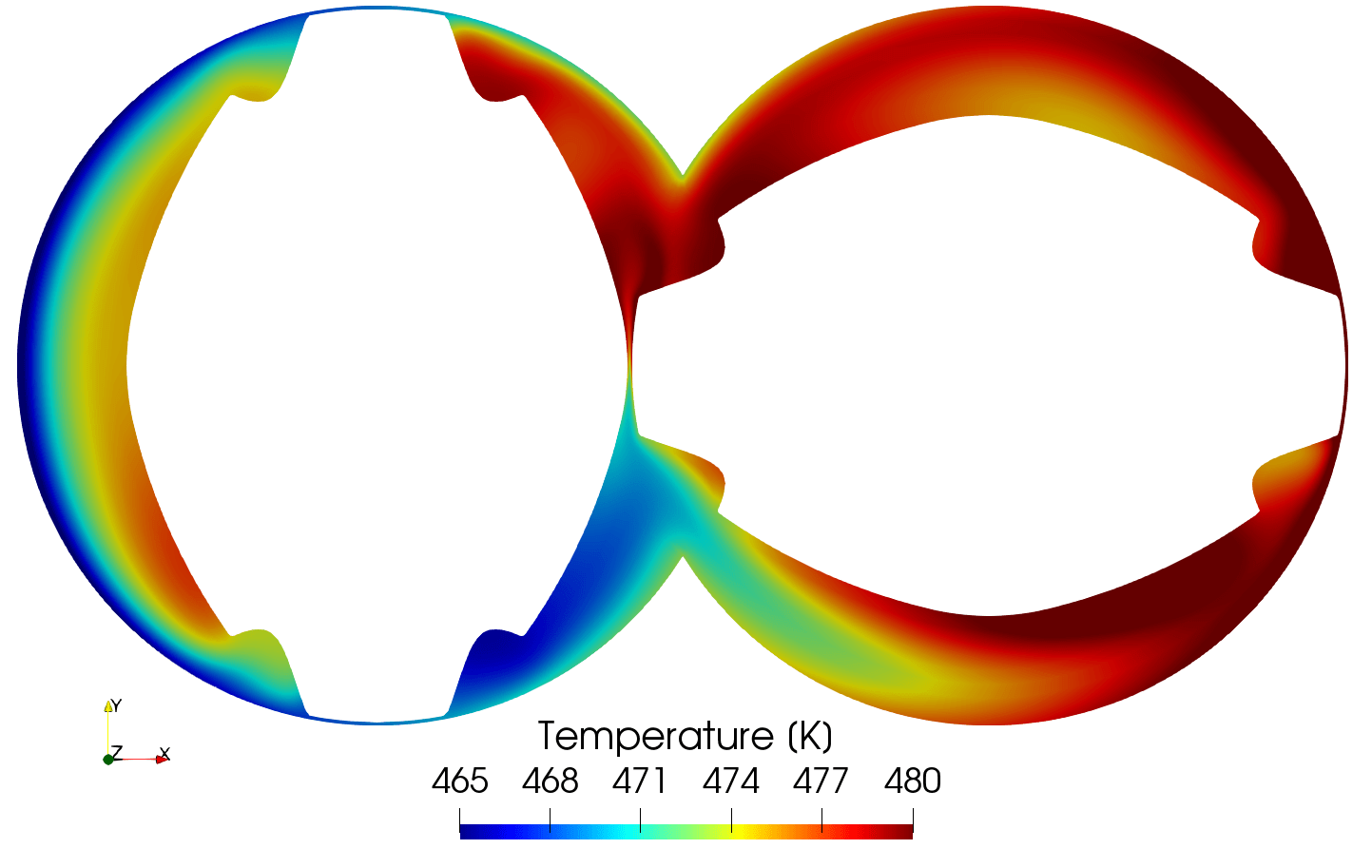}\label{fig:tempMix2Da}}
  \centering
  \subfigure[$t=1.55 s$]{\includegraphics[width=.38\linewidth]{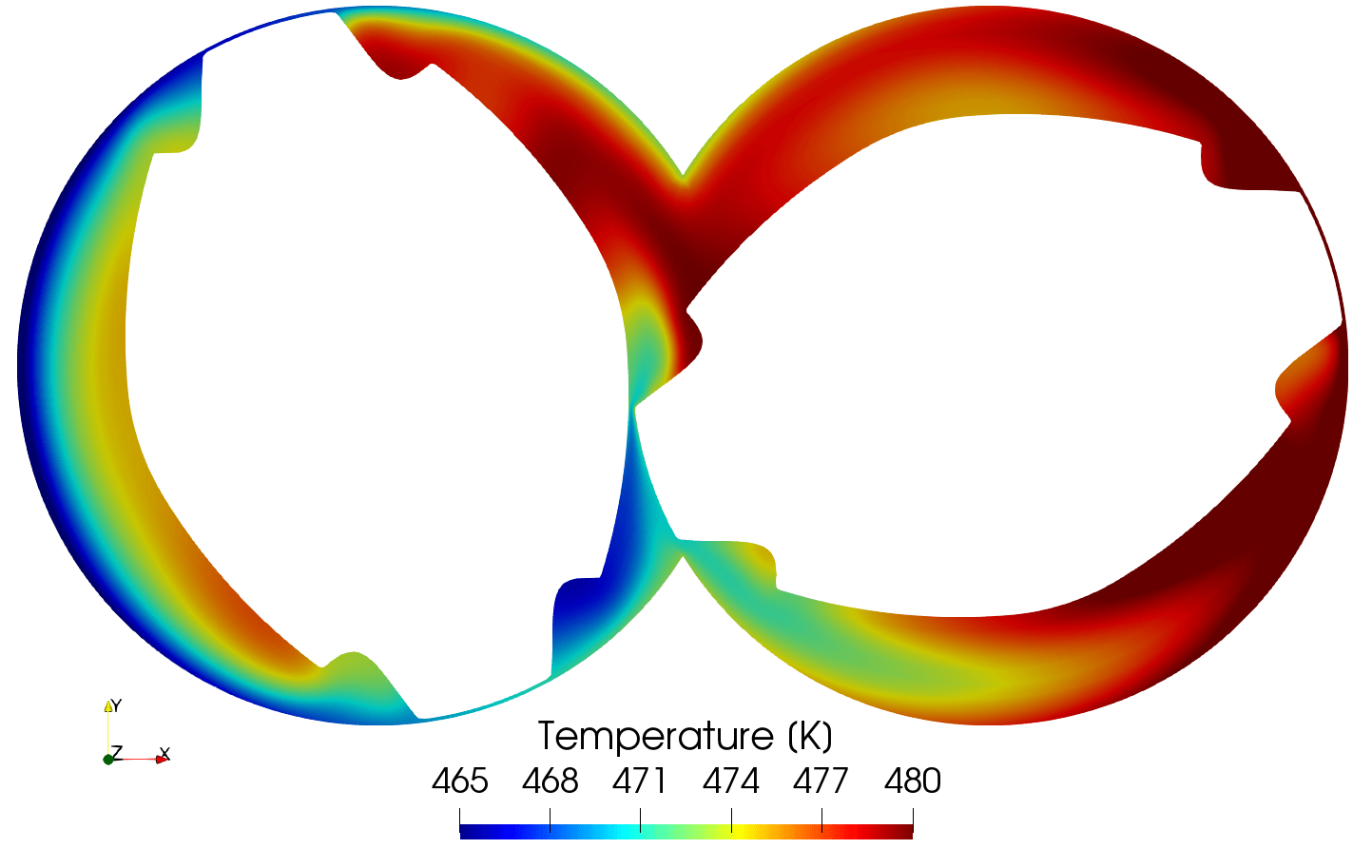}\label{fig:tempMix2Db}}
  \centering
  \subfigure[$t=1.6 s$]{\includegraphics[width=.38\linewidth]{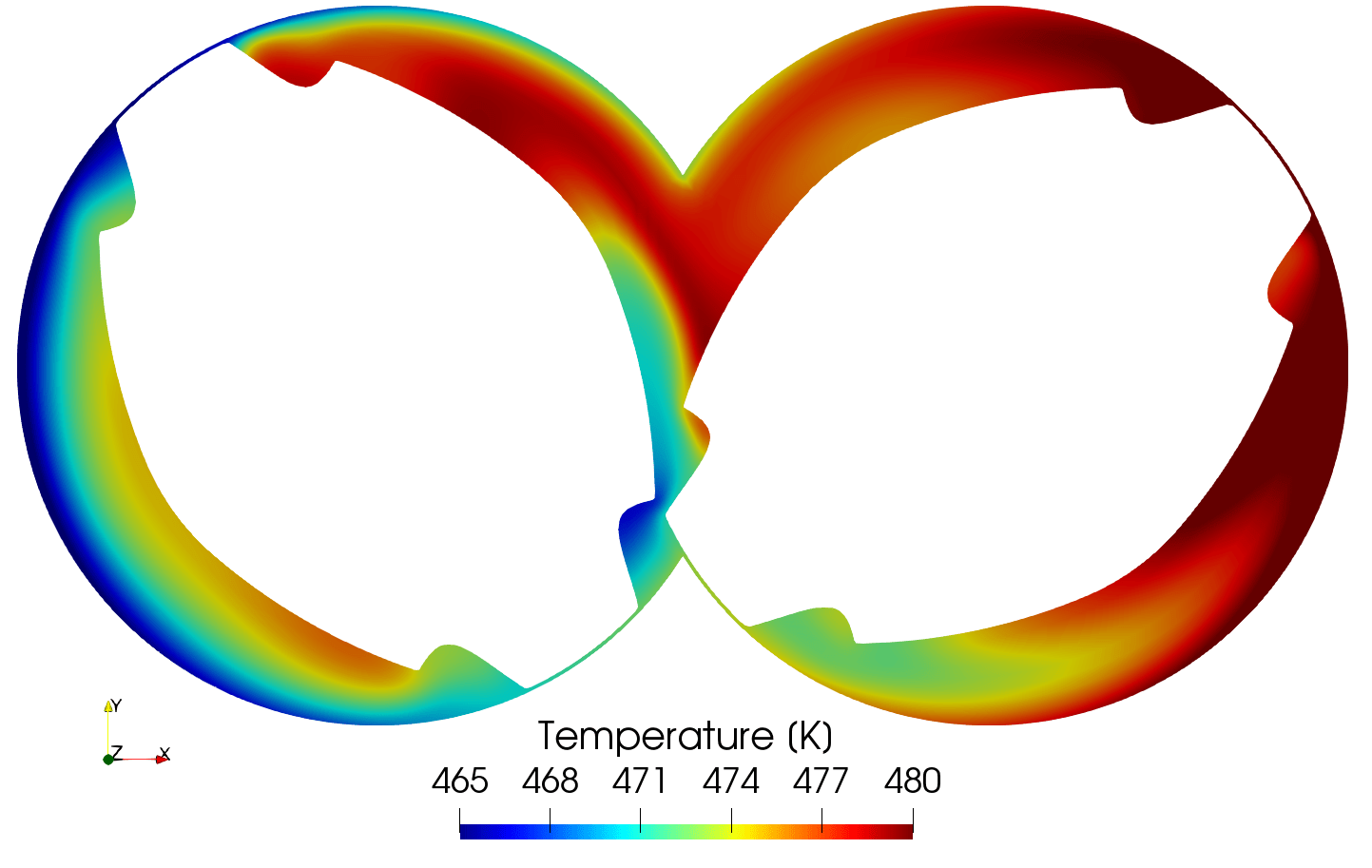}\label{fig:tempMix2Dc}}
  \centering
  \subfigure[$t=1.675 s$]{\includegraphics[width=.38\linewidth]{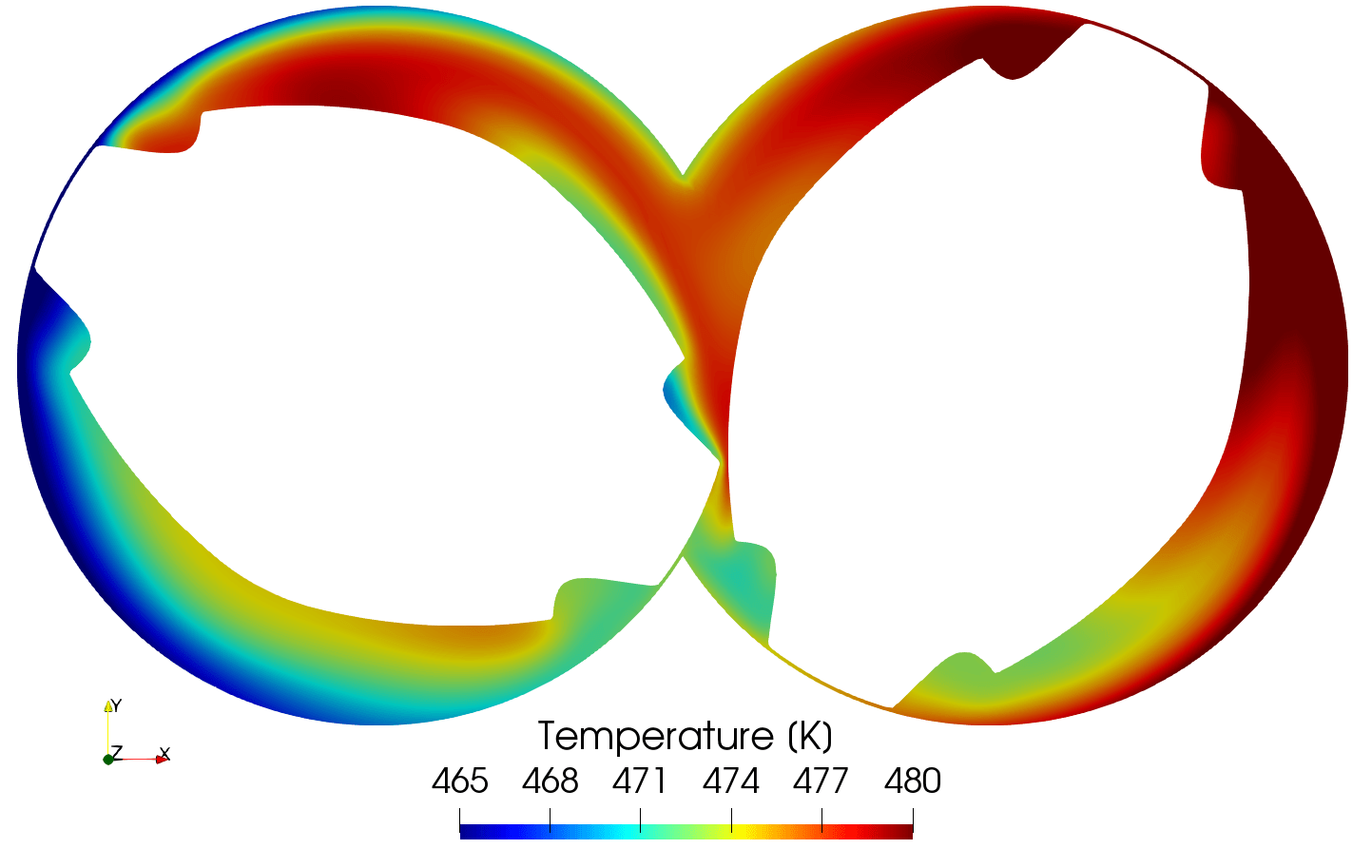}\label{fig:tempMix2Dd}}
  \centering
  \caption{Temperature field for Config. 2 for selected time steps demonstrating relevant physical effects including viscous dissipation.}
  \label{fig:tempMix2D}
\end{figure}

For Config. 2, we analyze the temperature and flow behavior at an earlier stage, see Figure \ref{fig:tempMix2D}. Based on the numerical experiments for Config. 1, we use a time step of $\Delta t = 0.003125 s$.
In this configuration, the cold melt is not transported inside the recessed portion, but is pushed into the intermeshing area due to a larger gap region between screw and barrel, see Figure \ref{fig:tempMix2Da}. Between time $t=1.55 \; s$ and $t=1.6 s$, the melt is pushed upwards out of the cold temperature part, which decreases the temperature in the upper part between the screws. At time $t=1.675 s$ we can observe an increase of temperature in the lower part of the intermeshing region. The direction of flow inside the screw-screw gap has changed sign and the melt is pushed downwards through the gap. Again, the melt is heated up especially inside this small gap due to high shear rates. All these observations are in line with what one would expect by purely looking at flow results. \\

The results for the two configurations show the potential of the presented spline-based meshing approach. We are able to generate high-quality meshes for extremely challenging moving domains.
Only generating a certain amount of slices and interpolating all other meshes seems to be a valid approach. Due to the space-time approach of the method, this unsteady 2D test case is already a proof of concept for 3D. Furthermore, it is noteworthy that the time spent for updating the mesh was less than $0.1 \; \%$. The simulations have been run on 24 cores on the Intel Xeon based RWTH cluster using an MPI parallelization.
For example, computing one revolution for Config. 1 on mesh 2 using a time step of $\Delta t = 0.00625 s$ takes $115 s$. The time spent for the mesh update was only $0.02 s$.

\begin{remark}
The appeal of the presented methodology becomes apparent when performing mesh refinement studies as well as time step size studies. The spline parameterization for the screw only has to be generated once at the beginning. For a mesh refinement study only the \textit{'scaffolding'} needs to be regenerated which is a simple evaluation of the spline parameterization. A time step refinement study is even simpler. Given a \textit{'scaffolding'} for instances of $\Theta$ in the interval $[0,\pi]$, we can simply compute the mesh at any time instance by interpolating between the generated instances. Thus, it is not necessary to generate any new mesh in case one aims to adapt the time step size.
\end{remark}

\subsection{3D Application Case}

Within this section, we aim to show the functioning of the spline-based parameterization technique for real 3D applications. Once again, we consider the temperature-dependent flow of a plastic melt through a complex screw geometry.
The plastic melt is modeled using the Cross-WLF model. The model parameters employed are based on a polypropylene from the product portfolio of a leading raw material manufacturer. The parameters are given in Table \ref{table:CrossWLF3D}. The plastic melt has density $\rho = 710 \; kg / m^3$, specific heat $c_p = 2400 \; J/(kg \;K)$ and thermal conductivity $\kappa_0 = 0.5 \; W/(m\;s)$.

\begin{figure}[h]
  \centering
  \subfigure[2D cross section.]{\includegraphics[width=.45\linewidth]{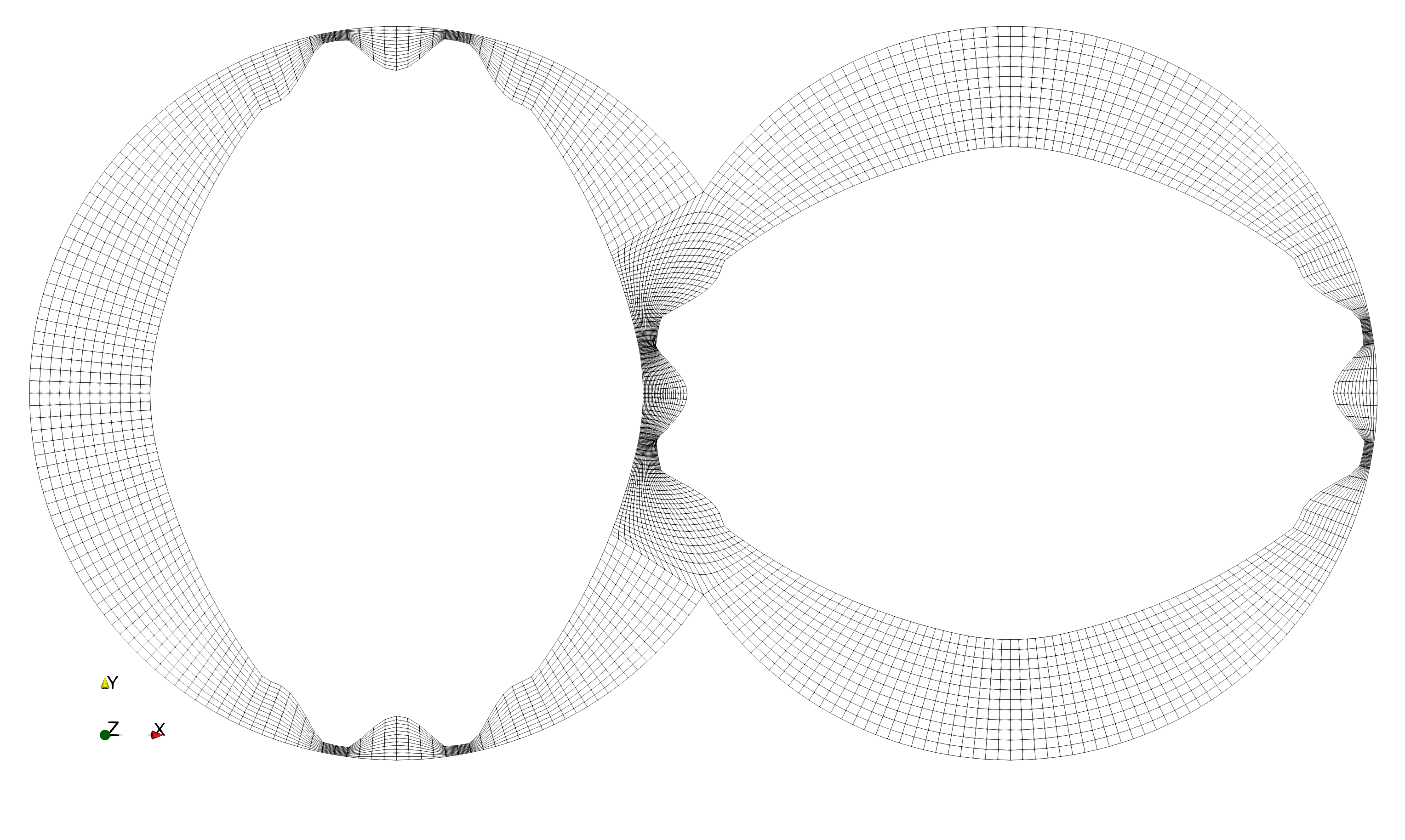}\label{fig:domain2D3D}}
  \centering
  \subfigure[3D sketch of extruder including in/outflow extension.]{\includegraphics[width=.4\linewidth]{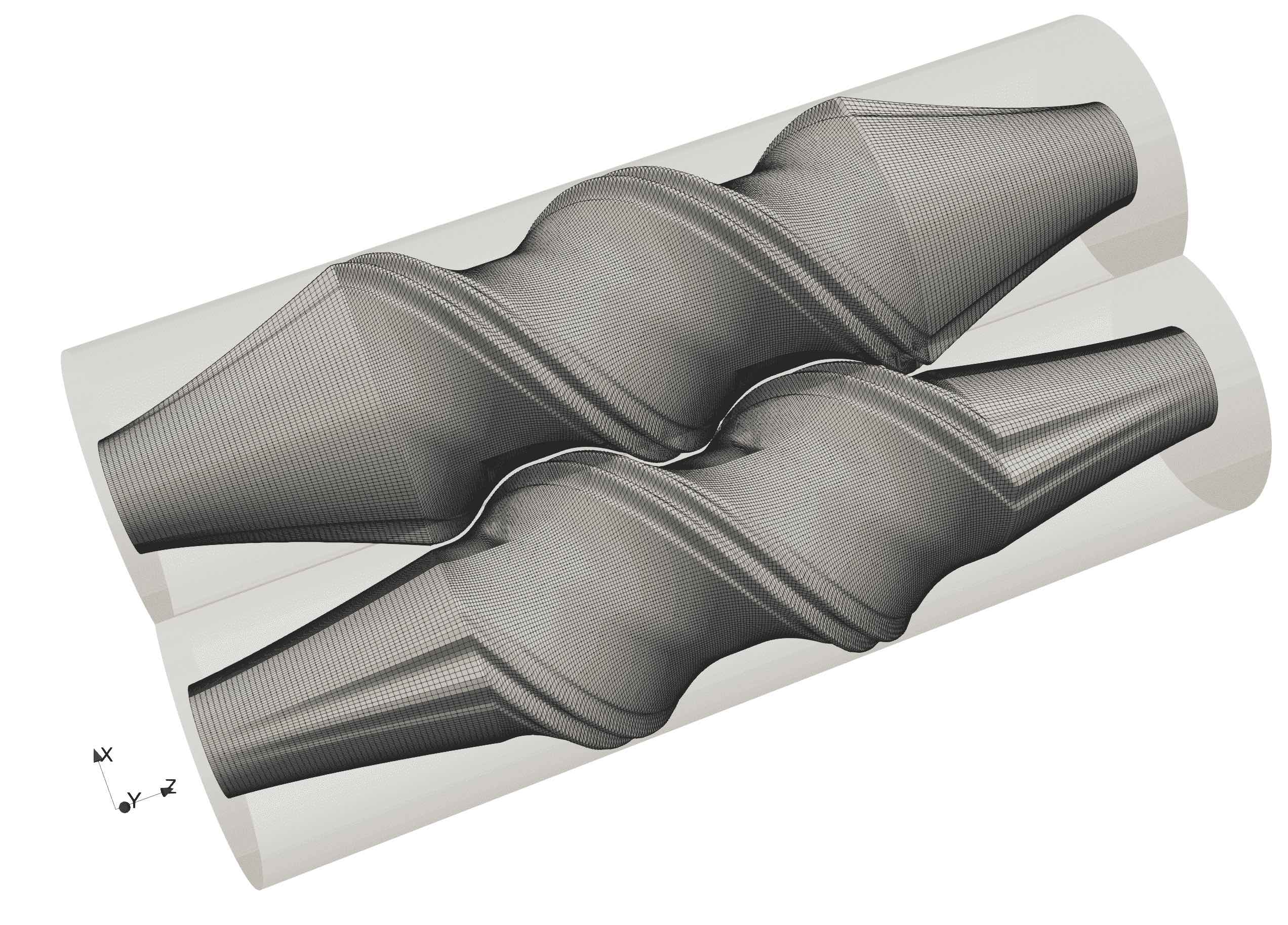}\label{fig:domain3D3D}}
  \caption{3D mixing element.}
  \label{fig:domainMixing3D}
\end{figure}

\begin{table}[h]
  \begin{minipage}[t]{0.48\linewidth}
      \centering
      \begin{tabular}{l r}
       \hline
       Screw radius $R_s$  & 0.156 $m$ \\
       Center line distance $C_l$ & 0.262 $m$ \\
       Screw-screw clearance $\delta _s$ &  0.004 $m$ \\
       Screw-barrel clearance $\delta _b$ & 0.004 $m$ \\
       Pitch length $p_l$ & 0.28 $m$ \\
       \hline
       \end{tabular}
       \caption{Geometry parameters of a 3D mixing screw element for temperature-dependent flow.}
       \label{table:screw3DTempUnsteady}
  \end{minipage}
  \begin{minipage}[t]{0.48\linewidth}
  \centering
  \begin{tabular}{l c c}
  \hline
  $D1$  & 1.2e+14 & $Pa \; s$ \\
  $\tau^{*}$ & 25680.0 & $Pa$ \\
  $n$ & 0.29 & - \\
  $T_{ref}$ & 263.15 & $K$ \\
  $A1$ & 28.32 & - \\
  $A2$ & 51.60 & $K$ \\
  \hline
\end{tabular}
\caption{Cross-WLF parameters.}
\label{table:CrossWLF3D}
\end{minipage}
\begin{minipage}[c]{1.0\linewidth}
\centering
\begin{tabular}{l c c c c c c c c}
 \hline
 $n_s$ C-grid & $n_s$ separator & $n_s$ total & $n_r$ & $n_a$ &$\#$ elements & $n_{\mu}$ & $n_{\nu}$ total & $n_{slices}$ $\pi$ \\
  150 & 70 & 220  & 12 & 300 & 1584000 & 220 & 16 & 101\\
 \hline
 \end{tabular}
 \caption{Mesh discretization for 3D mixing elements.}
 \label{table:mesh3DMixing}
\end{minipage}
\end{table}

The 2D screw geometry cross section is a simplified combination of the two 2D configurations of the previous section, see Figure \ref{fig:domain2D3D}. The screw geometry parameters are given in Table \ref{table:screw3DTempUnsteady}. The 3D setup is shown in Figure \ref{fig:domain3D3D}. For simplicity we only consider a single screw element. The screws rotate in mathematically positive direction with $\omega = 120$ rpm. The computational domain is extended at the inflow and outflow. By that, a solely positive velocity in $z$-direction is guaranteed at the outflow.
This is important since we set a natural boundary condition and thus, negative velocities might cause instabilities. Furthermore, this setting is also similar to a industry-like twin-screw extruder where the screw section leads to the die.
The extended inflow circumvents high shear rates in the inflow area, which would result in an unnaturally high temperature increase.
The extension is achieved by relaxing the 2D screw surface over a distance of $0.28 \; m$ to a circle with a radius of $0.06 \; m$. For the mesh, we simply interpolate between the original 2D mesh and a structured mesh between the two circles.
We use a background mesh with $n_s = 220$ elements on the screw, split into 70 elements inside the separator and 150 for the C-grid. $n_r = 12$ elements are used in radial direction.
A full pitch length is discretized using 200 elements in axial direction and the inflow and outflow extensions are discretized with 50 elements each. The
 \textit{'scaffolding'} generated by the spline parameterization consists of all points in circumferential direction, $n_s = n_{\nu}$ and $n_{\mu}=16$ inside the separator. In order to account for the rotation, the \textit{'scaffolding'} is evaluated for 101 equally distributed instances of $\Theta$ in the interval $[0,\pi]$. Storing only the \textit{'scaffolding'} instead of the full mesh for each instance results in  memory savings for the mesh of $75 \; \%$.
All mesh parameters are listed in Table \ref{table:mesh3DMixing}. \\

\begin{figure}[h]
  \centering
  \subfigure[Pressure results for the entire domain.]{\includegraphics[width=.45\linewidth]{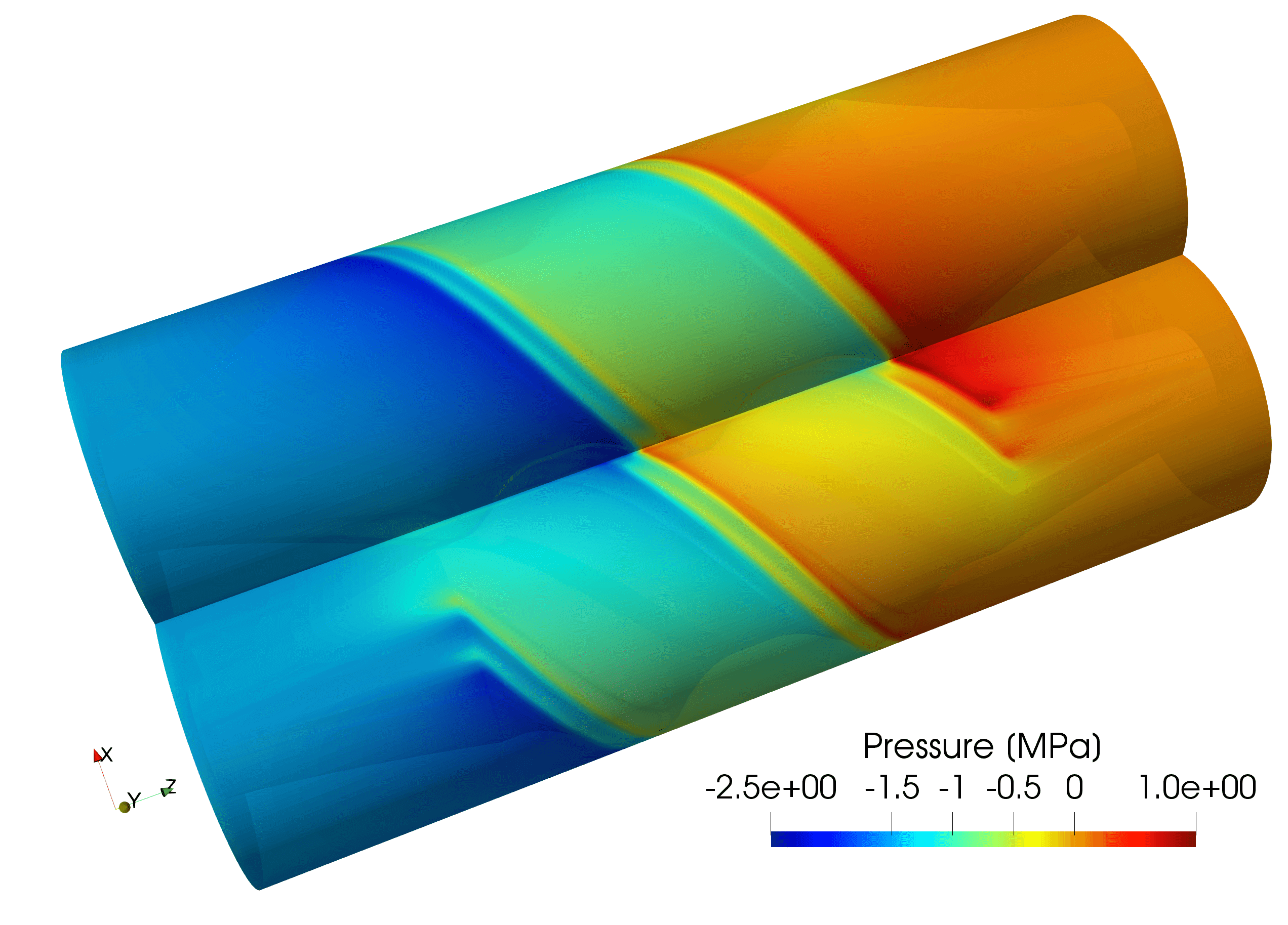}\label{fig:zvel3D}}
  \centering
  \subfigure[Contour for the velocity in $z$-direction at $y=0.074 \; m$.]{\includegraphics[width=.48\linewidth]{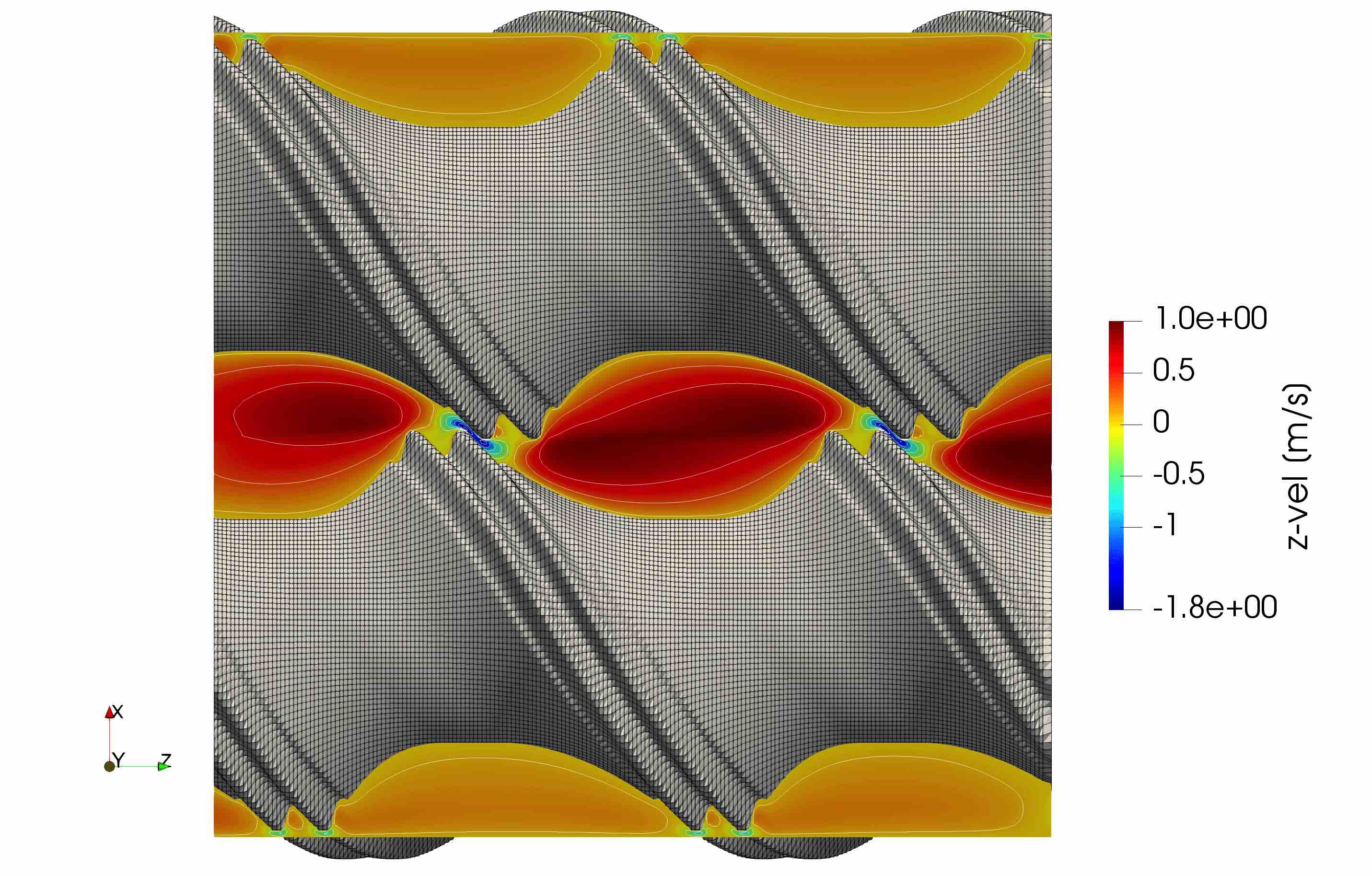}\label{fig:pressure3D}}
  \caption{Flow results for $t = 5.0 s$.}
  \label{fig:flowresultsMixing3D}
\end{figure}

We simulate a scenario with a mass flow rate of $\dot{m} = 25560 \; kg/h$. This is achieved by setting a uniform inflow velocity at the inlet.
The barrel is heated with $T_{barrel} = 473 K  + 10K \; * \; z/1.12m $ and the screws are considered to be adiabatic.
The inflow temperature is $T_{inflow} = 473 \; K$. In order to obtain a good initial condition for the temperature, we compute a steady solution where we neglect the viscous dissipation term and increase the conductivity by a factor of 10. The time step size is $\Delta t = 0.005 s$ which is chosen based on the results of the previous section. The results have been computed on 360 cores on the Intel Xeon based Juelich cluster. Computing one revolution takes $5832 s$, whereas the part spent for the mesh update is only $0.77 s$ which shows the efficiency of the presented meshing approach.

\begin{figure}[h]
  \centering
  \subfigure[$t=5.00 s$]{\includegraphics[width=.3\linewidth]{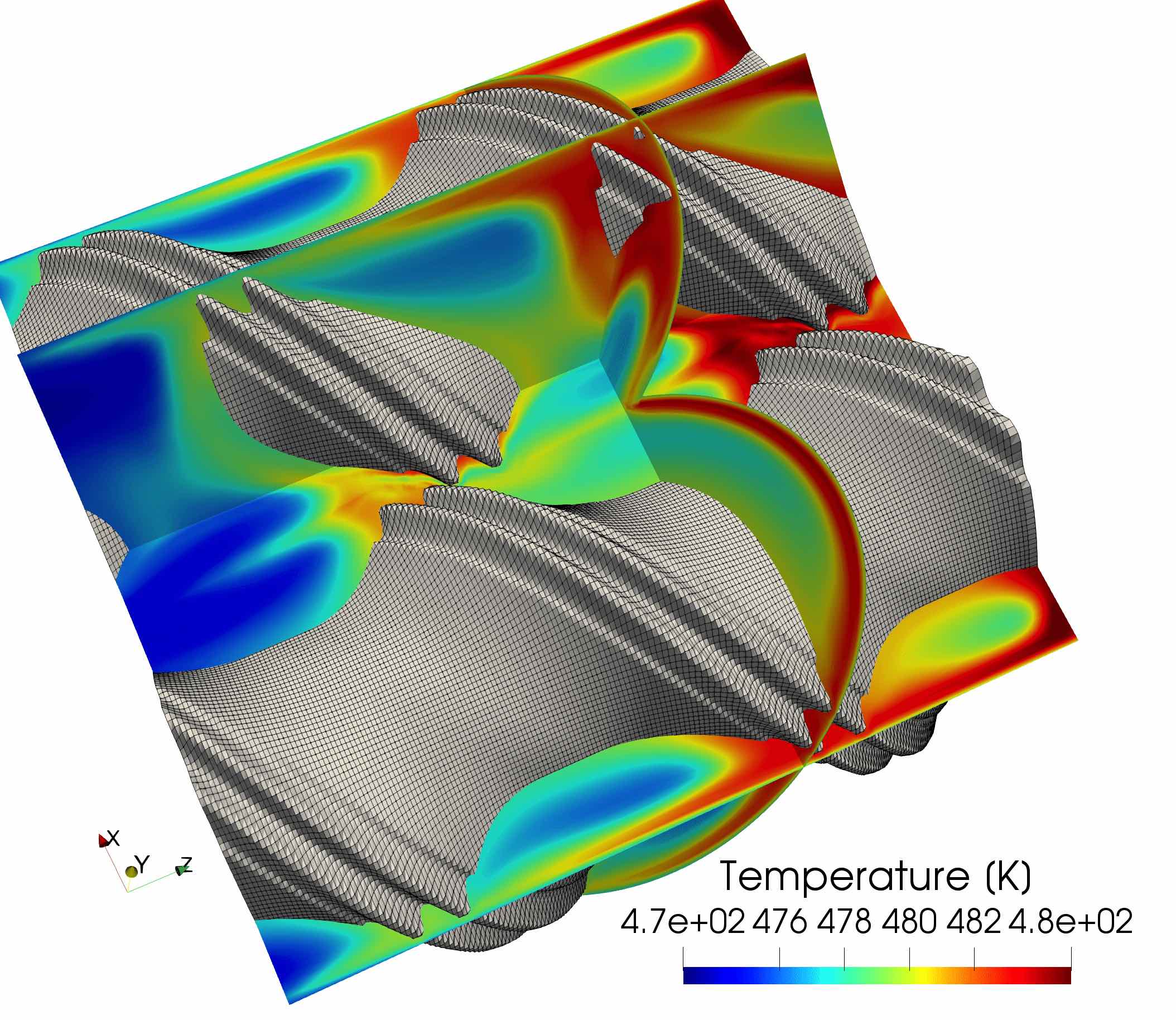}}
  \centering
  \subfigure[$t=5.04 s$]{\includegraphics[width=.3\linewidth]{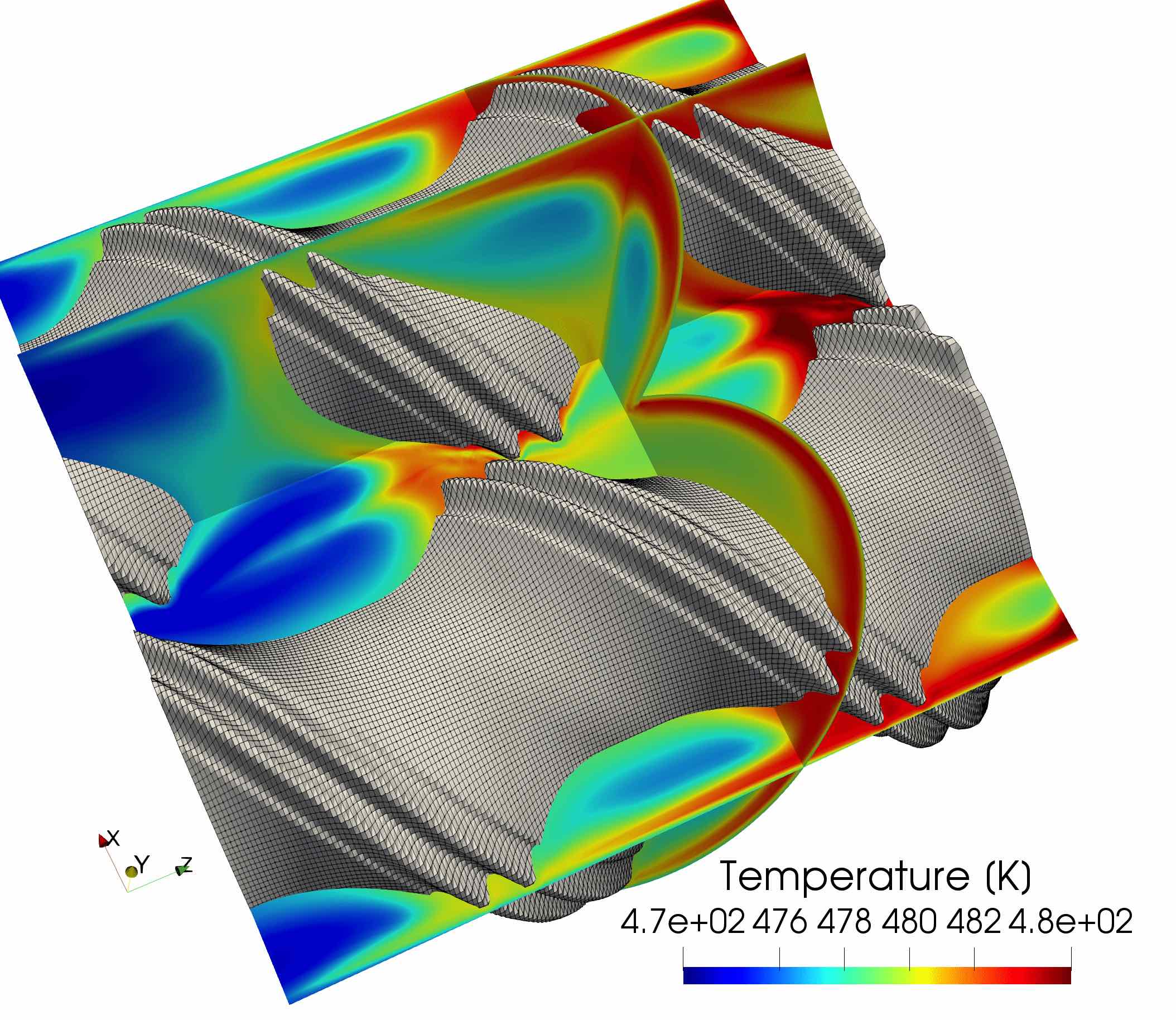}}
  \centering
  \subfigure[$t=5.08 s$]{\includegraphics[width=.3\linewidth]{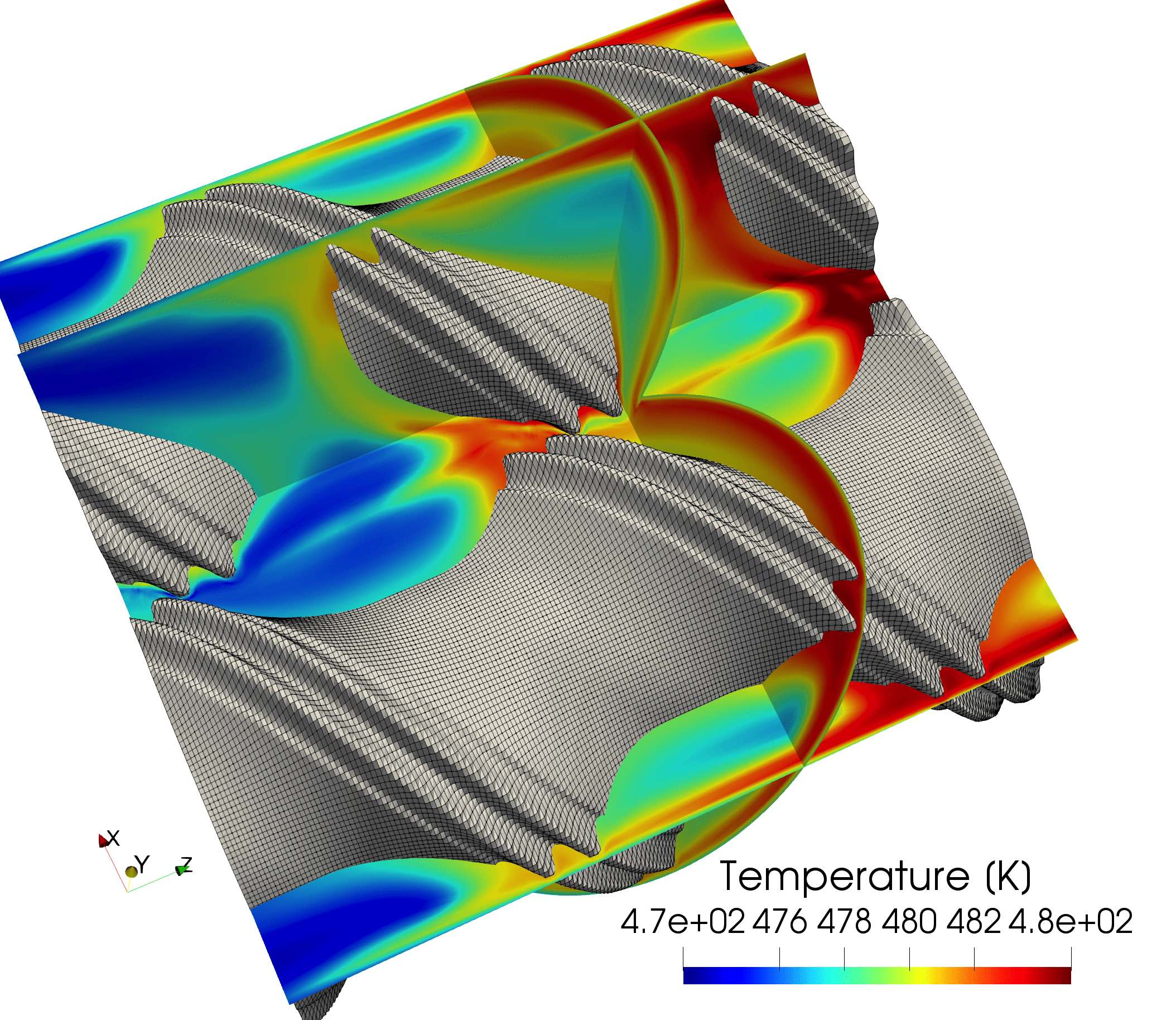}}
  \centering
  \subfigure[$t=5.12 s$]{\includegraphics[width=.3\linewidth]{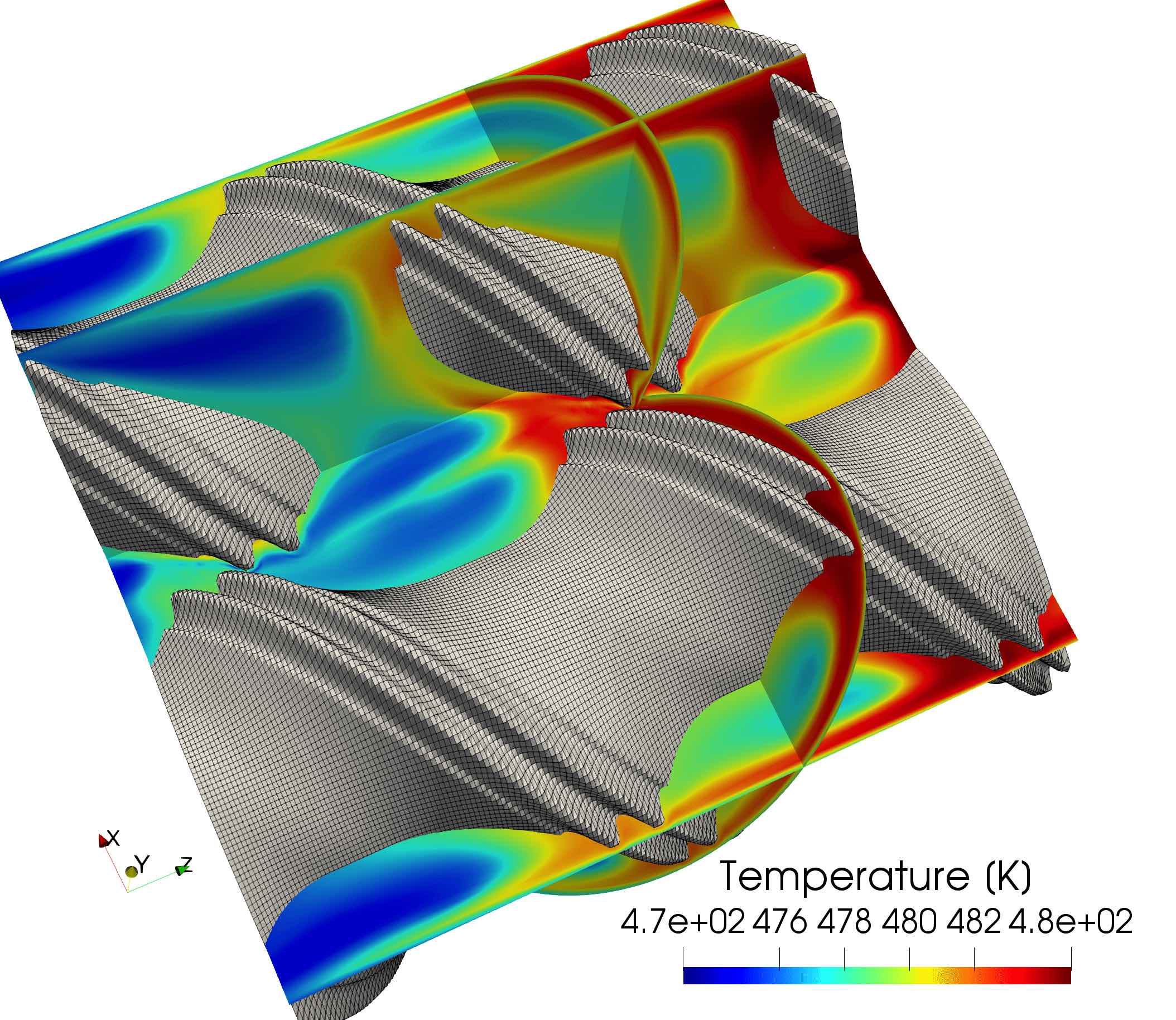}}
  \centering
  \subfigure[$t=5.16 s$]{\includegraphics[width=.3\linewidth]{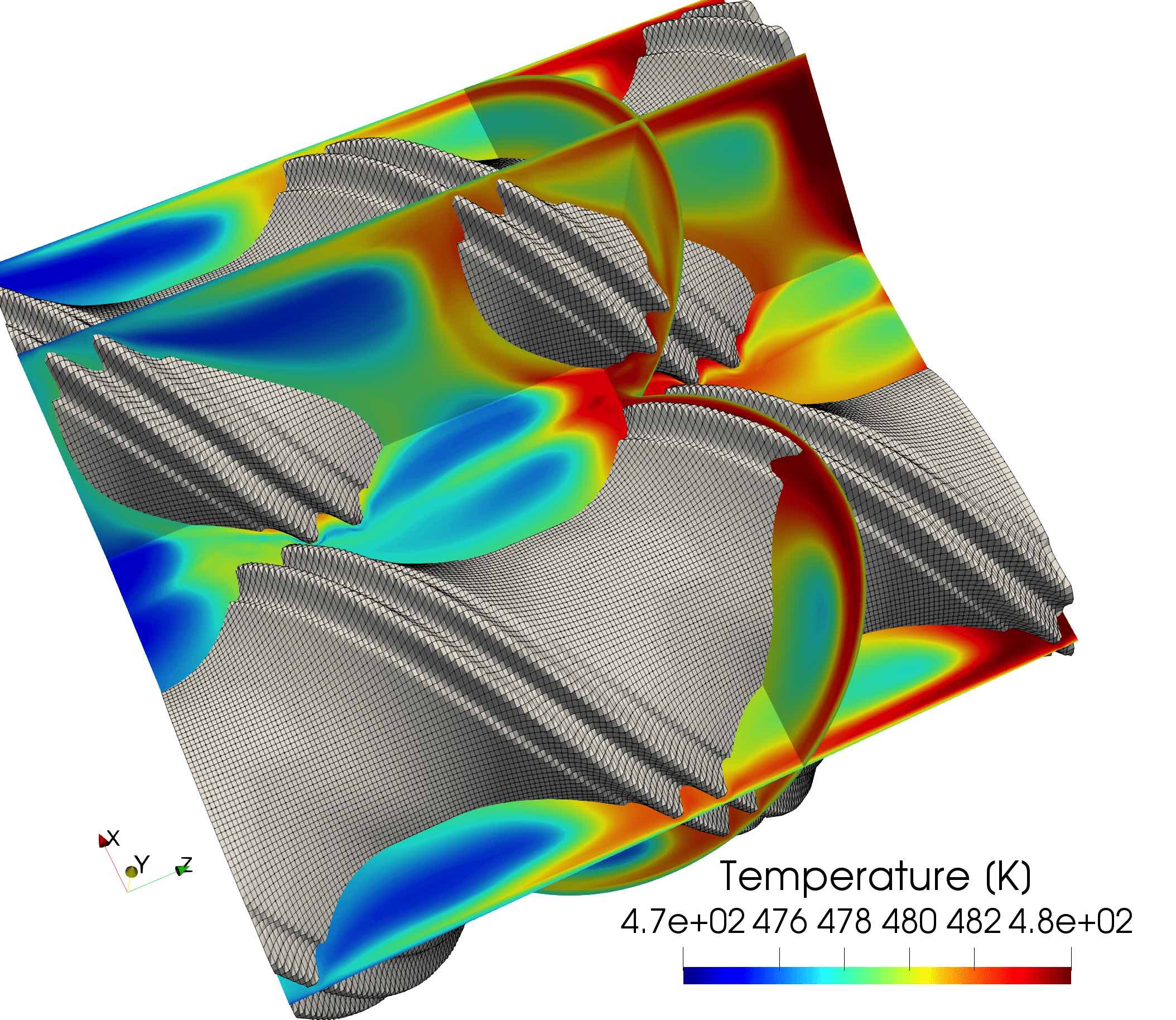}}
  \centering
  \subfigure[$t=5.20 s$]{\includegraphics[width=.3\linewidth]{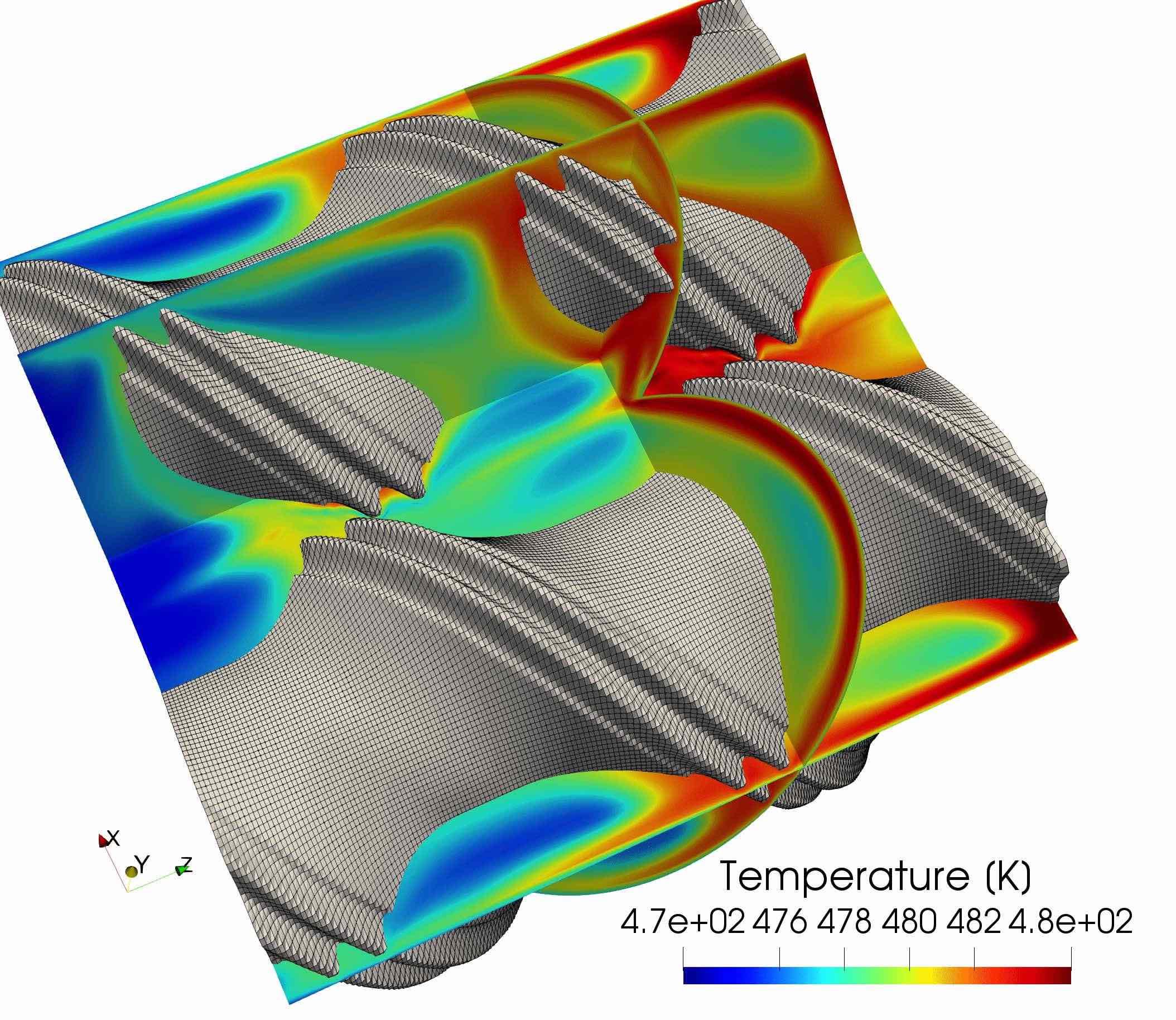}}
  \caption{Temperature distribution in mixing element displayed on three planes: xy-plane at $y=0.074\; m$; a plane with normal $\bold{n}$ = $(1.0,-1.5,0.0)$ and point $\bold{r}$ on the plane with $\bold{r}= (-0.07,0.0,0.0)$; a plane with $\bold{n} = (0.0,-1.0,1.0)$ and $\bold{r} = (0,0,0.56)$. }
  \label{fig:tempresultsMixing3D}
\end{figure}

The flow results are presented in Figure \ref{fig:flowresultsMixing3D}. The screw is able to transport a higher mass flow than the one set at the inflow. This results in an increase of pressure over the screw element, see Figure \ref{fig:pressure3D}. The velocity in $z$-direction is given in a plane at $y=0.074\; m$ in Figure \ref{fig:zvel3D}. A strong backflow behavior can be observed in the small gaps between the screws. This is desirable for twin-screw extruders since it increases their mixing capabilities. Inside the extruder we observe Reynolds numbers up to $Re = 0.1$, which justifies using Navier-Stokes equations.
The temperature field reaches a periodic state after roughly 10 revolutions. The temperature results for half a rotation are shown in Figure \ref{fig:tempresultsMixing3D}.
The melt is heated inside the high shear regions in the small gaps between the screws. Furthermore, we can observe how cold melt is pushed forward inside the high velocity regions, especially inside the wide parts of the intermeshing region.

\section{Conclusion}

Within this work, we presented an efficient and robust meshing strategy that allows to use boundary-conforming finite element methods to compute the unsteady flow of plastic melt inside co-rotating twin-screw extruders.
It is suitable for arbitrarily-shaped screw geometries. The method is a combination of spline-based meshing techniques based on EGG and SRMUM. 2D spline-based geometry descriptions are generated for a certain number of screw orientations. The spline description is evaluated at discrete points to obtain a point cloud that is used as a \textit{'scaffolding'} to adapt a structured background mesh to the actual screw configuration. In contrast to algebraic grid generation, the advantage of the proposed approach is two-fold: on the one hand, it allows for fine tuning of the mesh properties by employing a control mapping (see Section \ref{subsect:Applications_Geometry}). On the other hand, refinement studies are easily accomplished by evaluating the composite spline-mapping in a increasing number of points. Additionally, storing only a limited number of points instead of a full 3D mesh rotation saves a lot of memory. The actual mesh update at run-time is very cheap and requires only a fraction of the time spent for the actual solve of the flow solution. \\
A finite element method based on the Deformable-Spatial-Domain/Stabilized Space-Time (DSD/SST) finite element formulation was used. A 2D test case simulating the flow of an isothermal polymer in a twin-screw extruder cross section served as a validation case for the presented framework including the new meshing strategy.
Convergence of the solution has been demonstrated as well as accordance to results from literature.
Furthermore, two complex, mixing-element-like, screw shapes have been used to show the robustness of the meshing technique in 2D. We computed unsteady temperature-dependent flow of the plastic melt inside the aforementioned screw shapes.
The use of a space-time finite element method already proofs, that the presented method is capable of computing unsteady flow results using boundary-conforming meshes in 3D.
We additionally computed the flow of a temperature-dependent plastic melt in 3D for a complex screw shape in a single element twin-screw extruder to further demonstrate the great potential of the presented approach. \\
In the future, we aim to compute flow results for a larger variety of screw designs. This will be used in order to compare the quantities like residence time distributions or mixing behavior. The advantages of the method already shown for unsteady temperature results can also be exploited by computing solutions of advection-diffusion equations, in order to characterize mixing behavior of the individual screw design.

\section*{Acknowledgements}

The authors gratefully acknowledge the research funding which was partly provided by the MOTOR project that has received funding from the European Unions Horizon 2020 research and innovation programm under grant agreement No 678727. The computations were conducted on computing clusters supplied by the Juelich Aachen Research Alliance (JARA) and the RWTH IT Center.

\bibliographystyle{elsarticle-num}
\bibliography{bibliography}


\end{document}